\documentclass[12pt]{amsart}
\usepackage{amssymb}
\usepackage{mathrsfs}
\usepackage{bm}
\usepackage{graphicx}
\usepackage{color}
\usepackage[all]{xy}

\allowdisplaybreaks[4]

\newtheorem{thm}{Theorem}[section]
\newtheorem{lem}[thm]{Lemma}
\newtheorem{cor}[thm]{Corollary}
\newtheorem{prop}[thm]{Proposition}
\newtheorem{definition}[thm]{Definition}

\def \para{\refstepcounter{thm} \par\medskip\noindent
                \textbf{\thethm .} }

\def \remark{\refstepcounter{thm} \par\medskip\noindent
                \textbf{Remark \thethm .} }

\numberwithin{equation}{thm}

\renewcommand\AA{\mathbb A} 
\newcommand\BB{\mathbb B}
\newcommand\CC{\mathbb C}

\newcommand\KK{\mathbb K}

\newcommand\QQ{\mathbb Q}

\newcommand\TT{\mathbb T}

\newcommand\VV{\mathbb V}

\newcommand\ZZ{\mathbb Z}

\newcommand\bO{\mathbf O}

\newcommand\bQ{\mathbf Q}

\newcommand\bb{\mathbf b}
\newcommand\bc{\mathbf c}
\newcommand\bd{\mathbf d}

\renewcommand\bm{\mathbf m} 

\newcommand\bq{\mathbf q}

\newcommand\cB{\mathcal{B}}

\newcommand\cE{\mathcal{E}}

\newcommand\cI{\mathcal{I}}

\newcommand\cO{\mathcal{O}}

\newcommand\cT{\mathcal{T}}

\newcommand\fS{\mathfrak S}

\newcommand\fX{\mathfrak X}

\newcommand\fg{\mathfrak g}

\newcommand\sE{\mathscr E}
\newcommand\sF{\mathscr F}

\newcommand\sH{\mathscr H}

\newcommand\sS{\mathscr S}

\renewcommand\a{\alpha}  
\renewcommand\b{\beta}   
\newcommand\g{\gamma}  
\renewcommand\d{\delta}  

\newcommand\la{\lambda}

\newcommand\s{\sigma}

\newcommand\ve{\varepsilon}

\newcommand\vf{\varphi}

\newcommand\D{\Delta}

\renewcommand\Xi{\Xi}

\newcommand\vG{\varGamma}

\newcommand\vL{\varLambda}

\newcommand\Bb{\boldsymbol\beta}

\newcommand\Be{\boldsymbol\eta}

\newcommand{\dis}{\displaystyle}

\newcommand\wh{\widehat}
\newcommand\wt{\widetilde}
\newcommand\ol{\overline}

\newcommand\ra{\rightarrow}

\newcommand\LRa{\Leftrightarrow}

\newcommand\lan{\langle}
\newcommand\ran{\rangle}

\newcommand\Ker{\operatorname{Ker}}

\newcommand\End{\operatorname{End}}

\newcommand\rank{\operatorname{rank}}

\newcommand\id{\operatorname{id}}
\newcommand\Id{\operatorname{Id}}

\renewcommand\Im{\operatorname{Im}}

\newcommand\GL{\operatorname{GL}}

\newcommand\opp{\operatorname{opp}}

\newcommand\cmod{\operatorname{-mod}}

\newcommand\Fgl{\mathfrak{gl}}
\newcommand\Fsl{\mathfrak{sl}}

\newcommand\ev{\mathbf{{ev}}}

\newcommand{\isom}{\,\raise2pt\hbox{$\underrightarrow{\sim}$}\,}

\newcounter{ichi}
\newcommand{\roi}{\roman{ichi}}
\newcounter{ni}
\newcommand{\roii}{\roman{ni}}
\newcounter{san}
\newcommand{\roiii}{\roman{san}}
\newcounter{yon}
\newcommand{\roiv}{\roman{yon}}
\newcounter{go}
\newcommand{\rov}{\roman{go}}
\newcounter{roku}
\newcommand{\rovi}{\roman{roku}}
\newcounter{nana}
\newcounter{hachi}
\newcounter{kyu}
\begin{document}



\renewcommand{\theenumi}{\roman{enumi}}
\renewcommand{\labelenumi}{(\theenumi)}
\renewcommand{\thefootnote}{\fnsymbol{footnote}}
\renewcommand{\thefootnote}{\fnsymbol{footnote}}
\parindent=20pt


\setcounter{section}{-1}




\address{Department of Mathematics, Faculty of Science, Shinshu University, 
		Asahi 3-1-1, Matsumoto 390-8621, Japan}
		
\email{kwada@shinshu-u.ac.jp}



\medskip
\begin{center}
{\large \textbf{  Schur-Weyl dualities for shifted quantum affine algebras and Ariki-Koike algebras}}  
\footnote{Keywords : Quantum groups, Hecke algebras, $q$-Schur algebras, Schur-Weyl dualities, Highest weight covers. }
\\
\vspace{1cm}
Kentaro Wada 
\\[1em]
\end{center}


\title{} 
\maketitle

\markboth{Kentaro Wada}{  Schur-Weyl dualities for shifted quantum affine algebras and Ariki-Koike algebras }

\quad\\[-6em]

\begin{abstract}
We establish a Schur-Weyl duality between a shifted quantum affine algebra and an Ariki-Koike algebra. 
Then, we realize a cyclotomic $q$-Schur algebra in the context of the Schur-Weyl duality. 
\end{abstract}


\tableofcontents



\section{Introduction} 
\para 
Let $\sH_{n,r}$ be the Ariki-Koike algebra (cyclotomic Hecke algebra) associated with the complex reflection group 
$\fS_n \ltimes (\ZZ / r \ZZ)^n$ of type $G(r,1,n)$ over the field of complex numbers $\CC$ with parameters 
$q, Q_0,Q_1,\dots, Q_{r-1} \in \CC^{\times}$. 
We denote by $\sS_{n,r} (\bm)$ the cyclotomic $q$-Schur algebra associated with $\sH_{n,r}$ introduced in \cite{DJM98}, 
where $\bm=(m_1,\dots, m_r)$ is an $r$-tuple of positive integers. 

By \cite{DJM98}, it is known that $\sS_{n,r} (\bm)$ is a quasi-hereditary cover of $\sH_{n,r}$, 
equivalently $\sS_{n,r} \cmod$ is a highest weight cover of $\sH_{n,r}\cmod$ 
in the sense of \cite{Rou08} if $\bm$ is large enough. 
Moreover, in \cite{RSVV16} and \cite{Los16} independently,  
it is known that 
$\sS_{n,r} \cmod$ is equivalent to a certain full subcategory of an affine parabolic category 
$\bO$ of the affine Lie algebra $\wh{\Fgl}_N$, 
and also equivalent to the category $\cO$ of rational Cherednik algebra of type $G(r,1,n)$ with corresponding parameters 
as highest weight covers of $\sH_{n,r} \cmod$. 
In the degenerate cases,  
they are studied in \cite{BK08} 
which also contains some connection with finite $W$-algebras. 

In this paper, 
we establish a Schur-Weyl duality between 
a shifted quantum affine algebra and an Ariki-Koike algebra, 
then we realize the cyclotomic $q$-Schur algebra via this Schur-Weyl duality. 

\para 
In the case where $r=1$, 
the algebra $\sH_{n,1}$ is the Iwahori-Hecke algebra associated with the symmetric group $\fS_n$, 
and $\sS_{n,1}(m)$ is the classical $q$-Schur algebra given in \cite{DJ89}. 
It is well-known that the $q$-Schur algebra appears in the context of quantum Schur-Weyl duality established in \cite{Jim86}
(\cite{BLM90} and \cite{Du95} in the case where $q$ is a root of unity). 
In particular, the $q$-Schur algebra $\sS_{n,1}(m)$ is a quotient algebra of the quantum group $U_q(\Fgl_m)$ 
associated with the general linear Lie algebra $\Fgl_m$. 

\para 
In the case where $r>1$, the cyclotomic $q$-Schur algebra $\sS_{n,r} (\bm)$ is defined 
as the endomorphism ring of a certain $\sH_{n,r}$-module which is  a natural generalization 
of the classical $q$-Schur algebras. 
There is a nice survey \cite{Mat04} on the representation theory of $\sS_{n,r} (\bm)$. 
Through the study of representations of $\sS_{n,r}(\bm)$, 
it had been expected that there exists a certain quantum group 
in the background of the representation theory of cyclotomic $q$-Schur algebras. 

In \cite{DR99}, it is proven that the cyclotomic $q$-Schur algebra $\sS_{n,r}(\bm)$ 
has upper and lower Borel subalgebras, 
and $\sS_{n,r}(\bm)$ is a product of them. 
Moreover, 
the upper (resp. lower) Borel subalgebra of $\sS_{n,r}(\bm)$ is isomorphic to 
the upper (resp. lower) Borel subalgebra of the classical $q$-Schur algebra $\sS_{n,1}(m)$ 
which is  a quotient of the upper (resp. lower) Borel subalgebra of the quantum group $U_q(\Fgl_m)$ 
through the quantum Schur-Weyl duality if $\bm$ is large enough. 
A difference between the cyclotomic $q$-Schur algebra $\sS_{n,r}(\bm)$ ($r>1$) 
and the classical $q$-Schur algebra $\sS_{n,1}(m)$ appears in the commutation relations 
for upper and lower Borel subalgebras (see \cite{Wad11} for details). 
We will clarify that this difference corresponds to a shift of the shifted quantum affine algebra 
in this paper. 

In \cite{Wad16}, 
by using the structures of $\sS_{n,r}(\bm)$ explained in the above,  
the author introduces the deformed current Lie algebra $\Fgl_m^{\lan \bQ \ran}[x]$ and 
the $(q, \bQ)$-current algebra $U_q (\Fgl_m^{\lan \bQ \ran} [x])$ 
associated with the general linear Lie algebra $\Fgl_m$, 
and proves that the cyclotomic $q$-Schur algebra $\sS_{n,r}(\bm)$ is realized as a quotient of $U_q (\Fgl_m^{\lan \bQ \ran} [x])$. 
The deformed current Lie algebra $\Fgl_m^{\lan \bQ \ran}[x]$ is a certain deformation of the polynomial current Lie algebra $\Fgl_m [x] = \Fgl_m \otimes \CC [x]$ with parameters $\bQ = (Q_0,Q_1,\dots, Q_{r-1}) \in (\CC^{\times})^r$, 
and The $(q,\bQ)$-current algebra $U_q (\Fgl_m^{\lan \bQ \ran}[x])$ 
is a $q$-analogue of the universal enveloping algebra of $\Fgl_m^{\lan \bQ \ran}[x]$.  
The representations  of deformed current Lie algebras $\Fgl_m^{\lan \bQ \ran} [x]$ and also $\Fsl_m^{\lan \bQ \ran}[x]$ 
are studied in \cite{Wad18}. 

In \cite{FT19}, Finkelberg and Tsymbaliuk introduce the shifted quantum affine algebra $U_{q,[\bb]} = U_{q,[\bb]}(L\Fsl_m)$ as a certain shift of quantum loop algebra $U_q (L\Fsl_m)$. 
Then, in \cite{KW21}, the authors prove that the $(q,\bQ)$-current algebra $U_q (\Fsl_m^{\lan \bQ \ran}[x])$ 
is a subalgebra of the shifted quantum affine algebra $U_{q,[\bb_{\bm}]}$ 
with the special shift $\bb_{\bm}$, 
and study the representations  of $U_q (\Fsl_m^{\lan \bQ \ran}[x])$ 
by using the connection with the shifted quantum affine algebras. 
The representations of shifted quantum affine algebras in more general setting are studied in \cite{Her23}. 

\para 
In this paper, 
we establish the Schur-Weyl duality for the shifted quantum affine algebra $U_{q,[\bb_{\bm}]} $ 
with the special shift $\bb_{\bm}$ and the Ariki-Koike algebra $\sH_{n,r}$ 
(see \eqref{choice bi} for the shift $\bb_{\bm}$). 
Then, we realize the cyclotomic $q$-Schur algebra $\sS_{n,r}(\bm)$ through this Schur-Weyl duality. 
More precisely, 
we construct the $(U_{q,[\bb_{\bm}]}, \sH_{n,r})$-bimodule $V_{n,r}$ and its submodule $\cI_{V_{n,r}}$ 
via a shift of the quantum affine Schur-Weyl duality. 
We denote the representation corresponding to the action of $U_{q,[\bb_{\bm}]}$ 
on $V_{n,r}/ \cI_{V_{n,r}}$ by $\rho_{n,r} : U_{q,[\bb_{\bm}]} \ra \End_{\CC} (V_{n,r} / \cI_{V_{n,r}})$, 
and also denote the action of $\sH_{n,r}$ by $\s_{n,r} : \sH_{n,r}^{\opp} \ra \End_{\CC} (V_{n,r} / \cI_{V_{n,r}})$. 
Then, we have the following. 
\begin{thm}[{Theorem \ref{Thm SW SQA AK}}] 
\label{Thm SW SQA AK intro}
Assume that $q$ is not a root of unity and $m_k \geq n$ for  $k=1,2,\dots, r$. 
The shifted quantum affine algebra $U_{q,[\bb_{\bm}]}$ and the Ariki-Koike algebra $\sH_{n,r}$ 
satisfy the double centralizer property through the bimodule $V_{n,r} / \cI_{V_{n,r}}$, namely we have 
\begin{align*}
\Im \rho_{n,r} = \End_{\sH_{n,r}^{\opp}} (V_{n,r} / \cI_{V_{n,r}}), 
\quad 
\Im \s_{n,r} = \End_{U_{q,[\bb_{\bm}]}} (V_{n,r} / \cI_{V_{n,r}}). 
\end{align*}
In particular, we have the surjective homomorphism 
\begin{align*}
\rho_{n,r} : U_{q,[\bb_{\bm}]} \ra \End_{\sH_{n,r}^{\opp}} (V_{n,r} / \cI_{V_{n,r}}) \cong \sS_{n,r}(\bm). 
\end{align*}  

\end{thm}

In fact, the statements in this theorem also hold even the case where $q$ is a root of unity 
by taking a certain $\AA$-form $U_{q,[\bb_{\bm}]}^{\AA}$ of $U_{q,[\bb_{\bm}]}$ and its specialization, 
where $\AA = \ZZ [\hat{q}, \hat{q}^{-1}, \hat{Q}_0, \dots, \hat{Q}_{r-1}]$ with indeterminate variables 
$\hat{q}, \hat{Q}_0,\dots, \hat{Q}_{r-1}$. 
We give these treatments in the paragraph \ref{Par AA-form}. 
The $\AA$-form $U_{q,[\bb_{\bm}]}^{\AA}$ in here is enough to consider the cyclotomic $q$-Schur algebra at a root of unity, 
but we do not whether $U_{q,[\bb_{\bm}]}^{\AA}$ is the \lq\lq correct" integral form of the shifted quantum affine algebra or not, 
and it seems not. 
(see Remark \ref{Remark integral form}). 

By Theorem \ref{Thm SW SQA AK intro}, 
we have the surjective homomorphism 
$\rho_{n,r} : U_{q,[\bb_{\bm}]} \ra \sS_{n,r}(\bm)$ 
if $\bm$ is large enough.  
This surjection is an extension of the surjection from the $(q,\bQ)$-current algebra $U_q (\Fsl_m^{\lan \bQ \ran}[x])$ 
to $\sS_{n,r}(\bm)$ given in \cite{Wad16}. 
However, the argument in this paper is quite different from that in \cite{Wad16}. 
In \cite{Wad16}, we construct the surjection directly. 
On the other hand, in this paper, 
we obtain the surjection through the Schur-Weyl duality, 
and we clarify how the cyclotomic $q$-Schur algebra $\sS_{n,r} (\bm)$ 
appears in the representation theory of the shifted quantum affine algebra 
through the Schur-Weyl duality. 

\para 
In our argument, 
we identify the set $\vL_n(m)$ of compositions of $n$ having $m$-parts 
with the set $\vL_{n,r}(\bm)$ of $r$-compositions of $n$ having $\bm$-parts, 
and we regard them as a subset of the weight lattice of $\Fgl_m$ (see \S \ref{section notation}). 

In \S \ref{section CSA}, 
we review and study on cyclotomic $q$-Schur algebras. 
In particular, 
we consider the cyclic right $\sH_{n,r}$-modules 
$M^{\mu}$ and $\wt{M}^{\mu}$ for $\mu \in \vL_{n,r} (\bm)$ (see \eqref{def Mmu} for definitions). 
The cyclotomic $q$-Schur algebra  is defined by 
\begin{align*}
\sS_{n,r}(\bm)= \End_{\sH_{n,r}^{\opp}} (\bigoplus_{\mu \in \vL_{n,r}(\bm)} M^{\mu}) 
\end{align*}  
as the endomorphism ring of the direct sum of cyclic modules $M^{\mu}$ ($\mu \in \vL_{n,r} (\bm)$). 
The cyclic module $\wt{M}^{\mu}$ ($\mu \in \vL_{n,r} (\bm)$) has a supporting role to consider the connection with Schur-Weyl dualities explained in the below. 
We give the defining relations of these cyclic modules $M^{\mu}$ and $\wt{M}^{\mu}$ 
in Theorem \ref{Thm def rel mmu}, 
and we have the surjective homomorphism $\wt{M}^{\mu} \ra M^{\mu}$ as right $\sH_{n,r}$-modules 
by using these defining relations. 
 
It seems that these defining relations had been unknown, 
and may be useful in the study of Ariki-Koike algebras and cyclotomic $q$-Schur algebras themselves.

\para 
We briefly explain how to construct 
the $(U_{q, [\bb_{\bm}]}, \sH_{n,r})$-bimodule $V_{n,r}$ and its submodule $\cI_{V_{n,r}}$ 
used in Theorem \ref{Thm SW SQA AK intro}. 

We start with the quantum Schur-Weyl duality in \cite{Jim86}. 
Let $V$ be the vector representation of the quantum group $U_q (\Fgl_m)$, 
and we consider the tensor space $V^{\otimes n}$. 
Then, the quantum group $U_q (\Fgl_m)$ acts on $V^{\otimes n}$ through the coproduct of $U_q(\Fgl_m)$, 
and the Iwahori Hecke algebra $\sH_n$ associated with the symmetric group $\fS_n$ acts on $V^{\otimes  n}$ 
as a $q$-analogue of the place permutations of tensor products.  
Then, the algebras $U_q(\Fgl_m)$ and $\sH_{n}$ satisfy the double centralizer property, 
so-called quantum Schur-Weyl duality,  
through the tensor space $V^{\otimes n}$ 
by \cite[Proposition 3]{Jim86}. 

We consider the affinization of the quantum Schur-Weyl duality 
given in \cite{GRV94} and \cite{CP96} independently. 
Let $U_q (\wh{\Fsl}_m)$ be the quantum affine algebra of type $A_{m-1}^{(1)}$. 
For $\g \in \CC^{\times}$, 
there exists an algebra homomorphism, 
so-called evaluation homomorphism, 
$\ev_{\g} : U_q (\wh{\Fsl}_m) \ra U_q(\Fgl_m)$. 
We regard the vector representation $V$ of $U_q(\Fgl_m)$ as 
the $U_q (\wh{\Fsl}_m)$-module through the evaluation homomorphism $\ev_1$ at $\g=1$, 
and we consider the affinization $\wh{V} = V \otimes \CC [x,x^{-1}]$ of $V$ 
in the sense of \cite[\S 4.2]{Kas02}. 
Then, the quantum affine algebra $U_q(\wh{\Fsl}_m)$ acts on the tensor space $\wh{V}^{\otimes n}$ 
through the coproduct of $U_q(\wh{\Fsl}_m)$, 
and the affine Hecke algebra $\wh{\sH}_n$ of type $GL_n$ also acts on $\wh{V}^{\otimes n}$ 
from right. 
Then, these actions on $\wh{V}^{\otimes n}$ commute with each other, 
and the functor 
\begin{align*} 
\wh{V}^{\otimes n} \otimes_{\wh{\sH}_n} - :  \wh{\sH}_n \cmod \ra U_q (\wh{\Fsl}_m) \cmod, 
\quad 
M \mapsto \wh{V}^{\otimes n} \otimes_{\wh{\sH}_n} M
\end{align*}
gives an equivalence between the category of finite dimensional $\wh{\sH}_n$-modules 
and the corresponding subcategory of $U_q (\wh{\Fsl}_m)\cmod$ 
by \cite[Theorem 6.8]{GRV94} and \cite[Theorem 4.2]{CP96}. 
We call it the quantum affine Schur-Weyl duality. 

Note that the Ariki-Koike algebra $\sH_{n,r}$ is a quotient of the affine Hecke algebra $\wh{\sH}_n$, 
and we put $\wt{V}_{n,r} = \wh{V}^{\otimes n} \otimes_{\wh{\sH}_n} \sH_{n,r}$. 
Then, the space $\wt{V}_{n,r}$ has the $(U_q (\wh{\Fsl}_m), \sH_{n,r})$-bimodule structure. 
The $U_q (\wh{\Fsl}_m)$-module $\wt{V}_{n,r}$ has the weight space decomposition 
$\wt{V}_{n,r} = \bigoplus_{\mu \in \vL_n(m)} \wt{V}_{n,r}^{\mu}$, 
and this is also a decomposition as right $\sH_{n,r}$-modules. 
By using the defining relation of $\wt{M}^{\mu}$ and comparing dimensions, 
we see that $\wt{V}_{n,r}^{\mu}$ is isomorphic to $\wt{M}^{\mu}$ as right $\sH_{n,r}$-modules 
for each $\mu \in \vL_n(m)$,  
where we identify the set $\vL_n(m)$ with $\vL_{n,r}(\bm)$. 
Let $\cI_{\wt{V}_{n,r}}$ be the $\sH_{n,r}$-submodule of $\wt{V}_{n,r}$ 
such that $\cI_{\wt{V}_{n,r}} \cap \wt{V}_{n,r}^{\mu}$ coincides with the kernel of the surjection 
$\wt{M}^{\mu} \ra M^{\mu}$ under the isomorphism $\wt{V}_{n,r}^\mu \cong \wt{M}^{\mu}$ 
for each $\mu \in \vL_n(m)$. 
As a consequence, we have the isomorphism 
$\wt{V}_{n,r} / \cI_{\wt{V}_{n,r}} \cong \bigoplus_{\mu \in \vL_{n,r}(\bm)} M^{\mu}$ as right $\sH_{n,r}$-modules. 

Unfortunately, 
the space $\cI_{\wt{V}_{n,r}}$ is not a $U_q(\wh{\Fsl}_m)$-submodule of $\wt{V}_{n,r}$ (Proposition \ref{Prop fi wtvmu L notin IwtV}),  
and $U_q(\wh{\Fsl}_m)$ does not act on $\wt{V}_{n,r} / \cI_{\wt{V}_{n,r}}$. 
Then, we need to consider the shifts of $U_q(\wh{\Fsl}_m)$ and $\wt{V}_{n,r}$ as follows. 
Since the action of $U_q(\wh{\Fsl}_m)$ on $\wt{V}_{n,r}$ is level zero, 
the quantum loop algebra $U_q (L \Fsl_m)$ also acts on $\wt{V}_{n,r}$ through the isomorphism given in \cite{Dri87} and \cite{Bec94}. 
We consider the shifted quantum affine algebra $U_{q,[\bb_{\bm}]}$ with the shift $\bb_{\bm}$ corresponding to $\bm$ 
(see \eqref{choice bi} for the shift $\bb_{\bm}$), 
and consider the shift homomorphism $\iota_{[0, \bb_{\bm}]}^{\Be} : U_{q,[\bb_{\bm}]} \ra U_q (L \Fsl_m)$. 
(Correctly, we need to take a certain quotient of $U_{q,[\bb_{\bm}]}$ associated with $\Be$, 
or replace $U_q(L\Fsl_m)$ with $U_{q,[0]}$, see Remark \ref{Remark iota eta cd}. 
We do not care in this introduction since it is not so much trouble.)
By \cite{FT19}, it is known that $U_{q,[\bb_{\bm}]}$ is a coideal subalgebra of $U_q (L\Fsl_m)$ 
through the shift homomorphism $\iota_{[0,\bb_{\bm}]}^{\Be}$. 
In particular, we have the algebra homomorphism 
$\D_{\bb_{\bm},0} : U_{q,[\bb_{\bm}]} \ra U_{q,[\bb_{\bm}]} \otimes U_q (L\Fsl_m)$ 
induced from Drinfeld-Jimbo coproduct of $U_q (L\Fsl_m)$. 
We also take a one-dimensional $U_{q,[\bb_{\bm}]}$-module $L_{\Bb_{\bm}^{\bQ}}$ with respect to 
parameters $\bQ =(Q_0,Q_1,\dots, Q_{r-1})$ 
(see \eqref{act w}). 
Put $V_{n,r} = L_{\Bb_{\bm}^{\bQ}} \otimes \wt{V}_{n,r}$. 
Then, $U_{q,[\bb_{\bm}]}$ acts on $V_{n,r}$ through the algebra homomorphism $\D_{\bb_{\bm},0}$, 
and $\sH_{n,r}$ also acts on $V_{n,r}$ from right via the action on $\wt{V}_{n,r}$. 
As a consequence, we obtain the $(U_{q,[\bb_{\bm}]}, \sH_{n,r})$-bimodule $V_{n,r}$. 
We can also take an $\sH_{n,r}$-submodule $\cI_{V_{n,r}}$ of $V_{n,r}$ in a similar manner 
as in the case of $\wt{V}_{n,r}$, 
and we see that $\cI_{V_{n,r}}$ is also a $U_{q,[\bb_{\bm}]}$-submodule of $V_{n,r}$ 
(Corollary \ref{Cor cIVnr Uqbm sub}). 
Then, we obtain the Schur-Weyl duality in Theorem \ref{Thm SW SQA AK intro} 
through the $(U_{q,[\bb_{\bm}]}, \sH_{n,r})$-bimodule $V_{n,r} / \cI_{V_{n,r}}$. 

\para 
The construction of the $(U_{q,[\bb_{\bm}]}, \sH_{n,r})$-bimodule $V_{n,r} / \cI_{V_{n,r}}$ is a bit technical. 
However, in the case where $q=1$, we have more natural interpretation of the bimodule $V_{n,r} / \cI_{V_{n,r}}$ as follows. 

In \S \ref{Section shifted loop Lie alg}, 
we introduce the shifted loop Lie algebra $L_{[\bb]}^{\Be} \Fsl_m$ with a shift $\bb \in \ZZ_{\geq 0}^{m-1}$.
This is a shift of the loop Lie algebra $L \Fsl_m$ associated with $\Fsl_m$. 
In particular, we have the shift homomorphism 
$\iota_{[\bb]^{\sharp}}^{\Be} : L_{[\bb]}^{\Be} \Fsl_m \ra L\Fsl_m$ (Proposition \ref{Prop shift hom Lie}). 
We can regard the universal enveloping algebra $U(L_{[\bb]}^{\Be} \Fsl_m)$ of $L_{[\bb]}^{\Be} \Fsl_m$ 
as the shifted quantum affine algebra $U_{q,[\bb]}$ at $q=1$ by taking a certain form $\wt{U}_{q, \Be, [\bb]}$ of $U_{q,[\bb]}$. 

We denote by $\sH_{n,r}^{q=1}$ the Ariki-Koike algebra at $q=1$, 
and we denote by $U_{[\bb_{\bm}]}$ 
the universal enveloping algebra of the shifted loop Lie algebra $L_{[\bb_{\bm}]}^{\Be} \Fsl_m$ with the shift $\bb_{\bm}$. 
Then, we can construct the $(U_{[\bb_{\bm}]}, \sH_{n,r}^{q=1})$-bimodule $V_{n,r} / \cI_{V_{n,r}}$ in a similar way as in the above. 
On the other hand, the algebra $U_{[\bb_{\bm}]}$ is a Hopf subalgebra of $U (L \Fsl_m)$ through the shift homomorphism 
$\iota_{[\bb]^{\sharp}}^{\Be} : L_{[\bb]}^{\Be} \Fsl_m \ra L\Fsl_m$. 
Put $\VV = V_{1,r}/ \cI_{V_{1,r}}$ the bimodule in the case where $n=1$. 
We consider the tensor space $\VV^{\otimes n}$, 
and the algebra $U_{[\bb_{\bm}]}$ acts on $\VV^{\otimes n}$ through the coproduct of $U_{[\bb_{\bm}]}$. 
We can also define the right action of $\sH_{n,r}^{q=1}$ on $\VV^{\otimes n}$ in a natural way. 
Then, we can prove that $ \VV^{\otimes n} \cong V_{n,r} / \cI_{V_{n,r}} $ as $(U_{[\bb_{\bm}]}, \sH_{n,r}^{q=1})$-bimodules, 
and we obtain the statements corresponding to  Theorem \ref{Thm SW SQA AK intro} at $q=1$ 
through the tensor space $\VV^{\otimes n}$ 
(Theorem \ref{Thm SW SLA AK}). 

By the above argument at $q=1$,  
we can regard $(U_{q,[\bb_{\bm}]}, \sH_{n,r})$-bimodule $V_{n,r} / \cI_{V_{n,r}}$ 
as a $q$-analogue of the tensor space $\VV^{\otimes n}$. 
However, we can not consider the tensor space $\VV^{\otimes n}$ as the $U_{q,[\bb_{\bm}]}$-module directly 
since the shifted quantum affine algebra $U_{q,[\bb_{\bm}]}$ does not have a Hopf algebra structure although 
it is a coideal subalgebra of $U(L\Fsl_m)$.  
Then, as explained in the above, 
we need  a bit technical steps to construct the bimodule $V_{n,r} / \cI_{V_{n,r}}$, 
but it seems that the bimodule $V_{n,r} / \cI_{V_{n,r}}$ is natural one in the representation theory of shifted quantum affine algebras. 
\\

{\bf Acknowledgements:} 
The author was supported by JSPS KAKENHI Grant Number JP21K03178. 




\section{Notation}  
\label{section notation}
In this section, we give some notation used in this paper. 

\para 
For a condition $X$, 
we put $\d_{(X)} =1$ if $X$ is true, and $\d_{(X)} =0$ if $X$ is false. 
For integers $i$ and $j$, put $\d_{i,j} = \d_{(i=j)}$. 

For an integer $a \in \ZZ$, 
we put $\ZZ_{\geq a} = \{ b \in \ZZ \mid b \geq a\}$. 
It is similar for $\ZZ_{>a}$, $\ZZ_{\leq a}$ and $\ZZ_{<a}$. 

For an algebra $A$, we denote the opposite algebra of $A$ by $A^{\opp}$. 
We also denote the category of finitely generated left $A$-modules by $A \cmod$. 

\para 
For an $r$-tuple of positive integers $\bm=(m_1,m_2, \dots, m_r) \in \ZZ_{>0}^r$, 
put 
\begin{align*}
\vG(\bm) =\{(i,k) \in \ZZ^2 \mid 1 \leq i \leq m_k, \, 1 \leq k \leq r\}.
\end{align*}
We identify the set $\vG(\bm)$ with the set $\{1,2, \dots, m\}$, 
where $m=m_1+m_2+\dots +m_r$, 
by the bijection 
\begin{align*}
\xi : \vG(\bm) \ra \{1,2,\dots, m\}, 
\quad 
(i,k) \mapsto \sum_{p=1}^{k-1} m_p +i. 
\end{align*}
Then we can regard the set $\vG(\bm)$ as a partion of the set $\{1,2,\dots, m\}$ to $r$-parts. 


\para
For an $r$-tuple of positive integers $\bm=(m_1,m_2, \dots, m_r) \in \ZZ_{>0}^r$, 
and the integer $m=m_1+m_2+\dots+m_r$,  
put 
\begin{align*}
&\vL_n(m) = \{ \mu = (\mu_1, \mu_2,\dots, \mu_m) \in \ZZ_{\geq 0}^m \mid \sum_{i=1}^m \mu_i =n \}, 
\\
&\vL_{n,r}(\bm) = \left\{ \mu=(\mu^{(1)}, \mu^{(2)}, \dots, \mu^{(r)} ) \mid 
	\begin{array}{l}
	\mu^{(k)} =(\mu_1^{(k)}, \mu_2^{(k)}, \dots, \mu_{m_k}^{(k)}) \in \ZZ_{\geq 0}^{m_k} 
	\\
	\sum_{k=1}^r \sum_{i=1}^{m_k} \mu_i^{(k)} =n
	\end{array}
	\right\}. 
\end{align*}
An element of $\vL_n(m)$ is called a composition of $n$ with $m$ parts, 
and an element of $\vL_{n,r}(\bm)$ is called an $r$-composition.  
For $\mu =(\mu^{(1)}, \mu^{(2)}, \dots, \mu^{(r)}) \in \vL_{n,r}(\bm)$, 
we denote by $|\mu^{(k)}|$ the sum of integers $\sum_{i=1}^{m_k} \mu_i^{(k)}$.

We identify the set $\vL_{n,r}(\bm)$ with the set $\vL_n(m)$ by the bijection 
\begin{align}
\label{bij vLnrbm vLnm}
\begin{split}
&\vL_{n,r}(\bm) \ra \vL_n(m), 
\\
&(\mu^{(1)}, \mu^{(2)} \dots, \mu^{(r)}) \mapsto (\mu_1^{(1)}, \mu_2^{(1)}, \dots, \mu_{m_1}^{(1)}, 
\mu_1^{(2)}, \dots, \mu_{m_2}^{(2)}, \dots,  
\mu_{1}^{(r)}, \dots, \mu_{m_r}^{(r)}).
\end{split}
\end{align}

Let $P= \bigoplus_{i=1}^m \ZZ \ve_i$ be the weight lattice of $\Fgl_m$. 
Put $\a_i = \ve_i - \ve_{i+1}$ ($1\leq i \leq m-1$), 
then $\{\a_i \mid 1\leq i \leq m-1\}$ gives the set of simple roots of $\Fgl_m$ (resp. $\Fsl_m$). 
We regard the set $\vL_n(m)$ as a subset of the weight lattice $P$ 
by the injection 
\begin{align}
\label{inj vnm P}
\vL_n(m) \ra P,  \quad 
(\mu_1,\mu_2, \dots, \mu_m) \mapsto \sum_{i=1}^m \mu_i \ve_i. 
\end{align}
We also regard  the set $\vL_{n,r}(\bm)$ as a subset of the weight lattice $P$ 
through the bijection \eqref{bij vLnrbm vLnm} and the injection \eqref{inj vnm P}. 

\para 
For $\mu =(\mu_1, \mu_2, \dots, \mu_m) \in \vL_{n}(m) = \vL_{n,r}(\bm)$ 
and $1\leq i \leq m$, 
put 
\begin{align*}
N_i^{\mu} = \mu_1 + \mu_2 + \dots + \mu_i. 
\end{align*}
For $\mu \in \vL_{n,r}(\bm)$ and $0 \leq k \leq r$,  put 
\begin{align*}
a_k^{\mu} = \sum_{j=1}^k |\mu^{(j)}|. 
\end{align*}
For $\mu \in \vL_{n,r} (\bm)$ and $1 \leq j \leq n$, 
we define the integer $c_j^{\mu}$ by 
\begin{align}
\label{def cjmu} 
c_j^{\mu} = k \text{ if } a_{k-1}^{\mu} < j \leq a_k^{\mu}. 
\end{align}
From the definition, we see that 
\begin{align*}
c_{j}^{\mu} \leq c_{j'}^{\mu} \text{ if } j \leq j' 
\text{ and }
c_{N_{i}^{\mu} +1}^{\mu} = c_{N_i^{\mu} +2}^{\mu} = \dots = c_{N_i^{\mu} + \mu_{i+1}}^{\mu}.
\end{align*}


\section{Ariki-Koike algebras and affine Hecke algebras} 
\label{section AK}
In this section, 
we review the definition and some known basic properties of Ariki-Koike algebras. 
We also introduce some elements of the Ariki-Koike algebra which are used in the subsequence argument. 

\para 
Let $R$ be a commutative ring.  
We take parameters $q, Q_0, Q_1, \dots, Q_{r-1} \in R^{\times}$, 
where $R^{\times}$ is the set of invertible elements in $R$, 
and put $\bQ = (Q_0,Q_1,\dots, Q_{r-1})$. 
The Ariki-Koike algebra (cyclotomic Hecke algebra) $\sH_{n,r} = \sH_{n,r} (q, \bQ)$ 
associated with the complex reflection group $\fS_n \ltimes (\ZZ / r \ZZ)^n$ of type $G(r,1,n)$ 
is the associative algebra over $R$ generated by $T_0, T_1, \dots, T_{n-1}$ 
with the following defining relations: 
\begin{align*}
&(T_0 - Q_0) (T_0 - Q_1) \dots (T_0 - Q_{r-1}) =0, 
\quad 
(T_i - q) (T_i + q^{-1}) =0 \quad  (1 \leq i \leq n-1), 
\\
&T_0 T_1 T_0 T_1 = T_1 T_0 T_1 T_0, 
\quad 
T_i T_{i+1} T_i = T_{i+1} T_i T_{i+1} \quad (1 \leq i \leq n-2), 
\\
& T_i T_j = T_j T_i \quad ( |i-j| >1). 
\end{align*}

We remark that there exists an isomorphism of algebras 
\begin{align}
\label{iso sHnr sHnropp}
\sH_{n,r} \ra \sH_{n,r}^{\opp} 
\text{ such that } 
T_i \mapsto T_i \, (0 \leq i \leq n-1).
\end{align}

\para
The subalgebra of $\sH_{n,r}$ generated by $T_1, T_2,\dots, T_{n-1}$ is isomorphic to the Iwahori-Hecke algebra 
$\sH_n$ associated with the symmetric group $\fS_n$ of $n$ letters. 

For $i=1,2,\dots, n-1$, 
we denote the adjacent transposition by $s_i =(i,i+1)$. 
For $w \in \fS_n$, we denote the length of $w$ by $\ell (w)$. 
Let $\{ T_w \mid w \in \fS_n \}$ be the standard basis of $\sH_n$, 
namely $T_w = T_{i_1} T_{i_2} \dots T_{i_l}$ when $w=s_{i_1} s_{i_2} \dots s_{i_l}$ is a reduced expression. 

For $\mu =(\mu_1, \mu_2,\dots, \mu_m) \in \vL_n(m)$, let $\fS_{\mu}$ be the parabolic subgroup of $\fS_n$ corresponding to $\mu$,  
namely $\fS_{\mu} \cong \fS_{\mu_1} \times \fS_{\mu_2} \times \dots \times \fS_{\mu_m}$ 
is the subgroup of $\fS_n$ generated by 
\begin{align*}
\bigcup_{i=1}^m \{ s_{(\mu_1+\dots+ \mu_{i-1}) +1}, s_{(\mu_1+\dots+ \mu_{i-1}) +2}, \dots, 
	s_{(\mu_1+\dots+ \mu_{i-1}) +\mu_i -1} \}. 
\end{align*}
Put $\fS^{\mu} = \{ y \in \fS_n \mid \ell (x y ) = \ell (x) + \ell (y) \text{ for all } x \in \fS_{\mu}\}$, 
then it is well-known that $\fS^{\mu}$ gives the set of distinguished coset representatives of the coset $\fS_{\mu} \backslash \fS_n$. 
Moreover, for $w \in \fS_n$, we can uniquely write 
\begin{align}
\label{Tw Tx Ty}
T_w = T_x T_y \text{ for some } x \in \fS_{\mu} \text{ and } y \in \fS^{\mu}. 
\end{align}
 
\para 
For $1 \leq i \leq j \leq n-1$, put 
$\TT_{i,j}= T_i T_{i+1} \dots T_j$,  
$\wt{\TT}_{j,i} = T_j T_{j-1} \dots T_i$.  
For convenience, 
we also put $ \TT_{i,i-1} = \wt{\TT}_{i-1, i}=1$ for $i=1,2,\dots, n-1$. 

For $i=1,2,\dots,n$, 
we define an element $L_i \in \sH_{n,r}$ by $L_i = \wt{\TT}_{i-1,1} T_0 \TT_{1,i-1}$. 
We see that $T_i$ ($1\leq i \leq n-1$) and $L_j$ ($1\leq j \leq n$) are invertible elements of $\sH_{n,r}$ 
by definitions. 
The following relations are well-known, 
and one can check them by direct calculation using the defining relations of $\sH_{n,r}$.  
\begin{lem}
\label{Li Lj} 
We have the following.  
\begin{enumerate}
\item 
$L_i$ and $L_j$ commute with each other for any $1 \leq i,j \leq n$.  

\item 
$T_i$ and $L_j$ commute with each other if $j \not=i,i+1$. 
\\
$T_i L_i = L_{i+1} T_i - (q-q^{-1}) L_{i+1}$, 
$T_i L_{i+1} = L_i T_i + (q-q^{-1}) L_{i+1}$. 

\item 
$T_i$ commutes with both $L_i L_{i+1}$ and $L_i + L_{i+1}$. 

\item 
For any $a \in R$,  
$(L_1-a)(L_2-a) \dots (L_i-a)$ and $T_j$ commute with each other if $i \not=j$. 

\item 
$L_{i+1}^t  T_i  = T_i L_i^t + (q-q^{-1}) \sum_{s=0}^{t-1} L_i^{s} L_{i+1}^{t-s}$ 
for $1 \leq i \leq n-1$ and $t \in \ZZ_{>0}$. 
\\
$L_i^t T_i = T_i L_{i+1}^t - (q-q^{-1}) \sum_{s=1}^t L_i^{t-s} L_{i+1}^s$ 
for $1 \leq i \leq n-1$ and $t \in \ZZ_{>0}$. 

\item 
$L_i^k L_{i+1}^l T_i = \begin{cases} 
	T_i L_i^l L_{i+1}^k + (q-q^{-1}) \sum_{s=1}^{l-k} L_i^{l-s} L_{i+1}^{k+s} & \text{ if } k \leq l, 
	\\
	T_i L_i^l L_{i+1}^k - (q-q^{-1}) \sum_{s=1}^{k-l} L_i^{k-s} L_{i+1}^{l+s} & \text{ if } k >l 
	\end{cases}$
	\\
	for $1\leq i \leq n-1$. 
\end{enumerate}
\end{lem}

We have the following theorem for an $R$-free basis of $\sH_{n,r}$ proved  in \cite[Theomre 3.10]{AK94}.
\begin{thm}[{\cite[Theomre 3.10]{AK94}}]
\label{AK} 
The algebra $\sH_{n,r}$ is an $R$-free module with an $R$-free  basis 
$\{  T_w L_1^{k_1} L_2^{k_2} \dots L_n^{k_n}  \mid w \in \fS_n, \, 0 \leq k_i \leq r-1 \, (1 \leq i \leq n) \}$. 
The set $\{ L_1^{k_1} L_2^{k_2} \dots L_n^{k_n} T_w  \mid w \in \fS_n, \, 0 \leq k_i \leq r-1 \, (1 \leq i \leq n) \}$ 
also gives an $R$-free basis of $\sH_{n,r}$. 
\end{thm}

\para 
The affine Hecke algebra $\wh{\sH}_n$ of type $\GL_n$ is the associative algebra 
over $R$ generated by $T_1, T_2, \dots, T_{n-1}$ and $X_1^{\pm}, X_2^{\pm}, \dots, X_n^{\pm}$ 
with the following defining relations: 
\begin{align*}
&(T_i - q) (T_i + q^{-1}) =0 \quad (1\leq i \leq n-1),  
\\
& T_i T_{i+1} T_i = T_{i+1} T_i T_{i+1} \quad (1\leq i \leq n-2),  
\quad 
T_i T_j = T_j T_i \quad (|i-j| >1), 
\\
& X_i^+ X_i^- = X_i^- X_i^+=1, \quad X_i^+ X_j^+ = X_j^+ X_i^+ \quad (1\leq i,j \leq n), 
\\
& T_i X_i^+ T_i = X_{i+1}^+ \quad ( 1 \leq i \leq n-1),  
\quad 
T_i X_j^+ = X_j^+ T_i \quad (j \not=i,i+1). 
\end{align*}

By the defining relations of $\sH_{n,r}$ 
together with 
the definition of $L_j$ and Lemma \ref{Li Lj} (\roii), 
we see that there exists the surjective algebra homomorphism 
\begin{align}
\label{surj whH H}
\wh{\sH_n} \ra \sH_{n,r} 
\text{ such that } 
T_i \mapsto T_i, 
\quad 
X_j^+ \mapsto L_j, 
\quad 
X_j^- \mapsto L_j^{-1}.
\end{align}

\para 
For later discussion,  
we introduce some modifications of the elements $L_i^k \in \sH_{n,r}$ as follows.  
For $1\leq i \leq n$ and $ 1\leq k \leq r$,  
we define $L_i^{\lan k \ran} \in \sH_{n,r}$ by 
\begin{align*}
L_i^{\lan k \ran} = \wt{\TT}_{i-1,1} (T_0- Q_0) (T_0-Q_1) \dots (T_0- Q_{k-1}) \TT_{1,i-1}. 
\end{align*}
By the defining relations of $\sH_{n,r}$, 
we see that $L_i^{\lan r \ran} =0$ for any $i=1,2,\dots, n$. 
We have $T_i L_i^{\lan k \ran} T_i = L_{i+1}^{\lan k \ran}$ for $1\leq i \leq n-1$ and $1\leq k \leq n$.  
We also have 
\begin{align}
\label{Llan l ran L lan j ran}
L_i^{\lan k \ran} 
&= L_i^{\lan l \ran}  (\TT_{1,i-1})^{-1} ( T_0 - Q_l) (T_0 - Q_{l+1}) \dots (T_0- Q_{k-1}) \TT_{1,i-1}
\end{align}
for $1\leq i \leq n$ and $1 \leq l < k \leq r$. 
We note that $L_i^{\lan k \ran}$ and $L_j^{\lan l \ran}$ ($1\leq i,j \leq n$, $1\leq k,l \leq r$) do not commute in general.  


\remark  
In the case where $q=1$, we have $T_i^2=1$ for $i=1,2,\dots, n-1$. 
Then we have $\TT_{1,i-1} \wt{\TT}_{i-1,1} = \wt{\TT}_{i-1,1} \TT_{1,i-1}=1$. 
This implies that, if $q=1$, 
\begin{align*}
L_i^{\lan k \ran} = (L_i - Q_0) (L_i-Q_1) \dots (L_i - Q_{k-1})
\end{align*} 
for $1\leq i \leq n$ and $1\leq k \leq r$. 

In the last of this section, we give some technical lemmas for the elements $L_i^{\lan k \ran}$. 
We give a proof of these lemmas in Appendix \ref{Appendix proof technical lemmas}. 
\begin{lem}
\label{Li lan k ran}
For $1\leq i \leq n$ and $ 1 \leq k \leq r$, 
we have  
\begin{align*}
L_i^{\lan k \ran} = L_i^k + \sum_{0 \leq p_1,p_2,\dots, p_i \leq k-1} L_1^{p_1} L_2^{p_2} \dots L_i^{p_i} h_{(p_1,p_2,\dots, p_i)} 
\end{align*}
for some $h_{(p_1,p_2, \dots, p_i)} \in \sH_{[1,i]}$, where 
$\sH_{[1,i]}$ is the subalgebra of $\sH_{n,r}$ generated by $\{T_1, T_2,\dots, T_{i-1}\}$ 
with $\sH_{[1,1]} = \CC$. 
\end{lem}

\begin{lem}
\label{Lemma TT Lk}
For $1\leq i \leq j \leq n-1$, $1\leq l \leq n$ and $1 \leq k \leq r$, we have the following. 
\begin{enumerate}
\item  
If $l \not= i, i+1$, we have 
$T_i L_l^{\lan k \ran} = L_l^{\lan k \ran} T_i$. 

\item 
$T_i L_i^{\lan k \ran} = L_{i+1}^{\lan k \ran} T_i - (q-q^{-1}) L_{i+1}^{\lan k \ran}$.

\item 
$T_i L_{i+1}^{\lan k \ran} = L_i^{\lan k \ran}T_i + (q-q^{-1}) L_{i+1}^{\lan k \ran}$. 

\item 
$\TT_{i,j} L_{l}^{\lan k \ran} 
	= \begin{cases}
		L_l^{\lan k \ran} \TT_{i,j} 
			& \text{ if } l > j+1, 
		\\ \dis 
		L_i^{\lan k \ran} \TT_{i,j} + (q-q^{-1}) \sum_{p=1}^{j-i+1} L_{i+p}^{\lan k \ran} \TT_{i, i+p-2} \TT_{i+p, j} 
			& \text{ if } l=j+1, 
		\\ 
		L_{l+1}^{\lan k \ran} \big( \TT_{i,j} - (q-q^{-1}) \TT_{i, l-1} \TT_{l+1,j} \big) 
			& \text{ if } j+1 > l \geq i, 
		\\
		L_l^{\lan k \ran} \TT_{i,j} 
			& \text{ if } i >l.
	\end{cases} $ 

\item 
$ \wt{\TT}_{j,i} L_l^{\lan k \ran} 
	= \begin{cases}
		L_l^{\lan k \ran} \wt{\TT}_{j,i} 
			& \text{ if } l > j+1, 
		\\
		L_{l-1}^{\lan k \ran} \wt{\TT}_{j,i} + (q-q^{-1}) L_{j+1}^{\lan k \ran} (\TT_{l,j})^{-1} \wt{\TT}_{l-2,i} 
			& \text{ if } j +1 \geq l >i, 
		\\
		L_{j+1}^{\lan k \ran} (\TT_{i,j})^{-1} 
			& \text{ if } l=i, 
		\\
		L_l^{\lan k \ran} \wt{\TT}_{j,i} 
			& \text{ if } i >l.
	\end{cases} $ 
\end{enumerate}
\end{lem}


\begin{lem}
\label{Lemma Li Ljlan k ran} 
For $1 \leq i,j \leq n$ and $1\leq k \leq r$, we have  the following. 
\begin{enumerate}
\item 
\begin{enumerate} 
\item 
If $ i >j$, we have 
$L_i L_j^{\lan k \ran} = L_j^{\lan k \ran} L_i$.  
\item 
If $i=j$, we have  
\begin{align*}
L_j L_j^{\lan k \ran} 
= L_j^{\lan k \ran} L_j -  ( q-q^{-1}) \sum_{p=1}^{j-1} (L_j^{\lan k \ran} L_{j-p} - L_{j-p}^{\lan k \ran} L_j) 
	(\TT_{j-p+1, j-1})^{-1} \TT_{j-p, j-1}. 
\end{align*}
\item 
If $i <j$, we have 
\begin{align*}
L_i L_j^{\lan k \ran} = L_j^{\lan k \ran} L_i + (q-q^{-1}) (L_j^{\lan k \ran} L_i - L_i^{\lan k \ran} L_j) (\TT_{i+1, j-1})^{-1} \TT_{i, j-1}.  
\end{align*}
\end{enumerate} 

\item 
\begin{enumerate} 
\item 
If $i >j$, we have 
$L_i^{-1} L_j^{\lan k \ran} = L_j^{\lan k \ran} L_i^{-1}$. 

\item 
If $i=j$, we have 
\begin{align*}
L_j^{-1} L_j^{\lan k \ran} 
&= L_j^{\lan k \ran} L_j^{-1} \wt{\TT}_{j-1,1} \TT_{1, j-1} 
	- (q-q^{-1}) \sum_{p=1}^{j-1} L_{j-p}^{\lan k \ran} L_{j-p}^{-1} (\TT_{j-p+1, \, j-1})^{-1} \wt{\TT}_{j-p-1, \, 1} \TT_{1, \, j-1}. 
\end{align*}

\item 
If $ i < j$, we have 
\begin{align*}
L_i^{-1} L_j^{\lan k \ran} 
&= L_j^{\lan k \ran} L_i^{-1} - (q-q^{-1}) L_j^{\lan k \ran} L_j^{-1} \wt{\TT}_{j-1, \, i+1} \TT_{i, \, j-1} 
	\\ & \quad 
	+ (q-q^{-1}) L_i^{\lan k \ran} L_i^{-1} \wt{\TT}_{j-1, \, i+1} \wt{\TT}_{i-1, \, 1} \TT_{1, \, j-1}
	\\ & \quad 
	- (q-q^{-1})^2 \sum_{p=1}^{i-1} L_{i-p}^{\lan k \ran} L_{i-p}^{-1} \wt{\TT}_{j-1, \, i+1} \wt{\TT}_{i - p -1 , \, 1} 
		( \TT_{i - p +1, \, i-1})^{-1} \TT_{1, \, j-1}. 
\end{align*}
\end{enumerate} 
\end{enumerate}
\end{lem}

\begin{lem}
\label{Lemma Lj Qk LJlankran}
For $1 \leq j \leq n $ and $ 1 \leq k \leq r-1 $, 
we have 
\begin{align*}
(L_j - Q_k) L_j^{\lan k \ran} 
= L_j^{\lan k+1 \ran} + (q-q^{-1}) \sum_{p=1}^{j-1} L_p^{\lan k \ran} \wt{\TT}_{j-1, \, p+1} L_{p+1} \TT_{p, j-1}.
\end{align*}
\end{lem}


\section{Cyclotomic $q$-Schur algebras} 
\label{section CSA}

In this section, 
we review the definition and some known basic properties of cyclotomic $q$-Schur algebras.
Then, we give a defining relation of cyclic $\sH_{n,r}$-module $M^{\mu}$ (resp. $\wt{M}^{\mu}$) 
in Theorem \ref{Thm def rel mmu}. 
These defining relations are important to construct the Schur-Weyl duality in the subsequence sections.

\para 
For $\mu \in \vL_{n,r}(\bm)$,  we define the elements $x_\mu$ and $m_{\mu}$ of $\sH_{n,r}$ by 
\begin{align*}
x_{\mu} = \sum_{ w \in \fS_{\mu}} q^{\ell (w)} T_w, 
\quad 
m_{\mu} = \big( \sum_{w \in \fS_{\mu}} q^{\ell(w)} T_w \big) \big( \prod_{k=1}^{r-1} \prod_{i=1}^{a_k^{\mu}} (L_i - Q_k) \big).  
\end{align*}
It is well-known (e.g. \cite[Lemma 3.2]{Mat99}) that 
\begin{align}
\label{xmu Tw}
x_{\mu} T_w = q^{\ell (w)} x_{\mu} 
\text{ and }
m_{\mu} T_w = q^{\ell (w)} m_{\mu} 
\text{ for } w \in \fS_{\mu}. 
\end{align}
We denote the right ideal of $\sH_{n,r}$ generated by $m_{\mu}$ (resp. $x_{\mu}$) by $M^{\mu}$ (resp. $\wt{M}^{\mu}$),  
namely we have 
\begin{align}
\label{def Mmu}
M^{\mu} = m_{\mu} \sH_{n,r}  
\text{ and } 
\wt{M}^{\mu} = x_{\mu} \sH_{n,r}.
\end{align}  
Then, the cyclotomic $q$-Schur algebra $\sS_{n,r}(\bm)$ associated with $\sH_{n,r}$ is defined by  
\begin{align*}
\sS_{n,r}(\bm) = \End_{\sH_{n,r}^{\opp}} \big( \bigoplus_{\mu \in \vL_{n,r}(\bm)} M^{\mu} \big).  
\end{align*}
From the definitions, we see that $\bigoplus_{\mu \in \vL_{n,r}(\bm)} M^{\mu}$ is an $(\sS_{n,r}(\bm), \sH_{n,r})$-bimodule. 
The cyclotomic $q$-Schur algebra has the following property.  


\begin{thm}[{\cite[Theorem 6.12, Corollary 6.18]{DJM98}, \cite[Theorem 5.3]{Mat04}}]  
\label{Thm CSA SW}
Suppose that $R$ is a field, and $m_k \geq n$ for $k=1,2,\dots, r$, 
then we have the following. 
\begin{enumerate}
\item 
The cyclotomic $q$-Schur algebra $\sS_{n,r}(\bm)$ is a quasi-hereditary cellular algebra.  

\item 
The cyclotomic $q$-Schur algebra $\sS_{n,r}(\bm)$ and the Ariki-Koike algebra $\sH_{n,r}$ 
satisfy the double centralizer property through the bimodule $\bigoplus_{\mu \in \vL_{n,r}(\bm)} M^{\mu}$, namely we have 
\begin{align*}
\sS_{n,r}(\bm) = \End_{\sH_{n,r}^{\opp}} \big( \bigoplus_{\mu \in \vL_{n,r}(\bm)} M^{\mu} \big) 
\text{ and }
\sH_{n,r} \cong \End_{\sS_{n,r}(\bm)} \big( \bigoplus_{\mu \in \vL_{n,r}(\bm)} M^{\mu} \big). 
\end{align*}
\end{enumerate}
\end{thm}


\begin{lem}
\label{Lemma mmu L1k1 Lj-1kj-1 Lj lan li ran}
For $\mu \in \vL_{n,r}(\bm)$, $1\leq j \leq n$ and $k_1,k_2, \dots, k_{j-1} \geq 0$, 
we have  
\begin{align*}
m_{\mu} L_1^{k_1} L_2^{k_2} \dots L_{j-1}^{k_{j-1}} L_j^{\lan c_j^{\mu} \ran} =0. 
\end{align*}
In particular, 
we have 
$m_{\mu} L_j^{\lan c_j^{\mu} \ran} =0$ for $1 \leq j \leq n$. 

\begin{proof}
By the definition of $m_{\mu}$ together with Lemma \ref{Li Lj}  (\roi), 
we have 
\begin{align*}
&m_{\mu} L_1^{k_1} L_2^{k_2} \dots L_{j-1}^{k_{j-1}} L_j^{\lan c_j^{\mu} \ran} 
\\
&=x_{\mu}  L_1^{k_1} L_2^{k_2} \dots L_{j-1}^{k_{j-1}} \big( \prod_{k=1}^{c_j^{\mu}-1} \prod_{i=1}^{a_k^{\mu}} (L_i - Q_k) \big) 
	\big( \prod_{k= c_j^{\mu}}^{r-1} \prod_{i=1}^{a_k^{\mu}} (L_i - Q_k) \big) L_j^{\lan c_j^{\mu} \ran}. 
\end{align*}
Note that $k \geq c_j^{\mu}$ implies that $a_k^{\mu} \geq j$ by the definition \eqref{def cjmu}, 
and we have 
\begin{align*}
&\big( \prod_{k= c_j^{\mu}}^{r-1} \prod_{i=1}^{a_k^{\mu}} (L_i - Q_k) \big) L_j^{\lan c_j^{\mu} \ran} 
\\
&= \wt{\TT}_{j-1,1} \big( \prod_{k=c_j^{\mu}}^{r-1} \prod_{i=1}^{a_k^{\mu}} (L_i - Q_k) \big) 
	\big( \prod_{ p=0}^{c_j^{\mu} -1} (T_0-Q_p) \big) \TT_{1, j-1} 
\\
&= \wt{\TT}_{j-1,1} \big( \prod_{k=c_j^{\mu}}^{r-1} \prod_{i=2}^{a_k^{\mu}} (L_i - Q_k) \big) 
	\big( \prod_{k=c_j^{\mu}}^{r-1} (L_1 - Q_k ) \big) 
	\big( \prod_{ p=0}^{c_j^{\mu} -1} (T_0-Q_p) \big) \TT_{1, j-1} 
=0, 
\end{align*}
by Lemma \ref{Li Lj} (\roiv). 
Then, we obtain the claim of the lemma. 
\end{proof} 
\end{lem}



\begin{prop}
\label{Prop basis Mmu}
For $\mu \in \vL_{n,r}(\bm)$, we have the following. 
\begin{enumerate}
\item 
$\wt{\cB}^{\mu} = \{ x_{\mu} L_1^{k_1} L_2^{k_2} \dots L_n^{k_n} T_y 
	\mid 0 \leq k_i \leq r-1 \, (1\leq i \leq n), \, y \in \fS^{\mu}\} $ 
is an $R$-free basis of $\wt{M}^{\mu}$. 

\item 
$\cB^{\mu} = \{m_{\mu} L_1^{p_1} L_2^{p_2} \dots L_n^{p_n} T_y 
	\mid 0 \leq p_i \leq c_i^{\mu} -1 \, (1 \leq i \leq n), \, y \in \fS^{\mu}\} $ 
is an $R$-free basis of $M^{\mu}$. 
\end{enumerate}
\begin{proof}
(\roi). 
By \eqref{Tw Tx Ty} and \eqref{xmu Tw} together with Theorem \ref{AK}, 
we see that $\wt{M}^{\mu}$ is an $R$-free module with an $R$-free basis 
$\wt{\cB}^{\mu}_0 =\{ x_{\mu}T_y  L_1^{k_1} L_2^{k_2} \dots L_n^{k_n} 
	\mid 0 \leq k_i \leq r-1 \, (1\leq i \leq n), \, y \in \fS^{\mu}\} $. 
We also see that any element of $\wt{\cB}^{\mu}_0$ can be written as an $R$-linear combination 
of the elements of $\wt{\cB}^{\mu}$ by Lemma \ref{Li Lj} (\rovi). 
Thus $\wt{\cB}^{\mu}$ is an $R$-free basis of $\wt{M}^{\mu}$. 

(\roii). By Theorem \ref{AK},  we see that $M^{\mu}$ is spanned by 
$\cB_0^{\mu} = \{ m_{\mu} L_1^{k_1} L_2^{k_2} \dots L_n^{k_n} T_w \mid 0 \leq k_i \leq r-1, \, w \in \fS_n\}$ as an $R$-module. 
We claim that 
\begin{description}
\item[Claim A] 
$M^{\mu} $ is spanned by  
$\cB_1^{\mu} = \{ m_{\mu} L_1^{p_1} L_2^{p_2} \dots L_n^{p_n} T_w \mid 0 \leq p_i \leq c_i^{\mu} -1, \, w \in \fS_n\}$ 
as an $R$-module. 
\end{description}
Note that $c_i^{\mu} = c_{i+1}^{\mu}$ if $s_i \in \fS_{\mu}$ by \eqref{def cjmu}, 
then we see that $M^{\mu}$ is spanned by $\cB^{\mu}$ as an $R$-module from  
Claim A together with \eqref{Tw Tx Ty}, \eqref{xmu Tw} and Lemma \ref{Li Lj}. 

On the other hand, 
by \cite[Theorem 4.14]{DJM98}, 
we see that 
$M^{\mu}$ is an $R$-free module whose rank does not depend on a choice  of a ring $R$ and parameters. 
By considering the case where 
$\sH_{n,r}$ is isomorphic to the group ring of $\fS_n \ltimes (\ZZ/ r \ZZ)^n$, 
we see that $\rank_R M^{\mu} = | \cB^{\mu} |$ by \cite[Remarks 3.9. (\roii)]{DJM98}. 
As a consequence, 
we see that the set $\cB^{\mu}$ is an $R$-free basis of $M^{\mu}$.  
Therefore,  to complete a proof of (\roii), it is enough to show Claim A. 

For $m_{\mu} L_1^{k_1} L_2^{k_2} \dots L_n^{k_n} T_w \in \cB_0^{\mu}$,  
there exists an  integer $j \in \{ 0, 1, \dots,  n\}$  such that 
$k_i < c_i^{\mu}$ for all $i > j$ and $k_j \geq c_j^{\mu}$.  
If $j=0$, we see that $m_{\mu} L_1^{k_1} L_2^{k_2} \dots L_n^{k_n} T_w  \in \cB_1^{\mu}$.  
Suppose that $j >0$, and we have 
\begin{align*}
&m_{\mu} L_1^{k_1} L_2^{k_2} \dots L_n^{k_n} T_w 
\\
&= m_{\mu} L_1^{k_1} \dots L_{j-1}^{k_{j-1}} ( L_j^{c_j^{\mu}} L_j^{k_j - c_j^{\mu}}) L_{j+1}^{k_{j+1}} \dots L_n^{k_n} T_w 
\\
&= - \sum_{0 \leq p_1,p_2, \dots, p_j \leq c_j^{\mu} -1} m_{\mu} L_1^{k_1} \dots L_{j-1}^{k_{j-1}} 
	( L_1^{p_1} L_2^{p_2} \dots L_{j}^{p_j} h_{(p_1,p_2,\dots, p_j)} ) L_j^{k_j- c_j^{\mu}} L_{j+1}^{k_{j+1}} \dots L_n^{k_n} T_w 
\end{align*} 
by Lemma \ref{Li lan k ran} and Lemma \ref{Lemma mmu L1k1 Lj-1kj-1 Lj lan li ran}. 
Since $h_{(p_1,p_2,\dots, p_j)} \in \sH_{[1,j]}$, 
we also have 
\begin{align*}
h_{(p_1,p_2,\dots, p_j)} L_j^{k_j - c_j^{\mu}} 
= \sum_{ t_1,t_2,\dots, t_j \geq 0 \atop t_1+t_2+\dots+ t_j = k_j - c_j^{\mu}}  
	L_1^{t_1} L_2^{t_2} \dots L_j^{t_j} h_{(t_1,t_2, \dots, t_j)} 
\end{align*}
for some $h_{(t_1,t_2,\dots, t_j)} \in \sH_{[1,j]}$ by Lemma \ref{Li Lj}. 
Note that any element of $ \sH_{[1,j]}$ commute with $L_i $ for $ i >j$ by Lemma \ref{Li Lj} (\roii), 
and we see that 
\begin{align*}
m_{\mu} L_1^{k_1} L_2^{k_2} \dots L_n^{k_n} T_w 
= \sum_{ 0 \leq l_1, l_2,\dots, l_{j-1} } \sum_{0\leq l_j \leq k_j -1}
	 m_{\mu} L_1^{l_1} L_2^{l_2} \dots L_{j-1}^{l_{j-1}} L_j^{l_j} L_{j+1}^{k_{j+1}} \dots L_n^{k_n} 
	 h_{(l_1,\dots,l_j)} 
\end{align*}
for some $ h_{(l_1,\dots,l_j)}  \in \sH_{[1,n]}$. 
By repeating the above argument, we have Claim A. 
\end{proof}
\end{prop}

\para 
In order to describe the defining relations of the cyclic right $\sH_{n,r}$-modules $\wt{M}^{\mu}$ and $M^{\mu}$ 
($\mu \in \vL_{n,r}(\bm)$) respectively, 
we consider the following ideals of $\sH_{n,r}$. 

For $\mu \in \vL_{n, r}(\bm)$, 
let $\wt{\cI}^{\mu}$ be the right ideal of $\sH_{n,r}$ generated by 
\begin{align}
\label{Ti - q} 
\{T_i - q \mid s_i \in \fS_{\mu} \}, 
\end{align}  
and $\cI^{\mu}$ be the right ideal of $\sH_{n,r}$ generated by  
\begin{align}
\label{Ti - q Ljlan lj ran}
\{T_i - q \mid s_i \in \fS_{\mu}\} 
\cup 
\{ L_j^{\lan c_j^{\mu} \ran} \mid 1 \leq j \leq n \}. 
\end{align}
Then we have the following lemma corresponding to the relations \eqref{xmu Tw} 
and Lemma \ref{Lemma mmu L1k1 Lj-1kj-1 Lj lan li ran}. 

\begin{lem}
\label{rel Imu}
For $\mu \in \vL_{n,r}(\bm)$, we have the following. 
\begin{enumerate}
\item  
For $x \in \fS_{\mu}$, 
we have $T_x - q^{\ell (x)} \in \wt{\cI}^{\mu} \cap \cI^{\mu}$. 

\item  
For $1\leq j \leq n$ and $k_1,k_2,\dots, k_{n} \geq 0$,  
we have 
\begin{align*}
L_1^{k_1} L_2^{k_2} \dots L_{n}^{k_{n}} L_j^{\lan c_j^{\mu} \ran} \in \cI^{\mu}. 
\end{align*} 
\end{enumerate}

\begin{proof}
For $x  \in \fS_{\mu}$, 
let $x =s_{i_1} s_{i_2} \dots s_{i_l}$ be a reduced expression such that $s_{i_k} \in \fS_{\mu}$ ($1\leq k \leq l$). 
If $\ell (x) >1$,  we have  
\begin{align*}
T_x - q^{\ell (x)} 
&= T_{i_1} T_{i_2} \dots T_{i_l} - q^{ \ell (x)} 
= (T_{i_1} - q ) T_{i_2} T_{i_3} \dots T_{i_l} 
	+ q  ( T_{i_2} T_{i_3} \dots T_{i_l } - q^{\ell (x) -1} ), 
\end{align*}
and we can prove (\roi) by the induction on $\ell (x)$. 

We prove (\roii) by the induction on $k= k_1+k_2+\dots + k_{n}$.  
In the case where $k=0$, it is clear.  
Suppose that $k  >0$, 
and take $p$ such that $k_p \not=0$ and $k_i =0$ for all $ i >p$. 
Then we have 
\begin{align*}
&L_1^{k_1} L_2^{k_2} \dots L_{n}^{k_{n}} L_j^{\lan c_j^{\mu} \ran} 
\\
&= L_1^{k_1} L_2^{k_2} \dots L_{p-1}^{k_{p-1}} L_p^{k_p -1} 
	\\ & \quad \times 
	\begin{cases} 
	L_j^{\lan c_j^{\mu} \ran} L_p & \text{ if } p >j, 
	\\
	\dis 
	\big( L_j^{\lan c_j^{\mu} \ran} L_j   - (q-q^{-1}) \sum_{ z =1 }^{j-1} (L_j^{\lan c_j^{\mu} \ran} L_{j-z} - L_{j-z}^{\lan c_j^{\mu} \ran} L_j ) 
		(\TT_{j-z+1, j-1})^{-1} \TT_{j - z, j-1} \big) & \text{ if } p=j, 
	\\
	\big( L_j^{\lan c_j^{\mu} \ran} L_p  + (q-q^{-1}) (L_j^{\lan c_j^{\mu} \ran} L_p - L_p^{\lan c_j^{\mu} \ran} L_j ) 
		( \TT_{p+1, j-1})^{-1} \TT_{p, j-1}  \big) 
		& \text{ if } p < j 
	\end{cases} 
\end{align*}  
by Lemma \ref{Lemma Li Ljlan k ran}. 
Note that $ c_{z}^{\mu} \leq c_j^{\mu} $ if $z \leq j$ by the definition \eqref{def cjmu},  
and we see that the right-hand side of the above equation belong to the ideal $ \cI^{\mu}$ 
by the induction hypothesis together with \eqref{Llan l ran L lan j ran}. 
\end{proof}
\end{lem} 
In a similar way  as in the proof of Proposition \ref{Prop basis Mmu}, 
we have the following corollary by using Lemma \ref{rel Imu}. 


\begin{cor}
\label{Cor span Hnr Imu} 
For $\mu \in \vL_{n,r}(\bm)$, we have the following.  
\begin{enumerate}
\item 
$\sH_{n,r} / \wt{\cI}^{\mu}$ is spanned by 
\begin{align*}
\{ L_1^{k_1} L_2^{k_2} \dots L_n^{k_n} T_y + \wt{\cI}^{\mu} \mid 0 \leq k_i \leq r-1 \, (1\leq i \leq n), \, y \in \fS^{\mu}\} 
\end{align*}
as an $R$-module.  

\item 
$\sH_{n,r} / \cI^{\mu}$ is spanned by 
\begin{align*}
\{ L_1^{p_1} L_2^{p_2} \dots L_n^{p_n} T_y + \cI^{\mu} \mid 0 \leq p_i \leq c_i^{\mu} -1 \, (1\leq i \leq n), \, y \in \fS^{\mu} \} 
\end{align*} 
as an $R$-module. 
\end{enumerate}
\end{cor} 

Thanks to Proposition \ref{Prop basis Mmu} and Corollary \ref{Cor span Hnr Imu}, 
we obtain  defining relations of the cyclic right $\sH_{n,r}$-modules $\wt{M}^{\mu}$ and $M^{\mu}$ 
respectively as follows. 

\begin{thm}
\label{Thm def rel mmu}
For $\mu \in \vL_{n,r}(\bm)$, we have the following.  
\begin{enumerate}
\item 
We have $\wt{M}^{\mu} \cong \sH_{n,r} / \wt{\cI}^{\mu}$ as right $\sH_{n,r}$-modules. 
In particular, the set \eqref{Ti - q} gives a defining relation of the cyclic right $\sH_{n,r}$-module $\wt{M}^{\mu}$. 

\item 
We have $M^{\mu} \cong \sH_{n,r} / \cI^{\mu}$ as right $\sH_{n,r}$-modules. 
In particular, the set \eqref{Ti - q Ljlan lj ran} gives a defining relation of the cyclic right $\sH_{n,r}$-module $M^{\mu}$. 
\end{enumerate} 
\begin{proof}
We prove only (\roii). 
Since $M^{\mu}$  is the cyclic right $\sH_{n,r}$-module generated by $m_{\mu}$, 
there exists a natural surjective $\sH_{n,r}$-homomorphism 
$\vf : \sH_{n,r} \ra M^{\mu}$. 
By \eqref{xmu Tw} and  Lemma \ref{Lemma mmu L1k1 Lj-1kj-1 Lj lan li ran}, 
we see that $\Ker \vf$ include the right ideal $\cI^{\mu}$ of $\sH_{n,r}$, 
and the homomorphism $\vf$ induces the surjective $\sH_{n,r}$-homomorphism 
$\ol{\vf} : \sH_{n,r}/ \cI^{\mu} \ra M^{\mu}$. 
Then we conclude that the homomorphism $\ol{\vf}$ is an isomorphism 
by Proposition \ref{Prop basis Mmu} and Corollary \ref{Cor span Hnr Imu}. 
We can prove (\roi) in a similar way. 
\end{proof}
\end{thm}

\para 
Finally, we recall some generators of the cyclotomic $q$-Schur algebra $\sS_{n,r}(\bm)$ obtained in \cite{Wad11}. 

For $1 \leq i \leq m-1$, 
we define the elements $\sE_i, \sF_i \in \sS_{n,r}(\bm)$ by
\begin{align*}
&\sE_i (m_{\mu}) = \d_{(\mu_{i+1} \not=0)} q^{ - \mu_{i+1} +1} m_{\mu + \a_i} 
	\cdot \big( 1 + \sum_{p=1}^{\mu_{i+1} -1}  q^p \TT_{N_i^{\mu} +1, \, N_i^{\mu} +p} \big), 
\\
& \sF_i (m_{\mu}) = \d_{(\mu_i \not=0)} \begin{cases}
		\dis 
		q^{- \mu_i +1} m_{\mu - \a_i} \cdot \big( 1 + \sum_{p=1}^{\mu_i-1} q^p \wt{\TT}_{N_i^{\mu} -1, \, N_i^{\mu} -p} \big) 
			\hspace{-7em} 
			\\
			& \text{ unless } \xi^{-1}(i) = (m_k,k) \text{ for some } k, 
		\\ \dis 
		q^{- \mu_i+1} ( - Q_k^{-1}) m_{\mu - \a_i} \cdot ( L_{N_i^{\mu}} - Q_k) 
			\big( 1 + \sum_{p=1}^{\mu_i-1} q^p \wt{\TT}_{N_i^{\mu} -1, \, N_i^{\mu} -p} \big) 
			\hspace{-17em} 
			\\ 
			& \text{ if } \xi^{-1}(i) = (m_k,k) \text{ for some } k
	\end{cases}
\end{align*}
for $\mu \in \vL_{n,r}(\bm)$. 
We also define the element $1_{\mu} \in \sS_{n,r} (\bm) $  for $\mu \in \vL_{n,r}(\bm)$ as the projection to $M^{\mu}$. 
Then, we have the following proposition. 


\begin{prop}[{\cite[Proposition 7.7]{Wad11}}] 
\label{Prop gen sSnr}
Assume that $R$ is a filed, $q$ is not a root of unity and $m_i \geq n$ for all $i=1,2,\dots, r$, 
then the cyclotomic $q$-Schur algebra $\sS_{n,r}(\bm)$ is generated by 
$\sE_i$, $\sF_i$ ($1\leq i \leq m-1$) and $1_{\mu}$ ($\mu \in \vL_{n,r}(\bm) $). 
\end{prop}

\para 
\label{Para integral CSA}
We can also give  generators of $\sS_{n,r}(\bm)$ over an arbitrary ring $R$ with any parameters 
$q, Q_0, \dots, Q_{r-1} \in R^{\times}$ as follows. 
Let $\AA' =\ZZ [ \hat{Q}_0, \hat{Q}_1\dots, \hat{Q}_{r-1}]$ be the polynomial ring over $\ZZ$ with 
indeterminate variables $\hat{Q}_0, \hat{Q}_1,\dots, \hat{Q}_{r-1}$, 
and 
$\AA = \AA'[ \hat{q}, \hat{q}^{-1}] =\ZZ [ \hat{q}, \hat{q}^{-1}, \hat{Q}_0, \hat{Q}_1\dots, \hat{Q}_{r-1}] $ 
be the Laurent polynomial ring over $\AA'$ with the invertible variable $\hat{q}$. 
We also denote by $\KK = \QQ (\hat{q}, \hat{Q}_0, \dots, \hat{Q}_{r-1})$ the quotient field of $\AA$. 

We denote by $\sH_{n,r}^{\KK}$ (resp. $\sH_{n,r}^{\AA}$) and $\sS_{n,r}^{\KK}(\bm)$ (resp. $\sS_{n,r}^{\AA}(\bm)$)  
the Ariki-Koike algebra and the cyclotomic $q$-Schur algebra over $\KK$ (resp. over $\AA)$ with parameters 
$\hat{q}, \hat{Q}_0, \dots, \hat{Q}_{r-1}$.  
We also denote by $\sH_{n,r}^R$ and $\sS_{n,r}^R (\bm)$ the algebras over an arbitrary ring $R$ with parameters 
$q, Q_0, \dots, Q_{r-1} \in R^{\times}$. 
Then, the algebra $\sH_{n,r}^{\AA}$ (resp. $\sS_{n,r}^{\AA}$) is an $\AA$-subalgebra of $\sH_{n,r}^{\KK}$ (resp. $\sS_{n,r}^{\KK}$), 
and the algebra $\sH_{n,r}^R$ (resp. $\sS_{n,r}^{R} (\bm)$) is obtained as a specialized algebra $R \otimes_{\AA} \sH_{n,r}^{\AA}$ 
(resp. $ R \otimes_{\AA} \sS_{n,r}^{\AA}(\bm)$) 
through the ring homomorphism 
$\AA \ra R$ such that the parameters $\hat{q}, \hat{Q}_0, \dots, \hat{Q}_{r-1}$ send to 
$q, Q_0, \dots,  Q_{r-1}$ respectively. 

For $1\leq i \leq m-1$ and a positive integer $d \in \ZZ_{>0}$, 
put 
\begin{align*}
\sE_i^{(d)} = \frac{\sE_i^d }{ [d]_{\wh{q}} !}, 
\quad 
\sF_i^{(d)} = \frac{\sF_i^d }{ [d]_{\wh{q}} ! } 
\quad \in \sS_{n,r}^{\KK}(\bm), 
\end{align*} 
where $[d]_{\hat{q}} = (\hat{q}^d - \hat{q}^{-d}) / (\hat{q}- \hat{q}^{-1})$ and $[d]_{\hat{q}} ! = [d]_{\hat{q}} [d-1]_{\hat{q}} \dots [1]_{\hat{q}}$. 
Then the elements $\sE_i^{(d)}$ and $\sF_i^{(d)}$ belong to the $\AA$-subalgebra $\sS_{n,r}^{\AA} (\bm)$ 
of $\sS_{n,r}^{\KK}(\bm)$ 
(see also the proof of \cite[Theorem 8.1]{Wad16}). 
We also denote  the elements $1 \otimes \sE_i^{(d)}$ and $1 \otimes \sF_i^{(d)}$ 
of $\sS_{n,r}^R (\bm) = R \otimes_{\AA} \sS_{n,r}^{\AA} (\bm)$ 
by $\sE_i^{(d)}$ and $\sF_i^{(d)}$ simply. 
Then we have the following proposition. 

\begin{prop}[{\cite[Proposition 7.7]{Wad11}}] 
\label{Prop gen CSA AA}
Assume that  $m_i \geq n$ for all $i=1,2,\dots, r$, 
then the cyclotomic $q$-Schur algebra $\sS_{n,r}^R (\bm)$ is generated by 
$\sE_i^{(d)}$, $\sF_i^{(d)}$ ($1\leq i \leq m-1$, $d \in \ZZ_{>0}$) and $1_{\mu}$ ($\mu \in \vL_{n,r}(\bm) $). 
\end{prop}


\section{Quantum Schur-Weyl duality}  
\label{Section quantum SW}

In this section, we recall the quantum Schur-Weyl duality given in \cite{Jim86}, 
and also recall the connection with the classical $q$-Schur algebras. 

\para 
The quantum group $U_q (\Fgl_m)$ associated with the general linear Lie algebra $\Fgl_m$ 
is an associative algebra over $\CC$ with a parameter $q \in \CC \setminus \{0, \pm 1\}$  
generated by $E_i, F_i$ ($1 \leq i \leq m-1$) and $K_j^{\pm}$ ($1\leq j \leq m$) 
with the following defining relations:  
\begin{align*}
&K_i^+ K_i^- = K_i^- K_i^+ =1, 
\quad 
K_i^+ K_j^+ = K_j^+ K_i^+, 
\\ 
& K_i^+ E_j K_i^- = q^{\d_{i,j} - \d_{i,j+1}} E_j, 
\quad 
K_i^+ F_j K_i^- = q^{- (\d_{i,j} - \d_{i,j+1})} F_j,  
\\
& E_i F_j - F_j E_i = \d_{i,j} \frac{K_i^+ K_{i+1}^- - K_i^- K_{i+1}^+}{ q - q^{-1}},  
\\
&
E_{i \pm 1} E_i^2 - ( q + q^{-1} ) E_i E_{i \pm 1} E_i + E_i^2 E_{i \pm 1} =0, 
\quad 
E_i  E_j = E_j E_i  \text{ if } j \not=i \pm 1, 
\\
&
F_{i \pm 1} F_i^2 - ( q + q^{-1} ) F_i F_{i \pm 1} F_i + F_i^2 F_{i \pm 1} =0, 
\quad 
F_i  F_j = F_j F_i  \text{ if } j \not=i \pm 1. 
\end{align*}

The quantum group $U_q (\Fgl_m)$ has the coproduct $\D : U_q (\Fgl_m) \ra U_q (\Fgl_m) \otimes U_q (\Fgl_m)$ 
defined by 
\begin{align*}
&\D (E_i) = E_i \otimes K_i^+ K_{i+1}^- + 1 \otimes E_i, 
\quad 
\D (F_i) = F_i \otimes 1 + K_i^- K_{i+1}^+ \otimes F_i, 
\\
&\D( K_j^{\pm}) = K_j^{\pm} \otimes K_j^{\pm}.
\end{align*}

\para 
Let $V$ be an $m$-dimensional vector space over $\CC$ with a basis $\{v_1,v_2, \dots, v_m\}$. 
We define the left action of $U_q(\Fgl_m)$ on $V$ by 
\begin{align*}
E_i \cdot v_j = \d_{j, i+1} v_{j-1}, 
\quad 
F_i \cdot v_j = \d_{j, i} v_{j+1}, 
\quad 
K_i^{\pm} \cdot v_j = q^{\pm \d_{j,i}} v_j.
\end{align*}
The $U_q(\Fgl_m)$-module $V$ is called the $q$-vector representation of $U_q(\Fgl_m)$. 
The quantum group $U_q(\Fgl_m)$ acts on the tensor space $V^{\otimes n}$ through the coproduct $\D$. 
We denote by $\rho : U_q (\Fgl_m) \ra \End_{\CC} (V^{\otimes n})$ the representation corresponding to this action. 

On the other hand, 
we define a linear transformation $\cT \in \End_{\CC} (V \otimes V)$ by  
\begin{align*}
\cT (v_i \otimes v_j) 
= \begin{cases} 
	q v_{i} \otimes v_j & \text{ if } i=j, 
	\\
	v_j \otimes v_i & \text{ if } i <j, 
	\\
	v_j \otimes v_i + ( q -  q ^{-1}) v_i \otimes v_j & \text{ if } i >j, 
	\end{cases} 
\end{align*}
and we define the right action of the Iwahori-Hecke algebra $\sH_n$ associated with $\fS_n$ on $V^{\otimes n}$ by 
\begin{align*}
T_i = (\id_V)^{\otimes i-1} \otimes \cT \otimes (\id_V)^{\otimes (n-i-1)}
\end{align*} 
for $i=1,2,\dots, n-1$. 
We denote the representation corresponding to this action  by $\s : \sH_{n}^{\opp} \ra \End_{\CC} (V^{\otimes n})$. 
Then we have the quantum Schur-Weyl duality obtained in \cite{Jim86}. 


\begin{thm}[{\cite[Proposition 3]{Jim86}}] 
The actions of $U_{q}(\Fgl_m)$ and of $\sH_n$ on $V^{\otimes n}$ commute with each other. 
If $q$ is not a root of unity, 
 we have 
\begin{align*}
\Im \rho = \End_{\sH_n^{\opp}}( V^{\otimes n}) 
\text{ and } 
\Im \s = \End_{U_{q}(\Fgl_m)} (V^{\otimes n}). 
\end{align*}
Moreover, if $m \geq n$, the homomorphism $\s$ is injective, 
and we have the isomorphism $\sH_n \cong \End_{U_{q}(\Fgl_m)} (V^{\otimes n})$. 
\end{thm}


\para 
Suppose that $q$ is not a root of unity. 
We have the following weight space decomposition of $V^{\otimes n}$ as the  $U_{q} (\Fgl_m)$-module: 
\begin{align*}
V^{\otimes n} = \bigoplus_{\mu \in \vL_n(m)} V^{\otimes n}_{\mu}, 
\quad 
V^{\otimes n}_{\mu} =\{ v \in V^{\otimes n} \mid K_j^+ \cdot v = q^{\mu_j} v \,\, (1\leq j \leq m) \}. 
\end{align*}

Since the action of $\sH_n$ on $V^{\otimes n}$ commute with the action of $U_{q}(\Fgl_m)$,  
the weight space $V^{\otimes n}_{\mu}$ is an $\sH_{n}$-submodule of $V^{\otimes n}$.  
Moreover, from definitions, we see that 
\begin{align*}
\{ v_{j_1} \otimes v_{j_2} \otimes \dots \otimes v_{j_n} \mid \sharp \{ k \mid j_k = i \} = \mu_i \, ,1 \leq i \leq m\} 
\end{align*}
gives a $\CC$-basis of $V^{\otimes n}_{\mu}$, 
and  the right $\sH_{n}$-module $V_{\mu}^{\otimes n} $ is generated by 
\begin{align}
\label{def vmu}
v_{\mu} = \underbrace{v_{1} \otimes \dots \otimes v_1}_{\mu_1} \otimes \underbrace{v_2 \otimes \dots \otimes v_2}_{\mu_2} 
\otimes \dots \otimes \underbrace{v_m \otimes \dots \otimes v_m}_{\mu_m}. 
\end{align} 

We can easily check that $v_{\mu} T_i = q v_{\mu}$ if $s_i \in \fS_{\mu}$. 
Then, we have the surjective $\sH_n$-homomorphism 
$\vf: M^{\mu} \ra V^{\otimes n}_\mu$ such that $x_{\mu} \mapsto v_{\mu}$ by Theorem \ref{Thm def rel mmu}, 
where we note that $m_{\mu} = x_{\mu}$, thus $M^{\mu} = \wt{M}^{\mu}$, in the case of the Iwahori-Hecke algebra 
$\sH_n$ associated with $\fS_n$. 
We can also check that 
\begin{align}
\label{sharp smu}
\sharp \fS^{\mu} = 
\sharp \{ v_{j_1} \otimes v_{j_2} \otimes \dots \otimes v_{j_n} \mid \sharp \{ k \mid j_k = i \} = \mu_i \, ,1 \leq i \leq m\}  
\end{align}
(see, e.g.  \cite[\S9.1]{DDPW08}), 
and we see that the homomorphism $\vf: M^{\mu} \ra V^{\otimes n}_\mu$ is an isomorphism. 
As a consequence, 
we have the algebra isomorphism 
\begin{align*}
\Im \rho = \End_{\sH_n^{\opp}} (V^{\otimes n}) \cong \End_{\sH_n^{\opp}} \big( \bigoplus_{\mu \in \vL_n(m)} M^{\mu} \big) = \sS_n (\bm)
\end{align*}
which is a quotient of $U_{q} (\Fgl_m)$. 


\remark  
The results in this section also hold over an any field $F$ and any parameter $q \in F^{\times}$ 
by replacing $U_\bq(\Fgl_m)$ with the restricted quantum group associated  with $\Fgl_m$ over $F$ 
which is the specialized algebra of the Lusztig's integral form using divided powers. 
This fact is pointed out in \cite{Du95} by using the geometric realization of $U_q(\Fgl_m)$ given in \cite{BLM90}.



\section{Quantum affine Schur-Weyl duality} 
\label{Section QASW}

In this section, 
we recall the quantum affine Schur-Weyl duality given in 
\cite{GRV94} and \cite{CP96} independently. 

\para 
The quantum affine algebra $U_q (\wh{\Fsl}_m)$ is an associative algebra over $\CC$ with a parameter 
$q \in \CC \setminus \{0, \pm 1\}$ 
generated by $e_i, f_i, k_i^{\pm}$ ($0 \leq i \leq m-1$) with the following defining relations:  
\begin{align*}
&k_i^+ k_i^- = k_i^- k_i^+=1, 
\quad 
k_i^+ k_j^+ = k_j^+ k_i^+, 
\quad 
k_i^+ e_j k_i^- = q^{a_{ij}} e_j, 
\quad 
k_i^+ f_j k_i^- = q^{-a_{ij}} f_j, 
\\
& e_i f_j - f_j e_i = \d_{i,j} \frac{k_i^+ - k_i^-}{q-q^{-1}}, 
\\
& e_j e_i^2 - (q+q^{-1}) e_i e_j e_i + e_i^2 e_j =0 \text{ if } a_{ij} = -1, 
\quad 
e_i e_j = e_j e_i  \text{ if } a_{ij} =0,   
\\
& f_j e_i^2 - (q+q^{-1}) f_i f_j f_i + f_i^2 f_j =0 \text{ if } a_{ij} = -1, 
\quad 
f_i f_j = f_j f_i  \text{ if } a_{ij} =0,   
\end{align*}
where $(a_{ij})_{0 \leq i,j \leq m-1}$ is the Cartan matrix of type $A_{m-1}^{(1)}$. 

The quantum affine algebra $U_q( \wh{\Fsl}_m)$ has the coproduct 
$\D : U_q (\wh{\Fsl}_m) \ra U_q (\wh{\Fsl}_m)  \otimes U_q (\wh{\Fsl}_m) $ defined by 
\begin{align*}
\D (e_i) = e_i \otimes k_i^+ + 1 \otimes e_i, 
\quad 
\D(f_i) = f_i \otimes 1 +  k_i^- \otimes f_i, 
\quad 
\D(k_i^{\pm} ) = k_i^{\pm} \otimes k_i^{\pm}.
\end{align*}
We call it Drinfeld-Jimbo coproduct of $U_q (\wh{\Fsl}_m) $. 

\para 
For $\g \in \CC^{\times}$,  
there exists an algebra homomorphism 
$\ev_{\g} : U_q(\wh{\Fsl}_m) \ra U_q (\Fgl_m)$, 
so-called an evaluation homomorphism at $\g$,  
such that 
\begin{align*}
&e_i \mapsto E_i, \quad f_i \mapsto F_i, \quad k_i^+ \mapsto K_i^+ K_{i+1}^- \quad (1\leq i \leq m-1), 
\\
& k_0^+ \mapsto K_1^- K_m^+,  
\\
& e_0 \mapsto \g q^{-1} (K_1^+ K_m^+) [ F_{m-1}, [F_{m-2}, \dots, [F_2, F_1]_{q^{-1}} \dots ]_{q^{-1}} ]_{q^{-1}}, 
\\
& f_0 \mapsto (-1)^m \g^{-1} q^{m-1} (K_1^- K_m^-) [E_{m-1}, [E_{m-2}, \dots, [E_2, E_1]_{q^{-1}} \dots ]_{q^{-1}}]_{q^{-1}}.
\end{align*}

For a $U_q (\Fgl_m)$-module $M$, 
we regard $M$ as a $U_q ( \wh{\Fsl}_m)$-module through the homomorphism $\ev_{\g}$, 
and we denote it by $M^{\ev_{\g}}$. 

Recall that $V$ is the vector representation of $U_q(\Fgl_m)$ considered in the previous section. 
Let $\wh{V} = V^{\ev_1} \otimes \CC [x,x^{-1}]$ be the affinization of $V^{\ev_1}$ in the sense of \cite[\S4.2]{Kas02}, 
namely the generators $e_i $, $f_i $ and $k_i^{\pm}$  ($0\leq i \leq m-1$) act on  $\wh{V} $ 
by $e_i \otimes x^{\d_{i,0}}$, $f_i \otimes x^{- \d_{i,0}}$ and $k_i^{\pm} \otimes 1$ respectively. 
Then, we have that 
\begin{align}
\label{action whsl on whV}
\begin{split} 
&e_0 \cdot (v_j \otimes x^k) = \d_{j,1} v_m \otimes x^{k+1},  
\quad 
f_0 \cdot (v_j \otimes x^k ) = \d_{j,m} v_1 \otimes x^{k-1},  
\\
&k_0^{\pm} \cdot (v_j \otimes x^k) = q^{\pm (\d_{m,j} - \d_{1,j})} v_j \otimes x^k, 
\\
& e_i \cdot (v_j \otimes x^k)= \d_{i+1,j} v_{j-1} \otimes x^k, 
\quad 
f_i \cdot (v_j \otimes x^k) = \d_{i,j} v_{j+1} \otimes x^k, 
\\
& k_i^{\pm} \cdot ( v_j \otimes x^k) = q^{\pm (\d_{i,j} - \d_{i+1,j})} v_j \otimes x^k   
\quad (1\leq i \leq m-1).
\end{split}
\end{align}
The quantum affine algebra $U_q (\wh{\Fsl}_m)$ acts on the tensor space $\wh{V}^{\otimes n}$ through the coproduct $\D$. 

On the other hand, 
we define the right action of $\wh{\sH}_n$ on $\wh{V}^{\otimes n}$ as follows. 
We identify the $\CC$-vector space $\wh{V}^{\otimes n} = (V \otimes \CC[x,x^{-1}])^{\otimes n}$ 
with $V^{\otimes n} \otimes \CC [x_1^{\pm}, x_2^{\pm}, \dots, x_n^{\pm}]$ 
by the natural way. 
Then, we define a linear transformation $\wh{\cT} \in \End_{\CC} (\wh{V} \otimes \wh{V})$ by 
\begin{align*}
&\wh{\cT} ( (v_{j_1} \otimes v_{j_2}) \otimes x_1^{k_1} x_2^{k_2} ) 
\\
&= \cT (v_{j_1} \otimes v_{j_2}) \otimes x_1^{k_2} x_2^{k_1} 
	+ (q-q^{-1}) (v_{j_1} \otimes v_{j_2}) \otimes \frac{x_1^{k_2} x_2^{k_1} - x_1^{k_1} x_2^{k_2}}{x_1 x_2^{-1} -1}, 
\end{align*}
and  we define the right action of $\wh{\sH}_n$ on $\wh{V}^{\otimes n}$ by 
\begin{align*}
T_i = (\id_{\wh{V}} )^{\otimes i-1} \otimes \wh{\cT} \otimes (\id_{\wh{V}})^{\otimes (n-i-1)} 
\quad (1\leq i \leq n-1), 
\end{align*}
and the action of $X_j^{\pm} $ ($1\leq j \leq n$) is given as the multiplication by $x_j^{\pm}$. 
Then, we have the following quantum affine Schur-Weyl duality. 


\begin{thm}[{ \cite[Theorem 6.8]{GRV94}, \cite[Theorem 4.2]{CP96} }]
\label{Thm quantum affine SW}
The actions of $U_q (\wh{\Fsl}_m)$ and of $\wh{\sH}_n$ on $\wh{V}^{\otimes n}$ commute with each other.  
If $q$ is not a root of unity and $m \geq n$, then the functor 
\begin{align*}
\wh{V}^{\otimes n} \otimes_{\wh{\sH}_n} - : \wh{\sH}_n \cmod \ra U_q (\wh{\Fsl}_m) \cmod, 
\quad 
M \mapsto \wh{V}^{\otimes n} \otimes_{\wh{\sH}_n} M
\end{align*}
gives an equivalence between the category of finite dimensional $\wh{\sH}_n$-modules 
and the corresponding subcategory of $U_q (\wh{\Fsl}_m) \cmod$.  
\end{thm}


\section{The $(U_q(\wh{\Fsl}_m), \sH_{n,r})$-bimodule $\wt{V}_{n,r} = \wh{V}^{\otimes n} \otimes_{\wh{\sH}_n} \sH_{n,r}$} 
\label{Section wtVnr}

\para 
In this section, we denote $V^{\ev_1}$ by $V$ simply, and we identify the $\CC$-vector space $\wh{V}^{\otimes n} = (V \otimes \CC[x^{\pm}])^{\otimes n}$ 
with $V^{\otimes n} \otimes \CC [x_1^{\pm}, x_2^{\pm}, \dots, x_n^{\pm}] $ by the natural way. 

Put 
\begin{align*}
\wt{V}_{n,r} = \wh{V}^{\otimes n} \otimes_{\wh{\sH}_{n}} \sH_{n,r}, 
\end{align*}
where we regard $\sH_{n,r}$ as the left $\wh{\sH}_n$-module through the surjection \eqref{surj whH H}. 
Then the space $\wt{V}_{n,r}$ has a $(U_q(\wh{\Fsl}_m), \sH_{n,r})$-bimodule structure, 
and the functor 
\begin{align*}
\wt{V}_{n,r} \otimes_{\sH_{n,r}} - : \sH_{n,r} \cmod \ra U_q(\wh{\Fsl}_m) \cmod, 
\quad 
M \mapsto \wt{V}_{n,r} \otimes_{\sH_{n,r}} M
\end{align*} 
is the restriction of the functor $\wh{V}^{\otimes n} \otimes_{\wh{\sH}_{n}} -$ in Theorem \ref{Thm quantum affine SW} 
to the subcategory $\sH_{n,r} \cmod \subset \wh{\sH}_n \cmod$. 
Thus, we have the following corollary. 
\begin{cor}
If $q$ is not a root of unity and $m \geq n$, 
then the functor $\wt{V}_{n,r} \otimes_{\sH_{n,r}} -$ gives an equivalence 
 between the category of finite dimensional $\sH_{n,r}$-modules 
and the corresponding subcategory of $U_q (\wh{\Fsl}_m) \cmod$.  
\end{cor}

\para 
From the definitions, we have that 
\begin{align*}
\left\{ (v_{j_1}\otimes v_{j_2} \otimes \dots \otimes v_{j_n}) \otimes x_1^{k_1} x_2^{k_2} \dots x_n^{k_n} \otimes 1
	\mid \begin{matrix} 1 \leq j_1, \dots, j_n \leq m, \\ 0 \leq k_1,\dots, k_n \leq r-1 \end{matrix} \right\} 
\end{align*}
gives the $\CC$-basis of $\wt{V}_{n,r}$ (see \cite[Lemma 5.5]{LR06}). 
By \eqref{action whsl on whV}, we see that 
$\wt{V}_{n,r}$ has the weight space decomposition  
\begin{align}
\label{decom wtVnr}
\wt{V}_{n,r} = \bigoplus_{\mu \in \vL_n(m)} \wt{V}_{n,r}^{\mu}, 
\quad 
\wt{V}_{n,r}^{\mu} = \{ v \in \wt{V}_{n,r} \mid  k_i^+ \cdot v = q^{\mu_i - \mu_{i+1}} v \, \, (0 \leq i \leq m-1) \}, 
\end{align}
where we put $\mu_{0} := \mu_m$,  
and $\wt{V}_{n,r}^{\mu}$ has a $\CC$-basis 
\begin{align}
\label{basis wtVnrmu}
\left\{ 
(v_{j_1}\otimes v_{j_2} \otimes \dots \otimes v_{j_n}) \otimes x_1^{k_1} x_2^{k_2} \dots x_n^{k_n} \otimes 1
	\mid \begin{matrix} \sharp \{ k \mid j_k =i\} = \mu_i \,\, (1 \leq i \leq m)  
		\\
		0 \leq k_1, \dots, k_n \leq r-1
		\end{matrix}
\right\}. 
\end{align}

On the other hand, the decomposition \eqref{decom wtVnr} is also a decomposition as $\sH_{n,r}$-modules 
since the actions of $U_q(\wh{\Fsl}_m)$ and of $\sH_{n,r}$ on $\wt{V}_{n,r}$ commute with each other. 

For $\mu \in \vL_n(m)$, recall the element $v_{\mu} \in V^{\otimes n}$ given in \eqref{def vmu}, and put 
\begin{align*}
\wt{v}_{\mu} 
:= v_{\mu} \otimes 1 \otimes 1 \in  V^{\otimes n} \otimes \CC[x_1^{\pm}, \dots, x_n^{\pm}] \otimes_{\wh{\sH}_{n,r}} \sH_{n,r} 
= \wt{V}_{n,r}. 
\end{align*}
From the definition of the action of $\wh{\sH}_{n,r}$ on $\wh{V}^{\otimes n}$, 
we see that  $\wt{V}_{n,r}^{\mu}$ is generated by $\wt{v}_{\mu}$ as the $\sH_{n,r}$-module, 
and we also have 
\begin{align}
\label{wtvmu Ti}
\wt{v}_{\mu} \cdot T_i = q \wt{v}_{\mu} \text{ if } s_i \in \fS_{\mu}. 
\end{align}
Then, by Theorem \ref{Thm def rel mmu} (\roi), 
there exists a surjective $\sH_{n,r}$-homomorphism 
\begin{align}
\label{wtMmu to wtVnrmu}
\wt{M}^{\mu} \ra \wt{V}_{n,r}^{\mu}, \quad x_{\mu}  \mapsto \wt{v}_{\mu}, 
\end{align}
where we identify the set $\vL_{n,r}(\bm)$ with the set $\vL_{n}(m)$ by \eqref{bij vLnrbm vLnm}.  
Moreover, we see that 
$\dim \wt{M}^{\mu} = \dim \wt{V}_{n,r}^{\mu}$ by Proposition \ref{Prop basis Mmu}, \eqref{sharp smu} and \eqref{basis wtVnrmu}.
As a consequence, we have the following lemma. 
\begin{lem}
\label{Lemma iso wtMmu wt Vnrmu} 
For each $\mu \in \vL_{n,r}(\bm)$, there exists an isomorphism of $\sH_{n,r}$-modules 
from $\wt{M}^{\mu}$ to $\wt{V}_{n,r}^{\mu}$ 
given by \eqref{wtMmu to wtVnrmu}. 
\end{lem}

\para 
Let $\cI_{\wt{V}_{n,r}}$ be the $\sH_{n,r}$-submodule of $\wt{V}_{n,r}$ generated by 
\begin{align}
\label{gen wtVnr}
\{ \wt{v}_{\mu} \cdot L_j^{\lan c_j^{\mu} \ran} \mid \mu \in \vL_n(m), \, 1 \leq j \leq n \}. 
\end{align}
We remark that the $\sH_{n,r}$-submodule $\cI_{\wt{V}_{n,r}}$ has the decomposition 
\begin{align*}
\cI_{\wt{V}_{n,r}} = \bigoplus_{\mu \in \vL_n(m)} \cI_{\wt{V}_{n,r}} \cap \wt{V}_{n,r}^{\mu}, 
\end{align*}
where $\cI_{\wt{V}_{n,r}} \cap \wt{V}_{n,r}^{\mu}$ is generated by 
$\{\wt{v}_{\mu} \cdot L_j^{\lan c_j^{\mu} \ran} \mid 1\leq j \leq n\}$ as the $\sH_{n,r}$-module. 
Then, by Theorem \ref{Thm def rel mmu} and Lemma \ref{Lemma iso wtMmu wt Vnrmu},
 we see that 
there exists an $\sH_{n,r}$-isomorphism 
\begin{align*}
&\wt{V}_{n,r} / \cI_{\wt{V}_{n,r}} = \bigoplus_{\mu \in \vL_n(m)} \wt{V}_{n,r}^{\mu} / (\cI_{\wt{V}_{n,r}} \cap \wt{V}_{n,r}^{\mu}) 
\ra \bigoplus_{\mu \in \vL_{n,r}(\bm)} M^{\mu}, 
\\
&\wt{v}_{\mu} + \cI_{\wt{V}_{n,r}} \mapsto m_{\mu} \,\, (\mu \in \vL_n(m)). 
\end{align*}
However, the action of $U_q(\wh{\Fsl}_m)$ on $\wt{V}_{n,r}$ does not preserve the subspace $\cI_{\wt{V}_{n,r}}$, 
and the $U_q (\wh{\Fsl}_m)$-action on $\wt{V}_{n,r}$ does not induce the action on 
$\wt{V}_{n,r} / \cI_{\wt{V}_{n,r}}$. 
Hence, we need some shifts of $U_q (\wh{\Fsl}_m)$  and of $\wt{V}_{n,r}$ 
which are discussed in the subsequence sections. 

\para 
In the remaining this section, 
we see that the action of $U_q(\wh{\Fsl}_m)$ on $\wt{V}_{n,r}$ does not preserve the subspace $\cI_{\wt{V}_{n,r}}$. 
By  direct calculation (cf. \cite[Theorem 4.11]{GRV94}, \cite[Appendix A]{Wad11}), 
we have 
\begin{align}
\label{act Uqwhsl on wtVnr} 
\begin{split}
&k_i^{\pm} \cdot \wt{v}_{\mu} = q^{\pm (\mu_i - \mu_{i+1})} \wt{v}_{\mu} 
\\
& e_i \cdot \wt{v}_{\mu} = \d_{(\mu_{i+1} \not=0)} q^{- \mu_{i+1} +1} \wt{v}_{\mu + \a_i} 
	\cdot \big(1 + \sum_{p=1}^{\mu_{i+1} -1} q^p \TT_{N_i^{\mu}+1, \, N_i^{\mu} +p} \big), 
\\
& f_i \cdot \wt{v}_{\mu} = \d_{(\mu_i \not=0)} q^{- \mu_i +1} \wt{v}_{\mu - \a_i} \cdot 
	\big( 1 + \sum_{p=1}^{\mu_i-1} q^p \wt{\TT}_{N_i^{\mu} -1, N_i^{\mu} - p} \big), 
\\
\end{split}
\end{align}
for $1\leq i \leq m-1$. 

Suppose that $\xi^{-1}(i) = (m_k,k)$ for some $k$. 
Note that the actions of $U_q(\wh{\Fsl}_m)$ and of $\sH_{n,r}$ on $\wt{V}_{n,r}$ commute with each other, 
we have 
\begin{align*}
&f_i \cdot ( \wt{v}_{\mu} \cdot L_{N_i^{\mu}}^{\lan c_{N_i^{\mu}}^{\mu} \ran} )
\\
&=\d_{(\mu_i \not=0)} q^{- \mu_i +1} \wt{v}_{\mu - \a_i} \cdot 
	\big( 1 + \sum_{p=1}^{\mu_i-1} q^p \wt{\TT}_{N_i^{\mu} -1, N_i^{\mu} - p} \big)
	L_{N_i^{\mu}}^{\lan c_{N_i^{\mu}}^{\mu} \ran}
\\
&= \d_{(\mu_i \not=0)} q^{- \mu_i +1} \wt{v}_{\mu - \a_i} \cdot 
	\Big\{ L_{N_i^{\mu}}^{\lan c_{N_i^{\mu}}^{\mu} \ran} 
		\\ & \qquad 
		+ \sum_{p=1}^{\mu_i-1} q^p 
		\Big( L_{N_i^{\mu}-1}^{\lan c_{N_i^{\mu} }^{\mu} \ran} \wt{\TT}_{N_i^{\mu} -1, N_i^{\mu} - p} 
			+(q-q^{-1}) L_{N_i^{\mu}}^{\lan c_{N_i^{\mu}}^{\mu} \ran} (\TT_{N_i^{\mu},  N_i^{\mu}-1})^{-1} \wt{\TT}_{N_i^{\mu}-2, N_i^{\mu}-p} \Big)\Big\}
\\
&= \d_{(\mu_i \not=0)} q^{- \mu_i +1} \wt{v}_{\mu - \a_i} \cdot \Big\{ 
	L_{N_i^{\mu}}^{\lan c_{N_i^{\mu}}^{\mu} \ran} \Big( 1 + (q-q^{-1}) \sum_{p=1}^{\mu_i-1} q^p (\TT_{N_i^{\mu}, N_i^{\mu}-1})^{-1} \wt{\TT}_{N_i^{\mu}-2, N_i^{\mu}-p} \Big)
	\\ & \qquad 
	+ L_{N_i^{\mu}-1}^{\lan c_{N_i^{\mu} }^{\mu} \ran} \sum_{p=1}^{\mu_i-1} q^p   \wt{\TT}_{N_i^{\mu} -1, N_i^{\mu} - p} \Big\} 
\end{align*}
by \eqref{act Uqwhsl on wtVnr} and Lemma \ref{Lemma TT Lk}. 
Note that $\TT_{N_i^{\mu}, N_i^{\mu}-1} = 1$, and  that 
$L_{N_i^{\mu}}^{\lan c_{N_i^{\mu}} \ran}$ commutes with $\wt{\TT}_{N_i^{\mu}-2, N_i^{\mu}-p}$, 
by Lemma \ref{Lemma TT Lk} (\rov). 
Then, we have 
\begin{align*}
f_i \cdot ( \wt{v}_{\mu} \cdot L_{N_i^{\mu}}^{\lan c_{N_i^{\mu}}^{\mu} \ran} )
= \d_{(\mu_i \not=0)} q^{- \mu_i +1} \wt{v}_{\mu - \a_i} \cdot \big( 
	q^{2 \mu_i -2} L_{N_i^{\mu}}^{\lan c_{N_i^{\mu}}^{\mu} \ran} 
	+ L_{N_i^{\mu}-1}^{\lan c_{N_i^{\mu} }^{\mu} \ran} \sum_{p=1}^{\mu_i-1} q^p   \wt{\TT}_{N_i^{\mu} -1, N_i^{\mu} - p} \big), 
\end{align*}
where we use that $\wt{v}_{\mu-\a_i} \cdot \wt{\TT}_{N_i^{\mu}-2, \, N_i^{\mu}-p} = q^{p-1}$ 
($1\leq p \leq \mu_i-1$) by 
\eqref{wtvmu Ti}.
On the other hand, we see that 
\begin{align*}
c_{N_i^{\mu}}^{\mu - \a_i} = c_{N_i^{\mu}}^{\mu} +1 
\text{ and }
c_{N_i^{\mu}-1}^{\mu - \a_i} = c_{N_i^{\mu}}^{\mu} 
\text{ if } \mu_i \not=0 
\end{align*} 
by the definition \eqref{def cjmu} 
since $N_i^{\mu} = a_k^{\mu}$ and $a_k^{\mu -\a_i} = a_k^{\mu} -1$ in the case where  $\xi^{-1}(i) = (m_k,k)$. 
As a consequence, 
we have the following proposition. 

\begin{prop}
\label{Prop fi wtvmu L notin IwtV}
For $1\leq i \leq n-1$ and $\mu \in \vL_n(m)$, we have 
\begin{align*}
f_i \cdot ( \wt{v}_{\mu} \cdot L_{N_i^{\mu}}^{\lan c_{N_i^{\mu}}^{\mu} \ran} ) 
\not\in \cI_{\wt{V}_{n,r}} 
\text{ if } \mu_i \not=0 \text{ and } \xi^{-1}(i) =(m_k,k) \text{ for some } k. 
\end{align*}
\end{prop}


\remark 
In a similar way as in the proof of Corollary \ref{Cor cIVnr Uqbm sub}, 
we can prove that 
\begin{align*}
k_i^{\pm} \cdot ( \wt{v}_{\mu} \cdot L_j^{\lan c_j^{\mu} \ran} ) \in \cI_{\wt{V}_{n,r}}, 
\quad  
e_i \cdot ( \wt{v}_{\mu} \cdot L_j^{\lan c_j^{\mu} \ran} ) \in \cI_{\wt{V}_{n,r}}. 
\end{align*}
for $1\leq i \leq m-1$, $\mu \in \vL_n(m)$ and $1\leq j \leq n$. 
We can also prove that 
\begin{align*}
f_i \cdot ( \wt{v}_{\mu} \cdot L_j^{\lan c_j^{\mu} \ran} ) \in \cI_{\wt{V}_{n,r}} 
\text{ unless } \xi^{-1}(i) =(m_k,k) \text{ for some } k \text{ and } j=N_i^{\mu}. 
\end{align*}


\section{Shifted quantum affine algebras} 

In this section, 
we review the definition and fundamental properties of shifted quantum affine algebras given in \cite{FT19}. 

\begin{definition}[{\cite[\S 5.1]{FT19}}] 
For $q \in \CC^{\times} \setminus \{\pm 1\}$ and $\bb =(b_1,b_2,\dots, b_{m-1}) \in \ZZ^{m-1}$, 
we define the shifted quantum affine algebra 
$U_{q,[\bb]} = U_{q,[\bb]} (L \Fsl_m)$ by the following generators and defining relations: 
\begin{description}
\item[generators]  
$e_{i,t}, \, f_{i,t}, \, \psi_{i, s_i}^+, \, (\psi^+_{i, - b_i})^{-1}, \, \psi_{i,s}^-, \, (\psi_{i,0}^-)^{-1} $ 
($1\leq i \leq m-1$, $t \in \ZZ$, $s_i \in \ZZ_{\geq - b_i}$, $s \in \ZZ_{\leq 0}$) 

\item[defining relations]  
\begin{align*}
&\tag{U1} 
\psi_{i, - b_i}^+ (\psi_{i, - b_i}^+)^{-1} = 1 = (\psi_{i, -b_i}^+)^{-1} \psi_{i,-b_i}^+,  
	\quad 
	\psi_{i,0}^- (\psi_{i,0}^-)^{-1} = 1 = (\psi_{i,0}^-)^{-1} \psi_{i,0}^-, 
	\\
	& [\psi_{i}^{\ve} (z), \psi_j^{\ve'} (w)] =0 \quad (\ve, \ve' \in \{ \pm\}), 
\\
&\tag{U2} 
(z - q^{a_{ij}} w) e_i (z) e_j (w) = (q^{a_{ij}} z - w) e_j (w) e_i (z), 
\\
&\tag{U3} 
(q^{a_{ij}} z - w) f_i(z) f_j(w) = (z- q^{a_{ij}} w) f_j(w) f_i(z), 
\\
&\tag{U4} 
(z - q^{a_{ij}} w) \psi_i^{\ve} (z) e_j(w) = (q^{a_{ij}} z - w) e_j (w) \psi_i^{\ve} (z) 
	\quad (\ve \in \{\pm\}), 
\\
&\tag{U5}  
(q^{a_{ij}} z - w) \psi_i^{\ve} (z) f_j (w) = (z - q^{a_{ij}} w) f_j (w) \psi_i^{\ve} (z) 
	\quad ( \ve \in \{\pm\}), 
\\
&\tag{U6} 
[e_i(z), f_j(w)] 
	= \frac{\d_{i,j}}{q-q^{-1}} \d \left( \frac{z}{w} \right) (\psi_i^+ (z) - \psi_i^- (z)), 
\\
&\tag{U7} 
e_i(z) e_j(w) = e_j(w) e_i(z)  \text{ if } j \not= i,  i \pm 1, 
	\\ & 
	e_{i\pm 1}(w) \big( e_i(z_1) e_i(z_2) + e_i(z_2) e_i(z_1) \big) + \big( e_i(z_1) e_i(z_2) + e_i(z_2) e_i(z_1) \big) e_{i\pm 1}(w) 
	\\ & \quad 
	= [2] \big( e_i (z_1) e_{i \pm 1} (w) e_i (z_2) + e_i (z_2) e_{i \pm 1} (w) e_i (z_1) \big),   
\\
&\tag{U8}   
f_i(z) f_j(w) = f_j(w) f_i(z)  \text{ if } j \not= i,  i \pm 1, 
	\\ & 
	f_{i\pm 1}(w) \big( f_i(z_1) f_i(z_2) + f_i(z_2) f_i(z_1) \big) + \big( f_i(z_1) f_i(z_2) + f_i(z_2) f_i(z_1) \big) f_{i\pm 1}(w) 
	\\ & \quad 
	= [2] \big( f_i (z_1) f_{i \pm 1} (w) f_i (z_2) + f_i (z_2) f_{i \pm 1} (w) f_i (z_1) \big),   
\end{align*}
where $(a_{ij})_{1\leq i,j \leq m-1}$ is the Cartan matrix of type $A_{m-1}$, 
and we consider the generating series 
\begin{align*}
&e_i(z) = \sum_{t \in \ZZ} e_{i,t} z^{-t}, 
\quad 
f_i(z) = \sum_{t \in \ZZ} f_{i,t} z^{-t},  
\\
&\psi_i^+(z) = \sum_{t \geq - b_i} \psi_{i,t}^+ z^{-t}, 
\quad 
\psi_i^- (z) = \sum_{t \geq 0} \psi_{i,-t}^- z^t, 
\quad 
\d(z) = \sum_{t \in \ZZ} z^t.
\end{align*}
\end{description}
\end{definition}

\remark  
In the above definition, 
we consider a shift of only the positive part. 
In \cite[\S 5.1]{FT19},  
the definition of shifted quantum affine algebras are given for any shifts of both  positive and negative parts. 
However, they are isomorphic to the above algebras $U_{q,[\bb]}$ for some $ \bb \in \ZZ^{m-1}$ 
(see the comments before remark 5.2 in \cite[\S 5.1]{FT19} ).

\para  
We define the elements $h_{i, \pm t} \in U_{q,[\bb]}$ ($1\leq i \leq m-1$, $t \in \ZZ_{>0}$) by 
\begin{align*}
& (\psi_{i,- b_i}^+ z^{b_i} )^{-1} \psi_i^+(z)  = \exp \big( (q-q^{-1}) \sum_{t >0} h_{i,t} z^{-t} \big), 
\\
& (\psi_{i,0}^-)^{-1} \psi_i^-(z) = \exp \big( - (q-q^{-1}) \sum_{t >0} h_{i,-t}z^t \big). 
\end{align*}
Then, the relations (U4) and (U5) are replaced by 
\begin{align*}
\tag{U4'} 
& \psi_{i,-b_i}^+ e_{j,s} (\psi_{i,-b_i}^+)^{-1} = q^{a_{ij}} e_{j,s}, 
	\quad 
	\psi_{i,0}^- e_{j,s} ( \psi_{i,0} )^{-1} = q^{- a_{ij}} e_{j,s}, 
	\quad 
	[h_{i,t}, e_{j,s}] = \frac{ [ t a_{ij}] }{t} e_{j,s+t}, 
\\
\tag{U5'} 
& \psi_{i,-b_i}^+ f_{j,s} (\psi_{i,-b_i}^+)^{-1} = q^{ - a_{ij}} f_{j,s}, 
	\quad 
	\psi_{i,0}^- f_{j,s} ( \psi_{i,0} )^{-1} = q^{ a_{ij}} f_{j,s}, 
	\quad 
	[h_{i,t}, f_{j,s}] = - \frac{ [ t a_{ij}] }{t} f_{j,s+t} 
\end{align*}
respectively. 
In particular, we have 
\begin{align}
\label{eispm1 h e}
e_{i,s \pm 1} =  \frac{1}{[2]} [h_{i, \pm 1}, e_{i,s}], 
\quad 
f_{i, s \pm 1} = - \frac{1}{[2]} [h_{i, \pm 1}, f_{i, s }]  
\quad (s \in \ZZ). 
\end{align}
We also have 
\begin{align*}
&\psi_{i,s}^+ = (q-q^{-1}) [e_{i,s}, f_{i,0}] \quad (s >0), 
\quad 
\psi_{i,s}^- = - (q-q^{-1}) [e_{i,s}, f_{i,0}] \quad ( s < - b_i), 
\\
& \psi_{i,s}^+ - \psi_{i,s}^- = (q-q^{-1}) [e_{i, s}, f_{i,0}] \quad ( - b_i \leq s \leq 0)
\end{align*}
by the relation (U6). 
These relations imply the following lemma. 


\begin{lem}
\label{Lem gen Uqbb}
The algebra $U_{q,[\bb]}$ is generated by 
$e_{i,0}, f_{i,0}$, $\psi_{i,s_i}^+$ ($-b_i \leq s_i < 0$), $(\psi_{i, -b_i}^+)^{-1}$, 
$\psi_{i,0}^-$ , $(\psi_{i,0}^-)^{-1}$ and $h_{i,\pm 1}$ for $ i = 1,2 \dots, m-1$.  
\end{lem}


\para  
We can easily check that 
the elements $\psi_{i, - b_i}^+ \psi_{i,0}^-$ ($1\leq i \leq m-1$) are central elements of $U_{q,[\bb]}$. 

In the case where $\bb= \mathbf{0} = (0,0, \dots,0)$, 
we have 
\begin{align}
\label{iso Uq0 UqLsl}
U_{q,[\mathbf{0}]} / ( \psi_{i,0}^+ \psi_{i,0}^- -1 \mid 1\leq i \leq m-1 ) 
\cong U_q (L \Fsl_m) 
\end{align} 
as algebras, 
where $U_q (L \Fsl_m)$ is the quantum loop algebra associated with $\Fsl_m$. 
We denote corresponding generators of $U_q (L \Fsl_m)$ via the above isomorphism 
by same symbols. 

\para 
Put $c = k_0^+ k_1^+ \dots k_{m-1}^+ \in U_q (\wh{\Fsl}_m)$, 
then $c$ is the canonical central element of $U_q (\wh{\Fsl}_m)$. 
By \cite{Dri87} and \cite{Bec94}, it is known that 
there exists an isomorphism of algebras 
\begin{align}
\label{iso Uq whsl UqLsl}
\Psi : U_q (\wh{\Fsl}_m)/ (c-1) \ra U_q (L \Fsl_m)
\end{align} 
such that 
\begin{align*}
&\Psi (e_i) = e_{i,0}, \quad \Psi (f_i) = f_{i,0}, \quad \Psi (k_i^{\pm}) = \psi_{i,0}^{\pm} 
	\quad (1 \leq i \leq m-1),  
\\
& \Psi (k_0^{\pm}) = \psi_{1,0}^{\mp} \psi_{2,0}^{\mp} \dots \psi_{m-1,0}^{\mp}, 
\\
& \Psi (e_0) = [f_{m-1,0},[f_{m-2,0}, \dots, [f_{2,0}, f_{1,1}]_{q^{-1}} \dots ]_{q^{-1}}]_{q^{-1}} 
	(\psi_{1,0}^- \psi_{2,0}^- \dots \psi_{m-1,0}^-), 
\\
& \Psi (f_0) = (-q)^{m-2} ( \psi_{1,0}^+ \psi_{2,0}^+ \dots \psi_{m-1,0}^+) 
	[e_{m-1,0}, [e_{m-2,0}, \dots, [e_{2,0}, e_{1,-1}]_{q^{-1}} \dots ]_{q^{-1}} ]_{q^{-1}}. 
\end{align*}


\para
For $\bb = (b_1,\dots, b_{m-1}), \bc = (c_1,\dots, c_{m-1}), \bd = (d_1,\dots, d_{m-1}) \in (\ZZ_{\geq 0})^{m-1}$ 
such that $\bb = \bc + \bd$ and $ \Be = (\eta_1,\dots, \eta_{m-1}) \in (\CC^{\times})^{m-1}$, 
there is an algebra homomorphism, 
so-called a shift homomorphism,  
\begin{align}
\label{shift hom}
\iota_{[\bc, \bd]}^{\Be} : U_{q,[\bb]} \ra U_{q,[\mathbf{0}]} 
\end{align} 
such that 
\begin{align*}
e_i(z) \mapsto (1- \eta_i^{-1} z )^{c_i} e_i(z), 
\quad 
f_i(z) \mapsto ( 1 - \eta_i^{-1} z)^{d_i} f_i (z), 
\quad 
\psi_i^{\pm} (z) \mapsto ( 1 - \eta_i^{-1} z)^{c_i + d_i} \psi_i^{\pm} (z) 
\end{align*}
by  
\cite[Lemma 10.18]{FT19}, 
and the homomorphism $\iota_{[\bc, \bd]}^{\Be}$ is injective by \cite[Theorem 10.19]{FT19} 
(see also \cite[\S4.5]{Her23}). 
In particular, we have  
\begin{align}
\label{iota Be psi i -bi +}
\iota^{\Be}_{[\bc, \bd]} (\psi_{i,- b_i}^+) = (- \eta_i^{-1})^{b_i} \psi_{i,0}^+, 
\quad 
\iota^{\Be}_{[\bc, \bd]} (\psi_{i,0}^-) = \psi_{i,0}^-.
\end{align} 

\remark 
\label{Remark iota eta cd}
Recall that the elements $ \psi_{i, - b_i}^+ \psi_{i,0}^-$ ($1\leq i \leq m-1$) are central in $U_{q,[\bb]}$. 
By \eqref{iota Be psi i -bi +}, 
we see that  
the injective homomorphism $\iota_{[\bc, \bd]}^{\Be} : U_{q,[\bb]} \ra U_{q,[\mathbf{0}]} $ induces the 
injective homomorphism 
\begin{align*}
\iota^{\Be}_{[\bc, \bd]} : U_{q, [\bb]} / (  \psi_{i, - b_i}^+ \psi_{i,0}^- - (- \eta_i^{-1})^{ b_i} \mid 1 \leq i \leq m-1 )  \ra U_q (L \Fsl_m).
\end{align*}  

\para 
In \cite[Theorem 10.13]{FT19}, 
the coproduct $\D$ of $U_q (\wh{\Fsl}_m)$ is described by using the generators of $U_q ( L \Fsl_m)$  
through the isomorphism $U_q (\wh{\Fsl}_m)/ (c-1) \cong U_q (L \Fsl_m)$, 
and we can naturally lift this coproduct to the coproduct of $U_{q,[\mathbf{0}]}$. 
We denote this coproduct of $U_{q,[\mathbf{0}]}$ by the same symbol $\D$. 

For $\bb, \bc, \bd \in (\ZZ_{\geq 0})^{m-1}$ such that $\bb = \bc + \bd$ and $\Be \in (\CC^{\times})^{m-1}$,  
there exists an algebra homomorphism 
\begin{align*}
\D_{\bd, \bc} : U_{q,{[ \bb]}} \ra U_{q,{[ \bd ]}} \otimes U_{q, { [ \bc ]}}
\end{align*}
such that the diagram
\begin{align}
\label{comm diagram Ddc}
\xymatrix{
U_{q,[\bb]} \ar[r]^{\D_{\bd, \bc} \quad \quad }  \ar[d]_{\iota_{[\bc, \bd]}^{\Be}} 
	& U_{q,[\bd]} \otimes U_{q, [\bc]} \ar[d]^{\iota_{[\mathbf{0}, \bd]}^{\Be} \otimes \iota_{[\bc, \mathbf{0}]}^{\Be}}
\\
U_{q, [\mathbf{0}]} \ar[r]^{\D \qquad } & U_{q, [\mathbf{0}]} \otimes U_{q, [\mathbf{0}]  } 
}
\end{align}
is commute by \cite[Theorem 10.20]{FT19}. 

\para 
For $\bb \in (\ZZ_{\geq 0})^{m-1}$, put 
\begin{align*}
\BB_{[\bb]}
	&= \{ \Bb = (\b_{i,s}) \mid \b_{i,-b_i}, \b_{i,0} \in \CC^{\times}, \, \b_{i,s} \in \CC, \, 
	1 \leq i \leq m-1, \, - b_i < s < 0\}
	\\
	&= \bigcup_{1 \leq i \leq m-1} \big( \CC^{\times} \cup ( \bigcup_{- b_i < s < 0} \CC ) \cup \CC^{\times} \big). 
\end{align*}
For $\Bb \in \BB_{[\bb]}$, we can define the one-dimensional $U_{q,[\bb]}$-module 
$L_{\Bb} = \CC w$ by 
\begin{align*}
&e_{i,t} \cdot w= f_{i,t} \cdot w =0, 
\\
& \psi_{i,s_i}^+ \cdot w = 
	\begin{cases} 
		\b_{i,s_i} w  & \text{ if } - b_i \leq s_i \leq 0, 
		\\
		0 & \text{ if } s_i >0, 
	\end{cases}
\quad 
\psi_{i,s}^- \cdot w = 
	\begin{cases}
		\b_{i,s} w & \text{ if } - b_i \leq s \leq 0, 
		\\
		0 & \text{ if } s < -b_i
	\end{cases}
\end{align*}
for $1\leq i \leq m-1$, $t \in \ZZ$, $s_i \in \ZZ_{\geq -b_i}$ and  $s \in \ZZ_{\leq 0}$. 

\remark 
By \cite[Proposition 6.3]{Her23}, 
it is known that 
there is no finite dimensional $U_{q,[\bb]}$-module if $\bb \not\in (\ZZ_{\geq 0})^{m-1}$. 



\section{Schur-Weyl duality for the shifted quantum affine algebra and the Ariki-Koike algebra}
\label{Section SW}

In this section, 
we consider a shift of the 
$(U_q(\wh{\Fsl}_m), \sH_{n,r})$-bimodule $\wt{V}_{n,r} = \wh{V}^{\otimes n} \otimes_{\wh{\sH}_n} \sH_{n,r}$ 
studied in \S \ref{Section wtVnr}, 
and we establish the Schur-Weyl duality. 

\para 
By definition, 
we see that the canonical central element $c = k_0^+ k_1^+ \dots k_{m-1}^+ \in U_q (\wh{\Fsl}_m)$ acts trivially on $\wt{V}_{n,r}$. 
Then $U_q (L\Fsl_m)$ acts on $\wt{V}_{n,r}$ through the isomorphism \eqref{iso Uq whsl UqLsl}. 
We have the following proposition. 

\begin{prop}
\label{Prop act wtvmu}
For $\mu \in \vL_n(m)$ and $1\leq i \leq m-1$, we have the following. 
\begin{align*}
& e_{i,0} \cdot \wt{v}_{\mu} 
	= \d_{(\mu_{i+1} \not=0)} q^{- \mu_{i+1} +1} \wt{v}_{\mu + \a_i} \cdot 
		\big( 1 + \sum_{p=1}^{\mu_{i+1} -1} q^p \TT_{N_i^{\mu} +1, \, N_i^{\mu} +p} \big), 
\\
& f_{i,0} \cdot \wt{v}_{\mu} 
	= \d_{(\mu _i \not=0)} q^{- \mu_i +1} \wt{v}_{\mu - \a_i} \cdot 
		\big( 1 + \sum_{p=1}^{\mu_i-1} q^{p} \wt{\TT}_{N_i^{\mu} -1, \, N_i^{\mu} -p} \big), 
\\
& f_{i,1} \cdot \wt{v}_{\mu} 
	= \d_{(\mu_i \not=0)} q^{i - \mu_i +1} \wt{v}_{\mu - \a_i} \cdot 
	L_{N_i^{\mu}} \big( 1 + \sum_{p=1}^{\mu_i-1} q^{p} \wt{\TT}_{N_i^{\mu} -1, \, N_i^{\mu }-p} \big), 
\\
& \psi_{i,0}^{\pm} \cdot \wt{v}_{\mu} = q^{\pm (\mu_i - \mu_{i+1})} \wt{v}_{\mu}, 
\\
& h_{i,1} \cdot \wt{v}_{\mu} 
	= \wt{v}_{\mu} \cdot 
		\big( q^{i-1} \sum_{p=1}^{\mu_i} L_{N_{i-1}^{\mu} +p} - q^{i+1} \sum_{p=1}^{\mu_{i+1}} L_{N_i^{\mu} + p} \big), 
\\
& h_{i, -1} \cdot \wt{v}_{\mu} 
	= \wt{v}_{\mu} \cdot 
		\big( q^{-i+1} \sum_{p=1}^{\mu_i} L_{N_{i-1}^{\mu} +p}^{-1} - q^{ - i -1} \sum_{p=1}^{\mu_{i+1}} L_{N_i^{\mu} +p}^{-1} \big). 
\end{align*}
\begin{proof}
The formulas for the actions of $e_{i,0}$, $f_{i,0}$ and $\psi_{i,0}^{\pm}$ follows from \eqref{act Uqwhsl on wtVnr} 
through the isomorphism \eqref{iso Uq whsl UqLsl}. 

We consider the action of $h_{i, \pm 1}$. 
Through the isomorphism \eqref{iso Uq whsl UqLsl}, 
we can compute that 
\begin{align*}
\psi_{i, \pm s}^{\pm} \cdot ( v_j \otimes x^p ) = \pm (\d_{i,j} - \d_{i+1,j}) q^{ \pm i s}(q-q^{-1}) v_j \otimes x^{p \pm s} 
\end{align*}
in $U_q(L\Fsl_m)$-module $\wh{V}$. 
We note that 
$h_{i, \pm 1} = \pm (q-q^{-1})^{-1} (\psi_{i,0}^{\pm})^{-1} \psi_{i, \pm 1}^{\pm}$. 
Moreover, we can check that, 
for $1 \leq j_1 \leq j_2 \leq m$,  
\begin{align*}
h_{i, \pm 1} \cdot ( (v_{j_1} \otimes 1) \otimes (v_{j_2} \otimes 1)) 
= (h _{i, \pm 1} \otimes 1 + 1 \otimes h_{i, \pm 1}) \cdot ( (v_{j_1} \otimes 1) \otimes (v_{j_2} \otimes 1)) 
\end{align*}
in $\wh{V}^{\otimes 2}$ by \cite[Theorem 10.13]{FT19}. 
This implies that 
\begin{align*}
h_{i, \pm 1} \cdot \wt{v}_{\mu} = \sum_{k=0}^{n-1} 
	( \underbrace{1 \otimes \dots \otimes 1}_{k} \otimes h_{i, \pm 1} \otimes \underbrace{1 \otimes \dots \otimes 1}_{n-k-1} )
	\cdot \wt{v}_{\mu} 
\end{align*} 
in $\wt{V}_{n,r} = \wh{V}^{\otimes n} \otimes_{\wh{\sH}_n} \sH_{n,r}$. 
Then, we have the formula for the action of $h_{i, \pm 1}$.

We can also compute the action of $f_{i,1}$ by using the relation 
$f_{i,1}= - \frac{1}{[2]} [h_{i,1}, f_{i,0}]$ 
in \eqref{eispm1 h e}.
\end{proof}
\end{prop} 

\para 
Recall the $\sH_{n,r}$-submodule $\cI_{\wt{V}_{n,r}}$ of $\wt{V}_{n,r}$ generated by \eqref{gen wtVnr}. 
As seen in Proposition \ref{Prop fi wtvmu L notin IwtV}, 
the space $\cI_{\wt{V}_{n,r}}$ is not closed under the action of $U_q (L\Fsl_m)$. 
Then we consider the following shift. 

For $\bm =(m_1, m_2, \dots, m_r) \in (\ZZ_{>0})^r$ such that $m=m_1+m_2+\dots + m_r$,  
we take $\bb_{\bm} =(b_1, b_2,\dots, b_{m-1}) \in \{0,1\}^{m-1}$ as 
\begin{align}
\label{choice bi}
b_i = \begin{cases}
	1 & \text{ if } \xi^{-1} (i) = (m_k,k) \text{ for some } k, 
	\\
	0 & \text{ othewise}, 
	\end{cases}
\end{align}
and we consider the shifted quantum affine algebra $U_{q,[\bb_{\bm}]}$ with this shift $\bb_{\bm}$. 

For $\bQ=(Q_0, Q_1,\dots, Q_{r-1}) \in (\CC^{\times})^r$, 
take $\Bb_{\bm}^\bQ =(\b_{i,s}) \in \BB_{[\bb_{\bm}]}$ as 
\begin{align*}
\b_{i,s} = \begin{cases}
		1 & \text{ if } s=0, 
		\\
		- q^{-i} Q_k^{-1} & \text{ if } \xi^{-1} (i) = (m_k,k) \text{ for some } k \text{ and } s =-1, 
	\end{cases}
\end{align*} 
and we consider the one-dimensional $U_{q,[\bb_{\bm}]}$-module $L_{\Bb_{\bm}^{\bQ}} = \CC w$. 
Then we have that 
\begin{align}
\label{act w}
\begin{split} 
&e_{i,t} \cdot w = f_{i,t} \cdot w =0,
\\
& \psi_{i,s_i}^+ \cdot w = \begin{cases}
		- q^{-i} Q_k^{-1} w & \text{ if } \xi^{-1}(i) = (m_k,k) \text{ for some } k \text{ and } s_i=-1, 
		\\
		w & \text{ if } s_i=0, 
		\\
		0 & \text{ if } s_i >0, 
		\end{cases}
\\
& \psi_{i,s}^- \cdot w = \begin{cases}
		- q^{-i} Q_k^{-1} w & \text{ if } \xi^{-1}(i) = (m_k,k) \text{ for some } k \text{ and } s =-1, 
		\\
		w & \text{ if } s=0, 
		\\
		0 & \text{ if } s <-b_i 
		\end{cases}
\end{split} 
\end{align}
for $1\leq i \leq m-1$, $t \in \ZZ$, $s_i \in \ZZ_{\geq - b_i}$ amd $s \in \ZZ_{\leq 0}$.

\remark 
Put $\Be = (\eta_1, \eta_2,\dots, \eta_{m-1}) \in (\CC^{\times})^{m-1}$ as 
\begin{align*}
\eta_i = \begin{cases}
 		1 & \text{ unless } \xi^{-1}(i) = (m_k,k) \text{ for some } k,  
		\\
		q^i Q_k & \text{ if } \xi^{-1}(i) = (m_k,k) \text{ for some } k. 
	\end{cases} 
\end{align*}
Then, we have that 
\begin{align*}
\psi_{i, - b_i}^+ \psi_{i,0}^- \cdot w = ( - \eta_i^{-1})^{b_i} w 
\end{align*}
for $1\leq i \leq m-1$. 
See also Remark \ref{Remark iota eta cd}.  

\para 
Put 
\begin{align*}
V_{n,r} = L_{\Bb_{\bm}^{\bQ}} \otimes \wt{V}_{n,r}, 
\end{align*}
and the shifted quantum affine algebra $U_{q,[\bb_{\bm}]}$ acts on $V_{n,r}$
through the homomorphism 
$\D_{\bb_{\bm},0} : U_{q,[\bb_{\bm}]} \ra U_{q,[\bb_{\bm}]} \otimes U_{q, [0]}$, 
where we note the isomorphism \eqref{iso Uq0 UqLsl}. 
Then, $V_{n,r}$ has the $(U_{q,[\bb_{\bm}]}, \sH_{n,r})$-bimodule structure.

By \cite[Theorem 10.13]{FT19} together with the commutative diagram \eqref{comm diagram Ddc}, 
we can check that 
\begin{align}
\label{Dbm 0}
\begin{split}
&\D_{\bb_{\bm},0} ( \psi_{i, -b_i}^+ ) = \psi_{i, -b_i}^+ \otimes \psi_{i,0}^+, 
	\quad 
	\D_{\bb_{\bm},0} (\psi_{i,0}^-) = \psi_{i,0}^- \otimes \psi_{i,0}^-,  
\\
& \D_{\bb_{\bm},0} (e_{i,0}) 
	= e_{i,0} \otimes \psi_{i,0}^+ + 1 \otimes e_{i,0}, 
\\
& \D_{\bb_{\bm},0} (f_{i,0}) 
	\equiv f_{i,0} \otimes 1 + \psi_{i,0}^- \otimes f_{i,0} + \d_{(b_i=1)} \psi_{i, -1}^+ \otimes f_{i,1} 
		\mod \fX_{[\bb_{\bm}]}^+ \otimes \fX_{[0]}^-, 
\\
& \D_{\bb_{\bm},0} (h_{i,1}) 
	\equiv h_{i,1} \otimes 1 + 1 \otimes h_{i,1} 
	\mod \fX_{[\bb_{\bm}]}^+ \otimes \fX_{[0]}^-, 
\\
& \D_{\bb_{\bm},0} (h_{i,-1}) 
	\equiv h_{i, -1} \otimes 1 + 1 \otimes h_{i,-1} 
	\mod \fX_{[\bb_{\bm}]}^+ \otimes \fX_{[0]}^-,
\end{split}
\end{align}
where 
$\fX_{[\bb_{\bm}]}^+$ is the left ideal of $U_{q,[\bb_{\bm}]}$ generated by 
$\{ e_{i,t} \mid 1 \leq i \leq m-1, \, t \in \ZZ\}$ 
and 
$\fX_{[0]}^-$ is the left ideal of $U_{q,[0]}$ generated by  
$\{ f_{i,t} \mid 1 \leq i \leq m-1, \, t \in \ZZ \}$.  
Note that $\D_{\bb_{\bm},0} (\psi_{i,0}^-) = \psi_{i,0}^- \otimes \psi_{i,0}^-$, 
and we have the decomposition 
\begin{align*}
V_{n,r} = \bigoplus_{\mu \in \vL_n(m)} V_{n,r}^{\mu}, 
\quad 
V_{n,r}^{\mu} = \{ v \in V_{n,r} \mid \psi_{i,0}^- \cdot v = q^{- (\mu_i - \mu_{i+1})} v \, 1 \leq i \leq m-1\}. 
\end{align*}
This is also a decomposition of $\sH_{n,r}$-modules. 
By the previous arguments, 
we see that 
$V_{n,r}^{\mu}$ is generated by $v_\mu = w \otimes \wt{v}_{\mu}$ as the $\sH_{n,r}$-module,  
and there is an isomorphism of $\sH_{n,r}$-modules 
\begin{align}
\label{iso wtMmu Vmunr} 
\wt{M}^{\mu} \ra V_{n,r}^{\mu}, 
\quad 
x_{\mu} \mapsto v_{\mu}.
\end{align}
 
\para 
Let $\cI_{V_{n,r}}$ be the $\sH_{n,r}$-submodule of $V_{n,r}$ generated by 
\begin{align}
\{v_{\mu} \cdot L_j^{\lan c_j^{\mu} \ran} \mid \mu \in \vL_n(m), \, 1 \leq j \leq n \}, 
\end{align}
and we have the isomorphism of $\sH_{n,r}$-modules 
\begin{align}
\label{iso Vnr Inr Mmu} 
\begin{split}
&V_{n,r} / \cI_{V_{n,r}} = \bigoplus_{\mu \in \vL_n(m)} V_{n,r}^{\mu} / (\cI_{V_{n,r}} \cap V_{n,r}^{\mu}) 
\ra \bigoplus_{\mu \in \vL_{n,r}(\bm)} M^{\mu}, 
\\
&v_{\mu} + \cI_{V_{n,r}} \mapsto m_{\mu} \,\, (\mu \in \vL_n(m)) 
\end{split} 
\end{align}
by the isomorphism \eqref{iso wtMmu Vmunr} and Theorem \ref{Thm def rel mmu}. 

In fact, the space $\cI_{V_{n,r}}$ is also a $U_{q,[\bb_{\bm}]}$-submodule of $V_{n,r}$ as follows. 
By Proposition \ref{Prop act wtvmu}, \eqref{act w} and \eqref{Dbm 0}, 
we can easily check the following proposition. 
\begin{prop} 
\label{Prop action Uqbbmm Vnr}
For $\mu \in \vL_n(m)$ and $1\leq i \leq m-1$, 
we have the following.  
\begin{align*}
& e_{i,0} \cdot v_{\mu} = \d_{(\mu_{i+1} \not=0)} q^{- \mu_{i+1}+1} v_{\mu + \a_{i}} \cdot 
	\big( 1 + \sum_{p=1}^{\mu_{i+1} -1} q^p \TT_{N_i^{\mu} +1, \, N_i^{\mu} +p} \big), 
\\
& f_{i,0} \cdot v_{\mu} = \d_{(\mu_i \not=0)} 
	\begin{cases}
		\dis 
		q^{- \mu_i +1} v_{\mu - \a_i} \cdot \big( 1 + \sum_{p=1}^{\mu_i-1} q^p \wt{\TT}_{N_i^{\mu} -1, \, N_i^{\mu}-p} \big)
			\hspace{-15em} 
			\\
			& \text{ unless } \xi^{-1} (i) = (m_k,k) \text{ for some } k,  
		\\
		\dis 
		q^{- \mu_i+1} ( - Q_k^{-1}) v_{\mu -\a_i} 
			\cdot ( L_{N_i^{\mu}} - Q_k) \big( 1 + \sum_{p=1}^{\mu_i-1} q^p \wt{\TT}_{N_i^{\mu} -1, \, N_i^{\mu}-p} \big) 
			\hspace{-15em} 
			\\
			& \text{ if } \xi^{-1} (i) = (m_k,k) \text{ for some } k,  
	\end{cases}
\\
&\psi_{i, -b_i}^+ \cdot v_{\mu} = 
	\begin{cases}
		q^{\mu_i - \mu_{i+1}} v_{\mu} & \text{ unless }   \xi^{-1} (i) = (m_k,k) \text{ for some } k,  
		\\
		- q^{-i} Q_k^{-1} q^{\mu_i - \mu_{i+1}} v_{\mu} 
			& \text{ if } \xi^{-1} (i) = (m_k,k) \text{ for some } k,  
	\end{cases}
\\
& \psi_{i,0}^- \cdot v_{\mu} = q^{-(\mu_i - \mu_{i+1})} v_{\mu}, 
\\
& h_{i,1} \cdot v_{\mu} 
	= v_{\mu} \cdot 
		\big( q^{i-1} \sum_{p=1}^{\mu_i} L_{N_{i-1}^{\mu} +p} - q^{i+1} \sum_{p=1}^{\mu_{i+1}} L_{N_i^{\mu} + p} \big) 
		+ \d_{(\xi^{-1}(i) = (m_k,k))} \frac{ -q^i Q_k}{q-q^{-1}} v_{\mu}, 
\\
& h_{i, -1} \cdot v_{\mu} 
	= v_{\mu} \cdot 
		\big( q^{-i+1} \sum_{p=1}^{\mu_i} L_{N_{i-1}^{\mu} +p}^{-1} - q^{ - i -1} \sum_{p=1}^{\mu_{i+1}} L_{N_i^{\mu} +p}^{-1} \big) 
		+ \d_{(\xi^{-1}(i) = (m_k,k))} \frac{q^{-i} Q_k^{-1}}{q-q^{-1}} v_{\mu}.  
\end{align*}
\end{prop}

By using these formulas, we have the following corollary. 

\begin{cor}
\label{Cor cIVnr Uqbm sub} 
$\cI_{V_{n,r}}$ is a $U_{q,[\bb_{\bm}]}$-submodule of $V_{n,r}$. 

\begin{proof}
Throughout this proof, 
we represent a certain element of $\sH_{n,r}$ by $Y_{j}$ ($1\leq j \leq n$) 
if we do not need the explicit form of the element $Y_j$. 

By Lemma \ref{Lem gen Uqbb}, 
the algebra $U_{q,[\bb_{\bm}]}$ is generated by 
$e_{i,0}$ $f_{i,0}$, $\psi_{i, - b_i}^+$, $(\psi_{i, - b_i}^+)^{-1}$, $\psi_{i,0}^-$, $(\psi_{i,0}^-)^{-1}$ 
and $h_{i, \pm 1}$ \text{ for } $1 \leq i \leq m-1$.  
Then, it is enough to check the actions of these elements.

For $1 \leq i \leq m-1$, $1\leq j \leq n$ and $ \mu \in \vL_n(m)$, 
we have 
\begin{align*}
&e_{i,0} \cdot (v_{\mu} \cdot L_j^{\lan c_j^{\mu} \ran} )
\\
&= (e_{i,0} \cdot v_{\mu} ) \cdot L_j^{\lan c_j^{\mu} \ran}  
\\
&= \d_{(\mu_{i+1} \not=0)} q^{- \mu_{i+1}+1} v_{\mu + \a_{i}} \cdot 
	\big( 1 + \sum_{p=1}^{\mu_{i+1} -1} q^p \TT_{N_i^{\mu} +1, \, N_i^{\mu} +p} \big) 
	L_j ^{\lan c_j^{\mu} \ran} 
\\
&= \d_{(\mu_{i+1} \not=0)} q^{- \mu_{i+1}+1} v_{\mu + \a_{i}} \cdot 
	\begin{cases}
	\dis 
	L_j^{ \lan c_j^{\mu} \ran } 
		\big( 1 + \sum_{p=1}^{\mu_{i+1} -1} q^p \TT_{N_i^{\mu} +1, \, N_i^{\mu} +p} \big)  
		\quad 
		\text{ if } N_i^{\mu} \geq j \text{ or } j > N_{i+1}^{\mu}, 
	\\ \dis 
	\sum_{z=1}^{j- N_i^{\mu}} L_{N_i^{\mu}+z}^{\lan c_j^{\mu} \ran} Y_{N_i^\mu +z} 
		+ \d_{(N_{i+1}^{\mu} > j)} L_{j+1}^{\lan c_j^{\mu} \ran} Y_{j+1}
		\quad 
		\text{ if } N_{i+1}^{\mu} \geq j > N_i^{\mu}
	\end{cases}
\end{align*}
by Proposition \ref{Prop action Uqbbmm Vnr} and Lemma \ref{Lemma TT Lk} (\roiv). 
On the other hand, we have 
\begin{align*}
&c_j^{\mu + \a_i} = \begin{cases}
		c_j^{\mu} -1 
		& \text{ if } \xi^{-1}(i)=(m_k,k) \text{ for some } k \text{ and } j=N_i^{\mu} +1, 
		\\
		c_j^{\mu} & \text{otherwise}, 
	\end{cases}
\\
& c_{N_i^{\mu}+1}^{\mu + \a_i} = \begin{cases}
		c_{N_i^{\mu}+2}^{\mu + \a_i} -1 & \text{ if } \xi^{-1}(i)=(m_k,k) \text{ for some } k, 
		\\
		c_{N_i^{\mu} +2}^{\mu + \a_i} & \text{ otherwise}, 
	\end{cases}
\\
& c_{N_i^{\mu}+2}^{\mu + \a_i}  = c_{N_i^{\mu} +3}^{\mu + \a_i} = \dots = c_{N_{i+1}^{\mu}}^{\mu+\a_i}
\end{align*}
by the definition \eqref{def cjmu}. 
Then, we conclude that $e_{i,0} \cdot (v_{\mu} \cdot L_j^{\lan c_j^{\mu} \ran} ) \in \cI_{V_{n,r}}$. 

Next, we consider the action of $f_{i,0}$. 

Suppose $\xi^{-1} (i) \not= (m_k,k)$ for all $ 1 \leq k \leq r $, and we have 
\begin{align*}
&f_{i,0} \cdot (v_{\mu} \cdot L_j^{\lan c_j^{\mu} \ran}  ) 
\\
&= \d_{(\mu_i \not=0)}  q^{- \mu_i +1} v_{\mu - \a_i} \cdot 
	\big( 1 + \sum_{p=1}^{\mu_i-1} q^p \wt{\TT}_{N_i^{\mu} -1, \, N_i^{\mu}-p} \big) L_j^{\lan c_j^{\mu} \ran}
\\
&= \d_{(\mu_i \not=0)}  q^{- \mu_i +1} v_{\mu - \a_i} \cdot
	\begin{cases}
	\dis 
	L_j^{\lan c_j^{\mu} \ran} \big( 1 + \sum_{p=1}^{\mu_i-1} q^p \wt{\TT}_{N_i^{\mu} -1, \, N_i^{\mu}-p} \big)
		& \text{ if } j > N_i^{\mu} \text{ or } N_{i-1}^{\mu} +1 >j, 
	\\ 
	L_{j-1}^{\lan c_j^{\mu} \ran} Y_{j-1} + L_{N_i^{\mu}}^{\lan c_j^{\mu} \ran} Y_{N_i^{\mu}}
		& \text{ if } N_i^{\mu} \geq j > N_{i-1}^{\mu} 
	\end{cases}
\end{align*}	
by Proposition \ref{Prop action Uqbbmm Vnr} and Lemma \ref{Lemma TT Lk} (\rov). 
We also have 
\begin{align}
\label{cj mu - ai}
\begin{split}
&c_j^{\mu - \a_i} = \begin{cases}
	 c_j^\mu +1 & \text{ if } \xi^{-1} (i) = (m_k,k) \text{ for some } k 
	 	\text{ and } j = N_i^{\mu}, 
	\\
	c_j^{\mu} & \text{ otherwise}, 
	\end{cases}
\\
& c_{N_i^{\mu}}^{\mu - \a_i} = \begin{cases}
		c_{N_i^{\mu}-1}^{\mu - \a_i} +1 & \text{ if } \xi^{-1} (i) = (m_k,k) \text{ for some }k, 
		\\
		c_{N_i^{\mu}-1}^{\mu - \a_i} & \text{ otherwise}, 
	\end{cases}
\\
& c_{N_i^{\mu}-1}^{\mu -\a_i}  = c_{N_i^{\mu} -2}^{\mu - \a_i} = \dots = c_{N_{i-1}^{\mu} +1}^{\mu - \a_i} 
	\geq  c_{N_{i-1}^{\mu}}^{\mu - \a_i} 
\end{split} 
\end{align}
by the definition \eqref{def cjmu}. 
Then, we have $f_{i,0} \cdot (v_{\mu} \cdot L_j^{\lan c_j^{\mu} \ran} ) \in \cI_{V_{n,r}}$ 
if $\xi^{-1} (i) \not= (m_k,k)$ for all $ 1 \leq k \leq r $. 

Suppose $\xi^{-1} (i) = (m_k,k)$ for some $k$, and we have 
\begin{align*}
\begin{split}
&f_{i,0} \cdot (v_{\mu} \cdot L_j^{\lan c_j^{\mu} \ran}  ) 
\\
&= \d_{(\mu_i \not=0)}  q^{- \mu_i +1} ( - Q_k^{-1}) v_{\mu - \a_i} \cdot  
	\\ & \qquad 
	(L_{N_i^{\mu}} - Q_k) 
	\begin{cases}
	\dis 
	L_j^{\lan c_j^{\mu} \ran} \big( 1 + \sum_{p=1}^{\mu_i-1} q^p \wt{\TT}_{N_i^{\mu} -1, \, N_i^{\mu}-p} \big)
		& \text{ if } j > N_i^{\mu} \text{ or } N_{i-1}^{\mu} +1 >j, 
	\\ 
	L_{j-1}^{\lan c_j^{\mu} \ran} Y_{j-1} + L_{N_i^{\mu}}^{\lan c_j^{\mu} \ran} Y_{N_i^{\mu}}
		& \text{ if } N_i^{\mu} \geq j > N_{i-1}^{\mu} 
	\end{cases}
\end{split}
\end{align*}	
by Proposition \ref{Prop action Uqbbmm Vnr} and Lemma \ref{Lemma TT Lk} (\rov). 
Moreover, we see that 
\begin{align*}
&(L_{N_i^{\mu}} - Q_k) L_j^{\lan c_j^{\mu} \ran} 
	= \begin{cases} 
	L_j^{\lan c_j^{\mu} \ran}  (L_{N_i^{\mu}} - Q_k) & \text{ if } N_i^{\mu} >j, 
	\\
	L_j^{\lan c_j^{\mu} \ran}  Y_j + L_{N_i^{\mu}}^{\lan c_j^{\mu} \ran} Y_{N_i^{\mu}} 
		& \text{ if } j > N_i^{\mu}, 
	\\
	L_{N_i^{\mu}}^{\lan c_j^{\mu} \ran} Y_{N_i^{\mu}} + \sum_{z=1}^{N_i^{\mu}-1} L_{N_i^{\mu}  -z}^{\lan c_j^{\mu} \ran} Y_{N_i^{\mu} -z} 
		& \text{ if } j = N_i^{\mu}, 
	\end{cases} 
\\
&(L_{N_i^{\mu}} - Q_k) L_{j-1}^{\lan c_j^{\mu} \ran}  
	= L_{j-1}^{\lan c_j^{\mu} \ran}  (L_{N_i^{\mu}} - Q_k) 
	\text{ if } N_i^{\mu} \geq j 
\end{align*}
by Lemma \ref{Lemma Li Ljlan k ran}. 
We also have 
\begin{align*}
&(L_{N_i^{\mu}} - Q_k) L_{N_i^{\mu}}^{\lan c_j^{\mu} \ran} 
\\
&= L_{N_i^{\mu}}^{\lan c_j^{\mu} +1 \ran} + (q-q^{-1}) \sum_{z=1}^{N_i^{\mu} -1} L_z^{\lan c_j^{\mu} \ran} 
	\wt{\TT}_{N_i^{\mu} -1, z+1} L_{z+1} \TT_{z, N_i^{\mu}-1} 
	\text{ if } N_i^{\mu} \geq j > N_{i-1}^{\mu}
\end{align*}
by Lemma \ref{Lemma Lj Qk LJlankran} 
together with $c_{N_i^{\mu}}^{\mu} = c_{N_i^{\mu} -1}^\mu = \dots = c_{N_{i-1}^{\mu} +1}^\mu =k$ 
in the case where $\xi^{-1}(i)=(m_k,k)$ and $N_i^{\mu} \geq j > N_{i-1}^{\mu}$. 
Combining them, we have 
\begin{align}
\label{fi0 vmu Lj cjmu}
\begin{split} 
& f_{i,0} \cdot ( v_{\mu} \cdot L_j^{\lan c_j^{\mu} \ran}) 
\\
&= \d_{(\mu_i \not=0)}  q^{- \mu_i +1} ( - Q_k^{-1}) v_{\mu - \a_i} \cdot   
	\begin{cases}
		L_j^{\lan c_j^{\mu} \ran} Y_j + L_{N_i^{\mu}}^{\lan c_j^{\mu} \ran} Y_{N_i^{\mu}}
			& \text{ if } j > N_i^{\mu}, 
	\\
	L_j^{\lan c_j^{\mu} \ran} Y_j 
		& \text{ if } N_{i-1}^{\mu} +1 >j, 
	 \\ \dis 
	L_{N_i^{\mu} -1}^{\lan c_j^{\mu} \ran} Y_{N_i^{\mu}-1}	
	+ L_{N_i^{\mu}}^{\lan c_j^{\mu} +1 \ran} Y_{N_i^{\mu}} 
	+ \sum_{z=1}^{N_i^{\mu}-1} L_z^{\lan c_j^{\mu} \ran} Y_z 
		\hspace{-5em} 
		\\
		& \text{ if } N_i^{\mu} \geq    j > N_{i-1}^{\mu}. 
	\end{cases}
\end{split}
\end{align}

From the formula \eqref{fi0 vmu Lj cjmu} together with \eqref{cj mu - ai}, 
we see that $f_{i,0} \cdot ( v_{\mu} \cdot L_j^{\lan c_j^{\mu} \ran})  \in \cI_{V_{n,r}}$, 
where we note that $c_{N_i^{\mu} +1}^{\mu} = c_{N_i^{\mu}}^{\mu} +1$ if $\xi^{-1} (i) = (m_k,k)$, 
and $c_{j_1}^{\mu} \geq  c_{j_2}^{\mu}$ if $j_1 \geq j_2$. 

We can prove that 
$h_{i,\pm 1} \cdot (v_{\mu} \cdot L_j^{\lan c_j^{\mu} \ran}  )  \in \cI_{V_{n,r}}$ 
in a similar way using Proposition \ref{Prop action Uqbbmm Vnr} together with Lemma \ref{Lemma Li Ljlan k ran}. 
It is clear that $ \psi_{i, - b_i}^+ \cdot (v_{\mu} \cdot L_j^{\lan c_j^{\mu} \ran}  )  \in \cI_{V_{n,r}}$ 
and $\psi_{i,0}^- \cdot (v_{\mu} \cdot L_j^{\lan c_j^{\mu} \ran}  )  \in \cI_{V_{n,r}}$. 
\end{proof}
\end{cor}

Thanks to Corollary \ref{Cor cIVnr Uqbm sub}, 
the shifted quantum affine algebra 
$U_{q,[\bb_{\bm}]}$ acts on $V_{n,r} / \cI_{V_{n,r}}$ as the quotient module of $V_{n,r}$ by $\cI_{V_{n,r}}$,  
and the space $V_{n,r} / \cI_{V_{n,r}}$ has the $(U_{q,[\bb_{\bm}]}, \sH_{n,r})$-bimodule structure.  
We denote the representations corresponding to these actions by 
\begin{align*}
\rho_{n,r} : U_{q,[\bb_{\bm}]} \ra \End_{\CC}( V_{n,r}/ \cI_{V_{n,r}} ), 
\quad 
\s_{n,r} : \sH_{n,r}^{\opp} \ra \End_{\CC}( V_{n,r}/ \cI_{V_{n,r}} ). 
\end{align*}
Then, we have the following theorem. 

\begin{thm}
\label{Thm SW SQA AK} 
Assume that $q$ is not a root of unity and $m_k \geq n$ for  $k=1,2,\dots, r$. 
\begin{enumerate}
\item 
$V_{n,r} / \cI_{V_{n,r}} \cong \bigoplus_{\mu \in \vL_{n,r} (\bm)} M^{\mu}$ as right $\sH_{n,r}$-modules. 

\item 
$\Im \rho_{n,r} \cong \sS_{n,r}(\bm)$ as algebras, 
and these are quasi-hereditary algebras. 
We also have $ \Im \s_{n,r} \cong \sH_{n,r} $ as algebras. 

\item 
The shifted quantum affine algebra $U_{q,[\bb_{\bm}]}$ and the Ariki-Koike algebra $\sH_{n,r}$ 
satisfy the double centralizer property through the bimodule $V_{n,r} / \cI_{V_{n,r}}$, namely we have 
\begin{align*}
\Im \rho_{n,r} = \End_{\sH_{n,r}^{\opp}} (V_{n,r} / \cI_{V_{n,r}}), 
\quad 
\Im \s_{n,r} = \End_{U_{q,[\bb_{\bm}]}} (V_{n,r} / \cI_{V_{n,r}}). 
\end{align*}
\end{enumerate}

\begin{proof}
(\roi) has already proved in \eqref{iso Vnr Inr Mmu}, 
and we identify the space $V_{n,r}/ \cI_{V_{n,r}}$ with the space $\bigoplus_{\mu \in \vL_{n,r}(\bm)} M^{\mu}$ through this isomorphism. 
We will prove that $\Im \rho_{n,r} \cong \sS_{n,r} (\bm)$ as algebras, 
then the remaining statements follow from Theorem \ref{Thm CSA SW}. 

Since the action of $U_{q,[\bb_{\bm}]}$  on $V_{n,r} / \cI_{V_{n,r}}$ commutes with the action of $\sH_{n,r}$, 
we have 
\begin{align*}
\Im \rho_{n,r} \subset \End_{\sH_{n,r}^{\opp}} (V_{n,r}/ \cI_{V_{n,r}}) \cong \sS_{n,r}(\bm).
\end{align*}  
On the other hand, 
 the cyclotomic $q$-Schur algebra $\sS_{n,r}(\bm)$ is generated by 
$\sE_i$, $\sF_i$ ($1\leq i \leq m-1$) and $1_{\mu}$ ($\mu \in \vL_{n,r}(\bm)$) 
by Proposition \ref{Prop gen sSnr}. 
Thus, in order to prove that $\Im \rho_{n,r} \cong \sS_{n,r}(\bm)$, 
it is enough to prove that $\Im \rho_{n,r}$ contains the elements  
$\sE_i$, $\sF_i$ ($1\leq i \leq m-1$) and $1_{\mu}$ ($\mu \in \vL_{n,r}(\bm)$).

By Proposition \ref{Prop action Uqbbmm Vnr}, 
we have 
$\rho_{n,r} (e_{i,0}) = \sE_i$ and $\rho_{n,r} (f_{i,0}) = \sF_i$ for $1\leq i \leq m-1$. 
We show that the elements $1_{\mu}$ ($\mu \in \vL_{n,r}(\bm)$) belong to $\Im \rho_{n,r}$ 
by a similar argument to the proof of \cite[Theorem 3.4]{DG02} as follows. 

For $1 \leq i  \leq m-1$, $c \in \ZZ$ and $ t \in \ZZ_{\geq 0}$, 
put $K_i^+ = (\psi_{i,0}^-)^{-1}, K_i^{-} = \psi_{i,0}^- $ and 
\begin{align*}
\begin{bmatrix} 
K_i ; c \\ t 
\end{bmatrix} 
= \prod_{s=1}^t \frac{K_i^+ q^{c-s+1} - K_i^{-} q^{- c + s -1}}{q^s - q^{-s}}, 
\end{align*}
in $U_{q,[\bb_{\bm}]}$. 
Note that we identify the set $\vL_{n,r}(\bm)$ with the set $\vL_n(m)$ by the bijection \eqref{bij vLnrbm vLnm}. 
For $\la, \mu \in \vL_{n}(m)$ and $1\leq i \leq m-1$, we have 
\begin{align*}
\begin{bmatrix} K_i ; \la_i + \la_{i +1} + 2n \\ 2 \la_i + 2 n \end{bmatrix} m_{\mu} 
&= \prod_{s=1}^{ 2 \la_i + 2n} \frac{ q^{\mu_i - \mu_{i+1} + \la_i + \la_{i+1} + 2n - s+1} 
	- q^{  -(\mu_i - \mu_{i+1} + \la_i + \la_{i+1} + 2n - s+1)}}{q^s - q^{-s}} m_{\mu}
\\
&= \frac{ \prod_{s=1}^{2 \la_i + 2n} [\mu_i - \mu_{i+1} + \la_i + \la_{i+1} + 2n - (s-1)] }{ \prod_{s=1}^{2 \la_i + 2n} [s]}
	m_{\mu}. 
\end{align*}
We see that $\mu_i - \mu_{i+1} + \la_i + \la_{i+1} + 2n >0$ since $\la, \mu \in \vL_n(m)$. 
Thus we have 
\begin{align*}
\begin{bmatrix} K_i ; \la_i + \la_{i +1} + 2n \\ 2 \la_i + 2 n \end{bmatrix} m_{\mu}  \not=0 
&\LRa \mu_i - \mu_{i+1} + \la_i + \la_{i+1} + 2n - ( 2 \la_i +2n -1) > 0 
\\
&\LRa \la_{i+1 } - \la_i \geq \mu_{i+1} - \mu_i. 
\end{align*}
Moreover, we see that 
\begin{align*}
\begin{bmatrix} K_i ; \la_i + \la_{i +1} + 2n \\ 2 \la_i + 2 n \end{bmatrix} m_{\la} = m_{\la}.  
\end{align*}
Similarly, we have 
\begin{align*}
\begin{bmatrix} K_i ; \la_{i+1} - \la_i -1 \\ 2 n \la_{i+1} \end{bmatrix} m_{\mu} 
&= \frac{ \prod_{s=1}^{2 n \la_{i+1}} [\mu_i - \mu_{i+1} + \la_{i+1} - \la_i -1 - (s-1)]}{\prod_{s=1}^{2 n \la_{i+1}} [s] } m_{\mu}. 
\end{align*}
It is clear that 
$\begin{bmatrix} K_i ; \la_{i+1} - \la_i -1 \\ 2 n \la_{i+1} \end{bmatrix} =1$ if $\la_{i+1} =0$. 
In the case where  $\la_{i+1} \not=0$, we see that 
$\mu_i - \mu_{i+1} + \la_{i+1} - \la_i -1 - ( 2 n \la_{i+1} -1) <0$, 
and we have 
\begin{align*}
\begin{bmatrix} K_i ; \la_{i+1} - \la_i -1 \\ 2 n \la_{i+1} \end{bmatrix} m_{\mu} \not=0 
&\LRa \mu_i - \mu_{i+1} + \la_{i+1} - \la_i -1 <0 
\\
&\LRa \la_{i+1} - \la_i \leq \mu_{i+1} - \mu_i.  
\end{align*} 
Moreover, we see that 
\begin{align*}
\begin{bmatrix} K_i ; \la_{i+1} - \la_i -1 \\ 2 n \la_{i+1} \end{bmatrix} m_{\la} = m_{\la}. 
\end{align*}
For $1\leq i \leq m-1$ and $\la \in \vL_{n}(m)$, 
put 
\begin{align*}
K_i^{\la} = \begin{bmatrix} K_i ; \la_i + \la_{i +1} + 2n \\ 2 \la_i + 2 n \end{bmatrix} 
	\begin{bmatrix} K_i ; \la_{i+1} - \la_i -1 \\ 2 n \la_{i+1} \end{bmatrix}, 
\end{align*}
and $K^{\la} = K_1^{\la} K_2^{\la} \dots K_{m-1}^{\la}$. 
Then, by the above calculations, we have that 
\begin{align}
\label{Kla mmu}
\begin{split}
&K^{\la} \cdot m_{\la} = m_{\la}, 
\\
& K^{\la} \cdot m_{\mu} \not=0 
\LRa \begin{cases}
	\la_{i+1} - \la_i \geq \mu_{i+1} - \mu_i & \text{ if } \la_{i+1} =0 
	\\
	\la_{i+1} - \la_i  = \mu_{i+1} - \mu_i & \text{ if } \la_{i+1} \not=0 
	\end{cases}
	\text{ for all } 1 \leq i \leq m-1.
\end{split}
\end{align}

For $\la, \mu \in \vL_{n}(m)$, suppose that $\la \not=\mu$.  
Then,  there exists $k$ 
such that 
$\la_i = \mu_i$ for $ 1 \leq i \leq k$ and $\la_{k+1} \not= \mu_{k+1}$. 
This implies that $\la_{k+1} - \la_k \not= \mu_{k+1} - \mu_k$. 
Moreover we see that $\la_{k+1} - \la_k < \mu_{k+1} - \mu_k$ if $\la_{k+1} =0$. 
Thus, \eqref{Kla mmu} implies that 
$K^{\la} \cdot m_{\mu} = \d_{(\la=\mu)} m_{\la}$, and  we have $\rho_{n,r} (K^{\la} ) = 1_{\la}$. 
Now we conclude that 
$\Im \rho_{n,r} = \End_{\sH_{n,r}^{\opp}} (V_{n,r} / \cI_{V_{n,r}}) \cong \sS_{n,r}(\bm)$. 
\end{proof}
\end{thm}

\para 
\label{Par AA-form}
Recall $\AA = \ZZ[\hat{q}, \hat{q}^{-1}, \hat{Q_0}, \hat{Q}_1,\dots, \hat{Q}_{r-1}]$ and 
$\KK = \QQ(\hat{q}, \hat{q}^{-1}, \hat{Q_0}, \hat{Q}_1,\dots, \hat{Q}_{r-1})$ 
from the paragraph \ref{Para integral CSA}. 
We can define the shifted quantum affine algebra $U_{\hat{q},[\bb_{\bm}]}^{\KK}$ 
over $\KK$ with the parameter $\hat{q}$ 
in the same way. 
We can also define the $( U_{\hat{q},[\bb_{\bm}]}^{\KK} , \sH_{n,r}^{\KK})$-bimodule $V_{n,r}^{\KK} / \cI_{V_{n,r}^{\KK}}$ 
as in the above argument.  
Then, the statements of Theorem \ref{Thm SW SQA AK} over $\KK$ also hold.  
In particular, there exists the surjective homomorphism 
\begin{align*}
\rho_{n,r}^{\KK} : U_{q,[\bb_{\bm}]}^{\KK} \ra \End_{(\sH_{n,r}^{\KK})^{\opp}} (V_{n,r}^{\KK}/ \cI_{V_{n,r}^{\KK}}) 
	\cong \sS_{n,r}^{\KK}(\bm)
\end{align*} 
if $m_k \geq n$ for $k=1,2,\dots,r$. 

We can also define the right $\sH_{n,r}^{\AA}$-module $V_{n,r}^{\AA} / \cI_{V_{n,r}^{\AA}}$ in a similar manner, 
and we can regard $V_{n,r}^{\AA}/ \cI_{V_{n,r}^{\AA}}$ 
as an $\AA$-submodule of $V_{n,r}^{\KK}/ \cI_{V_{n,r}^{\KK}}$ in a natural way. 

Let $U_{\hat{q},[\bb_{\bm}]}^{\AA}$ be the $\AA$-subalgebra of $U_{\hat{q},[\bb_{\bm}]}^{\KK}$ generated by 
$e_{i,t}^d/ [d] !$, $f_{i,t}^d/ [d] !$ ($1\leq i \leq m-1$, $t \in \ZZ$, $d \in \ZZ_{>0}$) 
and 
$\begin{bmatrix} K_i ; c \\ t \end{bmatrix}$ ($1\leq i \leq m-1$, $c \in \ZZ$, $t \in \ZZ_{\geq 0}$), 
where we put $K_i^+=(\psi_{i,0}^-)^{-1}$ and $K_i^- = \psi_{i,0}^-$. 
Then we can check that the restriction of $\rho_{n,r}^{\KK}$ to $U_{q,[\bb_\bm]}^{\AA}$ 
induces the algebra homomorphism 
$\rho_{n,r}^{\AA} : U_{q,[\bb_{\bm}]}^{\AA} \ra \End_{(\sH_{n,r}^{\AA})^{\opp}} (V_{n,r}^{\AA}/ \cI_{V_{n,r}^{\AA}})$ 
in a similar way as in the proof of \cite[Theorem 8.1]{Wad16}. 
Moreover, we see that 
$\rho_{n,r}^{\AA} (e_{i,0}^d / [d]!) = \sE_i^d / [d] !$, $\rho_{n,r}^{\AA} (f_{i,0}^d / [d]!) = \sF_i^d / [d] !$ 
($1\leq i \leq m-1$, $d \in \ZZ_{\geq 0}$) 
and $\rho_{n,r}^{\AA} (K^{\la}) = 1_{\la}$ ($\la \in \vL_{n,r}(\bm)$). 
Thus, the homomorphism $\rho_{n,r}^{\AA}$ is surjective if $m_k \geq n$ for $k=1,2,\dots, r$ by Proposition \ref{Prop gen CSA AA}. 
Then, by taking a specialization, 
we have the statements of Theorem \ref{Thm SW SQA AK} over an arbitrary ring with any parameters. 
In particular, Theorem \ref{Thm SW SQA AK} holds even the case where $q$ is a root of unity 
by taking a specialization of $U_{q,[\bb_{\bm}]}^{\AA}$. 

\remark 
\label{Remark integral form}
In order to establish the Schur-Weyl duality and to consider the relationship with cyclotomic $q$-Schur algebras, 
it is enough to take the $\AA$-subalgebra $U_{\hat{q},[\bb_{\bm}]}^{\AA}$ of $U_{\hat{q},[\bb_{\bm}]}^{\KK}$ in the above. 
However, we do not known wheter $U_{\hat{q},[\bb_{\bm}]}^{\AA}$ is the \lq\lq correct" integral form of 
$U_{\hat{q}, [\bb_{\bm}]}^{\KK}$ or not. 
It seems that the \lq\lq correct" integral form is more complicated 
by looking at the argument to obtain the shifted loop Lie algebra from the shifted quantum affine algebra 
in \S \ref{Section shifted loop Lie alg}.



\section{Shifted loop Lie algebras} 
\label{Section shifted loop Lie alg}

In this section, we introduce the shifted loop Lie algebras. 
Then, we see that the universal enveloping algebra of the shifted loop Lie algebra 
can be regarded as the shifted quantum affine algebra at $q=1$ 
by taking a certain form of the shifted quantum affine algebra. 

We consider only the shifted loop Lie algebra associated with the special linear Lie algebra $\Fsl_m$ in this paper. 
We can also consider the shifted loop Lie algebra for an arbitrary finite dimensional simple Lie algebra 
in a similar manner.   

\begin{definition}
For $\Be =(\eta_1,\dots, \eta_{m-1}) \in (\CC^{\times})^{m-1}$ 
and $\bb =(b_1,  \dots, b_{m-1}) \in \ZZ_{\geq 0}^{m-1}$, 
we define the shifted loop Lie algebra $\fg_{\Be, [\bb]} = L_{[\bb]}^{\Be} \Fsl_m$ 
by the following generators and defining relations: 
\begin{description}
\item[generators] 
$E_{i,t}$, $F_{i,t}$, $H_{i,t}$  $(1\leq i \leq m-1, \, t \in  \ZZ)$
\item[defining relations] 
\begin{align*}
& \tag{L1} 
	[H_{i,t}, H_{j,s}] =0, 
\\ 
& \tag{L2} 
	[H_{i,t}, E_{j,s}] = a_{ij} E_{j,s+t}, 
	\quad 
	[H_{i,t}, F_{j,s}] = - a_{ij} F_{j,s+t}, 
\\
& \tag{L3} 
	[E_{i,t}, F_{j,s}] = \d_{i,j} \big( H_{i,s+t} + \d_{(b_i >0)} ( - \eta_i^{-1})^{b_i} H_{i, s+t+ b_i} \big), 
\\
& \tag{L4} 
	[E_{i,t}, E_{j,s}]= [F_{i,t}, F_{j,s}] =0 \text{ if } j \not= i \pm1, 
	\\ &
	[E_{i,t}, E_{i \pm 1, s}] = [E_{i, t+1}, E_{i \pm 1, s-1}], 
	\quad 
	[F_{i,t}, F_{i \pm 1, s}] = [F_{i, t+1}, F_{ i \pm 1, s-1}],  
\\
& \tag{L5} 
	[E_{i,s}, [E_{i,t}, E_{i\pm 1,u}]] = [F_{i,s}, [F_{i,t}, F_{i\pm 1, u}]] =0. 
\end{align*}
\end{description}
\end{definition}

\remark 
In the case where $\bb =0$, we have $\fg_{\Be, [0]} \cong L \Fsl_m$. 

\para 
For $1\leq i, j \leq m$ and $t \in \ZZ$, 
we define the element $\cE_{(i,j),t}^{[\bb]} \in \fg_{\Be, [\bb]}$ by 
\begin{align*}
\cE_{(i,j),t}^{[\bb]} = \begin{cases}
		H_{i,t} & \text{ if } i=j, 
		\\
		[E_{i,0}, [E_{i+1,0}, \dots, [E_{j-2,0}, E_{j-1,t}] \dots ]] & \text{ if } i < j, 
		\\
		[F_{i-1,0}, [F_{i-2,0}, \dots, [F_{j+1,0}, F_{j,t}] \dots ]] & \text{ if } i > j.  
	\end{cases}
\end{align*}
In particular, we have $\cE_{(i,i+1),t}^{[\bb]} = E_{i,t}$ and $\cE_{(i+1,i), t}^{[\bb]} = F_{i,t}$. 
Then, we can prove the following proposition by the standard argument (see e.g. \cite[Proposition 2.6]{Wad16}). 

\begin{prop}
\label{Prop basis fg eta b}
The set 
$\{\cE_{(i,j),t}^{[\bb]} \mid 1 \leq i,j \leq m, \, t \in \ZZ \}$ 
is a basis of $\fg_{\Be, [\bb]}$. \end{prop}

\begin{prop}
\label{Prop shift hom Lie}
For $\Be \in (\CC^{\times})^{n-1}$ and $\bb \in \ZZ_{\geq 0}^{m-1}$, 
there exists an injective homomorphism of Lie algebras 
\begin{align*}
\iota_{[\bb]^\sharp}^{\Be} : \fg_{\Be, [\bb]} \ra \fg_{\Be, [0]} \cong L\Fsl_m
\end{align*}
such that 
\begin{align*}
\iota_{[\bb]^\sharp}^{\Be} ( E_{i,t}) = E_{i,t}, 
\quad 
\iota_{[\bb]^\sharp}^{\Be} (F_{i,t}) = F_{i,t} + \d_{(b_i >0)} ( - \eta_i^{-1})^{b_i} F_{i, t+b_i}, 
\quad 
\iota_{[\bb]^\sharp}^{\Be} (H_{i,t}) =H_{i,t}.
\end{align*}

\begin{proof} 
We can check that the homomorphism $\iota_{[\bb]^\sharp}^{\Be}$ is well-defined by direct calculations. 
By definitions, it is clear that 
\begin{align}
\label{iotab0 cE 1}
\iota_{[\bb]^\sharp}^{\Be} (\cE_{(i,j),t}^{[\bb]}) = \cE_{(i,j),t}^{[0]} \text{ if } i \leq j.
\end{align}
On the other hand, for $1 \leq j < i \leq m -1$ and $s, t \in \ZZ$, 
we can prove that 
$[F_{i, s}, \cE_{(i,j),t}^{[0]}] = \cE_{(i+1,j),s+t}^{[0]}$ 
in a similar way as in the proof of \cite[Lemma 2.5]{Wad16}, and we have 
\begin{align}
\label{iotab0 cE 2}
\iota_{[\bb]^\sharp}^{\Be} (\cE_{(i,j),t}^{[\bb]}) = \cE_{(i,j),t}^{[0]} + \sum_{k >0} n_k \cE_{(i,j), t +k}^{[0]}  
\quad (n_k \in \ZZ) \text{ if } i >j. 
\end{align}
Then the equations \eqref{iotab0 cE 1} and \eqref{iotab0 cE 2} together with Proposition \ref{Prop basis fg eta b} 
imply that $\iota_{[\bb]^\sharp}^{\Be}$ is injective. 
\end{proof}
\end{prop}


\para 
We denote by $U_{[\bb]}$ the universal enveloping algebra of $\fg_{\Be, [\bb]}$, 
and we denote by $\D$ the usual coproduct of $U_{[\bb]}$. 
Through the injection $\iota_{[\bb]^\sharp}^{\Be}$, 
we can regard $U_{[\bb]}$ is a subalgebra of $U_{[0]}$. 
We also see that 
\begin{align*}
\D \circ \iota_{[\bb]^\sharp}^{\Be} (E_{i,t}) 
	&= E_{i,t} \otimes 1 + 1 \otimes E_{i,t}, 
\quad 
\D \circ \iota_{[\bb]^\sharp}^{\Be} (H_{i,t}) 
	= H_{i,t} \otimes 1 + 1 \otimes H_{i,t}, 
\\
\D \circ \iota_{[\bb]^\sharp}^{\Be} (F_{i,t}) 
	&= \D ( F_{i,t} + \d_{(b_i>0)} (-\eta_i^{-1})^{b_i} F_{i,t +b_i} )
	\\
	&= F_{i,t} \otimes 1 + 1 \otimes F_{i,t} + \d_{(b_i >0)} (-\eta_i^{-1})^{b_i} ( F_{i, t+b_i} \otimes 1 + 1 \otimes F_{i, t+b_i} )
	\\
	&= ( F_{i,t} + \d_{(b_i>0)} (-\eta_i^{-1})^{b_i} F_{i,t +b_i} ) \otimes 1 
		+ 1 \otimes ( F_{i,t} + \d_{(b_i>0)} (-\eta_i^{-1})^{b_i} F_{i,t +b_i} ). 
\end{align*}
These formulas imply the following proposition.

\begin{prop}
For $\Be \in (\CC^{\times})^{n-1}$ and $\bb \in \ZZ_{\geq 0}^{m-1}$, we have the following.  

\begin{enumerate} 
\item 
We have 
$\D \circ \iota_{[\bb]^\sharp}^{\Be} = (\iota_{[\bb]^\sharp}^{\Be} \otimes \iota_{[\bb]^\sharp}^{\Be}) \circ \D$, 
namely, the following  diagram is commute. 
\begin{align*}
\xymatrix{
U_{[\bb]} \ar[r]^{\D \quad \quad }  \ar[d]_{\iota_{[\bb]^\sharp}^{\Be}} 
	& U_{ [\bb]} \otimes U_{ [\bb]} \ar[d]^{ \iota_{[\bb]^\sharp}^{\Be}  \otimes \iota_{[\bb]^\sharp}^{\Be} }
\\
U_{ [\mathbf{0}]} \ar[r]^{\D  \quad \quad } & U_{ [\mathbf{0}]} \otimes U_{ [\mathbf{0}]  } 
}
\end{align*}
Moreover, we can regard $U_{[\bb]}$ as a Hopf subalgebra of $U_{[0]}$ through the injection $\iota_{[\bb]^\sharp}^{\Be}$. 

\item 
There exists an algebra homomorphism 
$\D_{\bb}^{\sharp} : U_{[\bb]} \ra U_{[\bb]} \otimes U_{[0]}$ such that 
\begin{align*}
&\D_{\bb}^\sharp (E_{i,t}) = E_{i,t} \otimes 1 + 1 \otimes E_{i,t}, 
\quad 
\D_{\bb}^\sharp (H_{i,t}) = H_{i,t} \otimes 1 + 1 \otimes H_{i,t}, 
\\
&\D_{\bb}^{\sharp} (F_{i,t}) = F_{i,t} \otimes 1 + 1 \otimes ( F_{i,t} + \d_{(b_i>0)} (-\eta_i^{-1})^{b_i} F_{i,t +b_i} ).
\end{align*} 
Moreover, we have $\D \circ \iota_{[\bb]^\sharp}^{\Be} = (\iota_{[\bb]^\sharp}^{\Be} \otimes \Id) \circ \D_{\bb}^{\sharp}$, 
namely, the following diagram is commute. 
\begin{align*}
\xymatrix{
U_{[\bb]} \ar[r]^{\D_{\bb}^{\sharp} \quad \quad }  \ar[d]_{\iota_{[\bb]^\sharp}^{\Be}} 
	& U_{ [\bb]} \otimes U_{ [ \mathbf{0} ]} \ar[d]^{ \iota_{[\bb]^\sharp}^{\Be}  \otimes \Id }
\\
U_{ [\mathbf{0}]} \ar[r]^{\D  \quad \quad } & U_{ [\mathbf{0}]} \otimes U_{ [\mathbf{0}]  } 
}
\end{align*}
\end{enumerate}
\end{prop}

\para 
In the rest of this section, 
we give some connection for the shifted loop Lie algebras and the shifted quantum affine algebras. 

For $\bb =(b_1,\dots, b_{m-1}) \in \ZZ_{\geq 0}^{m-1}$ and 
$\Be = (\eta_1, \dots, \eta_{m-1}) \in (\CC^{\times})^{m-1}$, 
we denote by $U_{q,\Be, [\bb]}$ the quotient algebra of $U_{q,[\bb]}$ 
by the two-sided ideal generated by $\{\psi_{i, - b_i}^+ \psi_{i,0}^- - (- \eta_i^{-1})^{ b_i} \mid 1 \leq i \leq m-1\}$, 
namely, we have 
\begin{align*}
U_{q,\Be, [\bb]} 
=U_{q, [\bb]} / (  \psi_{i, - b_i}^+ \psi_{i,0}^- - (- \eta_i^{-1})^{ b_i} \mid 1 \leq i \leq m-1). 
\end{align*} 

We also conisder  an injective shift homomorphism 
\begin{align}
\label{shift hom bbsharp}
\iota_{[\bb]^{\sharp}}^{\Be} : U_{q,[\bb]} \ra U_{q, [0]}
\end{align}
defined by 
\begin{align*}
&e_i(z) \mapsto e_i(z), 
\quad 
f_i (z) \mapsto \big( 1 + \d_{(b_i>0)} (- \eta_i^{-1} z)^{b_i} \big) f_i(z), 
\\
&
\psi_i^{\pm}(z) \mapsto \big( 1 + \d_{(b_i>0)} (- \eta_i^{-1} z)^{b_i} \big) \psi_i^{\pm} (z).
\end{align*}
We have 
\begin{align*}
&\iota_{[\bb]^{\sharp}}^{\Be} (e_{i,t}) = e_{i,t}, 
\quad 
\iota_{[\bb]^{\sharp}}^{\Be} (f_{i,t}) = f_{i,t} + \d_{(b_i >0)} ( -\eta_i^{-1})^{b_i} f_{i, t+ b_i}, 
\\
& \iota_{[\bb]^{\sharp}}^{\Be} (\psi_{i,t}^+) 
	=\begin{cases}
		( -\eta_i^{-1})^{b_i} \psi_{i, t+b_i}^+ & \text{ if } 0 > t \geq -b_i, 
		\\
		\psi_{i,t}^+ + \d_{(b_i>0)} ( -\eta_i^{-1})^{b_i} \psi_{i, t+b_i}^+ 
			& \text{ if } t \geq 0, 
	\end{cases}
\\
& \iota_{[\bb]^{\sharp}}^{\Be} (\psi_{i,t}^- ) 
	= \begin{cases}
		\psi_{i,t}^- & \text{ if } 0 \geq t > - b_i, 
		\\
		\psi_{i,t}^- + \d_{(b_i >0)} ( -\eta_i^{-1})^{b_i} \psi_{i, t + b_i}^- & \text{ if } - b_i \geq t, 
	\end{cases}
\end{align*}
and we see that the shift homomorphism $\iota_{[\bb]^{\sharp}}^{\Be}$ induces the injective  homomorphism 
$\iota_{[\bb]^{\sharp}}^{\Be} : U_{q,\Be, [\bb]} \ra U_q (L\Fsl_m)$. 

\para 
For $1\leq i \leq m-1$ and $t \in \ZZ$, 
we define the elements 
$K_i^{\pm}$ and $H_{i,t}$ of $U_{q, \Be, [\bb]}$ by 
\begin{align*}
&K_i^+ = ( \psi_{i,0}^-)^{-1}, \quad K_i^- = \psi_{i,0}^-, 
\\
& H_{i,t} = \frac{1}{q-q^{-1}} 
	\begin{cases}
	\dis 
	\sum_{k=0}^{\zeta_{b_i} (t)}  (-1)^k ( - \eta_i^{-1})^{-(k+1) b_i} \psi_{i, t- (k+1) b_i}^+ 
		& \text{ if } t >0, 
	\\ \dis 
	( -\eta_i^{-1})^{-b_i} \psi_{i, -b_i}^+ - \psi_{i,0}^- 
		& \text{ if } t=0, 
	\\ \dis 
	- \sum_{k=0}^{\zeta_{b_i} (t)} (-1)^k ( - \eta_i^{-1})^{k b_i} \psi_{i, t+ k b_i}^- 
		& \text{ if } t <0,
	\end{cases}
\end{align*}
where the integer $\zeta_{b_i}(t) \in \ZZ_{\geq 0}$ is determined by 
\begin{align*}
 t = \begin{cases}
 		\zeta_{b_i} (t) \cdot b_i + r & \text{ if } t \geq 0, 
		\\
		- \zeta_{b_i} (t) \cdot b_i - r & \text{ if } t <0
 	\end{cases}
\end{align*}
with $ 0 \leq r \leq b_i-1$ in the case where $b_i >0$, 
and we put $\zeta_{b_i}(t)=0$ if $b_i=0$. 

For an example, 
in the case where $b_i=3$, 
we have 
\begin{align*}
&H_{i,6} = \frac{ (- \eta_i^{-1})^{-3} \psi_{i,3}^+ - (- \eta_i^{-1})^{-6} \psi_{i,0}^+ + (- \eta_i^{-1})^{-9} \psi_{i,-3}^+}{q-q^{-1}}, 
\\
& H_{i,5} = \frac{ (- \eta_i^{-1})^{-3} \psi_{i,2}^+ - (- \eta_i^{-1})^{-6} \psi_{i,-1}^+}{q-q^{-1}}, 
	\quad 
	H_{i,4} = \frac{ (- \eta_i^{-1})^{-3} \psi_{i,1}^+ - (- \eta_i^{-1})^{-6} \psi_{i,-2}^+}{q-q^{-1}}, 
	\\
	& \qquad 
	H_{i,3} = \frac{ (- \eta_i^{-1})^{-3} \psi_{i,0}^+ - (- \eta_i^{-1})^{-6} \psi_{i,-3}^+}{q-q^{-1}}, 
\\
& H_{i,2} = \frac{ (-\eta_i^{-1})^{-3} \psi_{i,-1}^+}{q-q^{-1}}, 
	\quad 
	 H_{i,1} = \frac{ (-\eta_i^{-1})^{-3} \psi_{i,-2}^+}{q-q^{-1}}, 
\\
& H_{i,0} = \frac{ (-\eta_i^{-1})^{-3} \psi_{i,-3}^+ - \psi_{i,0}^-}{q-q^{-1}}, 
\\
& H_{i,-1} = \frac{ - \psi_{i,-1}^-}{q-q^{-1}}, 
	\quad 
	H_{i,-2} = \frac{ - \psi_{i,-2}^-}{q-q^{-1}}, 
\\
& H_{i,-3} = \frac{ - \psi_{i,-3}^- + (-\eta_i^{-1})^3 \psi_{i,0}^-} {q-q^{-1}}, 
	\quad 
	 H_{i,-4} = \frac{ - \psi_{i,-4}^- + (-\eta_i^{-1})^3 \psi_{i,-1}^-} {q-q^{-1}},
	 \\ & \qquad 
	  H_{i,-5} = \frac{ - \psi_{i,-5}^- + (-\eta_i^{-1})^3 \psi_{i,-2}^-} {q-q^{-1}},
\\
& H_{i, -6} = \frac{ - \psi_{i,-6}^- + (-\eta_i^{-1})^3 \psi_{i, -3}^- - (-\eta_i^{-1})^6 \psi_{i,0}^-}{q-q^{-1}}.
\end{align*}

By direct calculation, we obtain the following lemma. 

\begin{lem}
\label{Lemma rel wtUqBebb}
For $\bb \in \ZZ_{\geq 0}^{m-1}$ and $\Be  \in (\CC^{\times})^{m-1}$, 
We have the following relations in $U_{q, \Be, [\bb]}$. 
\begin{align*}
\tag{$\wt{U}1$} 
	& K_i^+ K_i^- = K_i^- K_i^+ =1, 
		\quad K_i^+ - K_i^- = (q-q^{-1}) H_{i,0}, 
	\\
	& [K_i^+, K_j^+] = [K_i^+, H_{j,s}] = [ H_{i,t}, H_{j,s}] =0, 
\\
\tag{$\wt{U}2$}
	&e_{i, t+1} e_{j,s} - q^{a_{ij}} e_{j,s} e_{i,t+1}  =   q^{a_{ij}} e_{i,t} e_{j, s+1} - e_{j,s+1} e_{i,t}, 
\\
\tag{$\wt{U}3$}
	&f_{i,t+1} f_{j,s}  - q^{- a_{ij}} f_{j,s} f_{i,t+1} = q^{- a_{ij}} f_{i,t} f_{j,s+1}  -  f_{j,s+1} f_{i,t}, 
\\
\tag{$\wt{U}4$}  
	& K_i^+ e_{j,s} K_i^- = q^{a_{ij}} e_{j,s}, 
	\\
	& H_{i,0} e_{j,s} - q^{a_{ij}} e_{j,s} H_{i,0} = [a_{ij}] e_{j,s} K_i^-, 
	\\
	& H_{i,t+1} e_{j,s} - q^{a_{ij}} e_{j,s} H_{i,t+1} = q^{a_{ij}} H_{i,t} e_{j, s+1} - e_{j,s+1} H_{i,t}, 
\\
\tag{$\wt{U}5$} 
	& K_i^+ f_{j,s} K_i^- = q^{ - a_{ij}} e_{j,s}, 
	\\
	& H_{i,0} f_{j,s} - q^{-a_{ij}} f_{j,s} H_{i,0} = [-a_{ij}] f_{j,s} K_i^-, 
	\\
	& H_{i,t+1} f_{j,s} - q^{ - a_{ij}} f_{j,s} H_{i,t+1} = q^{ - a_{ij}} H_{i,t} f_{j, s+1} - f_{j,s+1} H_{i,t}, 
\\
\tag{$ \wt{U}6$} 
	& [e_{i,t}, f_{j,s}] = \d_{ij} \big( H_{i,s+t} + \d_{(b_i >0)} (- \eta_i^{-1})^{b_i} H_{i, s+t +b_i} \big), 
\\ 
\tag{$\wt{U}7$}
	& [e_{i,t}, e_{j,s}] =0 \quad  \text{ if } j \not=i, i \pm1, 
	\\
	&e_{i \pm 1,u} (e_{i,s} e_{i,t} + e_{i,t} e_{i,s} ) + (e_{i,s} e_{i,t} + e_{i,t} e_{i,s} ) e_{i \pm 1, u} 
	\\ & \quad 
	= (q+q^{-1}) (e_{i,s} e_{i\pm1 ,u} e_{i,t} + e_{i,t} e_{i \pm 1, u} e_{i,s}), 
\\
\tag{$\wt{U}8$}
	& [f_{i,t}, f_{j,s}] =0 \quad  \text{ if } j \not=i, i \pm1,  
	\\
	&f_{i \pm 1,u} (f_{i,s} f_{i,t} + f_{i,t} f_{i,s} ) + (f_{i,s} f_{i,t} + f_{i,t} f_{i,s} ) f_{i \pm 1, u} 
	\\ & \quad 
	= (q+q^{-1}) (f_{i,s} f_{i\pm1 ,u} f_{i,t} + f_{i,t} f_{i \pm 1, u} f_{i,s}). 
\end{align*}
\end{lem}

We can also check the following proposition on the shift homomorphism 
$\iota_{[\bb]^\sharp}^{\Be}$ directly. 

\begin{prop}
\label{Prop iotabbsharp}
For $\bb \in \ZZ_{\geq 0}^{m-1}$ and $\Be  \in (\CC^{\times})^{m-1}$, 
we have 
\begin{align*}
&\iota_{[\bb]^\sharp}^{\Be} (e_{i,t}) = e_{i,t}, 
\quad 
\iota_{[\bb]^\sharp}^{\Be} (f_{i,t}) = f_{i,t} + \d_{(b_i >0)} ( -\eta_i^{-1})^{b_i} f_{i, t + b_i}, 
\\
&\iota_{[\bb]^\sharp}^{\Be} (H_{i,t}) = H_{i,t}, 
\quad 
\iota_{[\bb]^\sharp}^{\Be} (K_i^{\pm}) = K_i^{\pm}.
\end{align*}
\end{prop}

\begin{definition}
For $q \in \CC^{\times}$, $\bb =(b_1,\dots, b_{m-1}) \in \ZZ^{m-1}$ 
and $ \Be =(\eta_1,\dots, \eta_{m-1}) \in (\CC^{\times})^{m-1}$, 
we define the algebra 
$\wt{U}_{q, \Be, [\bb]} $ by the generators 
$e_{i,t}$, $f_{i,t}$, $H_{i,t}$, $K_i^{\pm}$ ($1\leq i \leq m-1$, $ t \in \ZZ$) 
and the defining relations 
$(\wt{U}1)$ - $(\wt{U}8)$. 
\end{definition}

\begin{prop}
Assume that $q \not= \pm 1$, 
then there exists an algebra isomorphism 
$\Phi : \wt{U}_{q,\Be, [\bb]} \ra U_{q, \Be, [\bb]}$ such that 
\begin{align*}
&\Phi  (e_{i,t}) = e_{i,t}, 
\quad 
\Phi  (f_{i,t} ) = f_{i,t}, 
\quad 
\Phi  (K_i^+) = (\psi_{i,0}^-)^{-1}, 
\quad 
\Phi (K_i^-) = \psi_{i,0}^-, 
\\
& \Phi (H_{i,t}) 
	= \frac{1}{q-q^{-1}} 
	\begin{cases}
	\dis 
	\sum_{k=0}^{\zeta_{b_i} (t)}  (-1)^k ( - \eta_i^{-1})^{-(k+1) b_i} \psi_{i, t- (k+1) b_i}^+ 
		& \text{ if } t >0, 
	\\ \dis 
	( -\eta_i^{-1})^{-b_i} \psi_{i, -b_i}^+ - \psi_{i,0}^- 
		& \text{ if } t=0, 
	\\ \dis 
	- \sum_{k=0}^{\zeta_{b_i} (t)} (-1)^k ( - \eta_i^{-1})^{k b_i} \psi_{i, t+ k b_i}^- 
		& \text{ if } t <0.	
\end{cases}
\end{align*}
\begin{proof}
The well-definedness of $\Phi$ follows from Lemma \ref{Lemma rel wtUqBebb}. 
The inverse homomorphism $\Psi : U_{q,\Be, [\bb]} \ra \wt{U}_{q, \Be, [\bb]}$ 
is given by 
\begin{align*}
& \Psi(e_{i,t}) = e_{i,t}, 
\quad 
\Psi(f_{i,t}) = f_{i,t}, 
\\
&
\Psi( \psi_{i, -b_i}^+) = ( - \eta_i^{-1})^{b_i} K_i^+, 
\quad 
\Psi( (\psi_{i, -b_i}^+)^{-1} ) = ( - \eta_i^{-1})^{- b_i} K_i^-, 
\\
&\Psi( \psi_{i,0}^- ) = K_i^-, 
\quad 
\Psi( (\psi_{i,0}^-)^{-1} ) = K_i^+, 
\\
& \Psi (\psi_{i,t}^+ ) = 
	\begin{cases}
	(q-q^{-1}) \big( H_{i,t} + \d_{(b_i >0)} (- \eta_i^{-1})^{b_i} H_{i, t+b_i} \big) 
		& \text{ if }t >0, 
	\\
	K_i^+ + \d_{(b_i >0)} (q-q^{-1}) (- \eta_i^{-1})^{b_i} H_{i, b_i} 
		& \text{ if } t=0, 
	\\
	(q-q^{-1}) ( - \eta_i^{-1})^{b_i} H_{i, t+ b_i} 
		&\text{ if } 0 > t > -b_i, 
	\end{cases}
\\
& \Psi( \psi_{i,t}^- ) =
	\begin{cases}
	- (q-q^{-1}) H_{i,t}  
		& \text{ if } 0 > t > -b_i, 
	\\
	(- \eta_i^{-1})^{b_i} K_i^- - \d_{(b_i>0)} (q-q^{-1}) H_{i, - b_i} 
		&\text{ if } t = -b_i, 
	\\
	- (q-q^{-1}) \big( H_{i,t} + \d_{(b_i>0)} ( -\eta_i^{-1})^{b_i} H_{i, t+b_i} \big) 
		& \text{ if } - b_i >t. 
\end{cases}
\end{align*}
\end{proof}
\end{prop}

\para 
We consider the algebra 
$\wt{U}_{q=1,\Be, [\bb]}$ in the case where  $q=1$, 
then we have $K_i^{+} = \pm 1$ by the relation ($\wt{U}1$). 
In this case, we easily see that 
the universal enveloping algebra $U (\fg_{\Be, [\bb]})$ of the shifted loop Lie algebra $\fg_{\Be, [\bb]}$ 
is isomorphic to the algebra $\wt{U}_{q=1, \Be, [\bb]} / (K_i^+ -1 \mid 1 \leq i \leq m-1)$ 
by the isomorphism 
\begin{align*}
&U (\fg_{\Be, [\bb]}) \ra \wt{U}_{q=1, \Be, [\bb]}/ (K_i^+ -1 \mid 1 \leq i \leq m-1), 
\\
&E_{i,t} \mapsto e_{i,t}, 
\quad 
F_{i,t} \mapsto f_{i,t}, 
\quad 
H_{i,t} \mapsto H_{i,t}.
\end{align*}
Thus,  we can regard the universal enveloping algebra $U (\fg_{\Be, [\bb]})$ 
as the shifted quantum affine algebra $U_{q,\Be, [\bb]}$ at $q=1$. 
Moreover, the compatibility of the shift homomorphism $\iota_{[\bb]^\sharp}^{\Be}$ 
follows from Proposition \ref{Prop iotabbsharp}. 

\remark 
\label{Remark bbm wtU}
We take $\bb_{\bm}$ as in \eqref{choice bi}, 
then we see that the shift homomorphism 
$\iota_{[\bb_{\bm}]^{\sharp}}^{\Be} $ in \eqref{shift hom bbsharp} 
coincides with the shift homomorphism 
$\iota_{[0, \bb_{\bm}]}^{\Be} $ in  \eqref{shift hom}. 
Moreover, for the one-dimensional $U_{q,[\bb_{\bm}]}$-module  
$L_{\Bb_{\bm}^{\bQ}} = \CC w$ given by \eqref{act w}, 
we see that 
$e_{i,t} \cdot w= f_{i,t} \cdot w  = H_{i,t} \cdot w=0$ 
for $1\leq i \leq m-1$ and $t \in \ZZ$. 


\section{Schur-Weyl duality for the shifted loop Lie algebra and the Ariki-Koike algebra at $q=1$}

In this section, we prove that the Schur-Weyl duality is obtained through a certain tensor space $\VV^{\otimes n}$ 
in the case where $q=1$. 

\para 
The affine Hecke algebra $\wh{\sH}_n^{q=1}$ over $\CC$ at $q=1$ is isomorphic to the semidirect product 
$\CC[\fS_n] \ltimes \CC[X_1^{\pm}, X_2^{\pm}, \dots , X_n^{\pm}]$, where $\CC[\fS_n]$ is the group algebra of $\fS_n$.  
We denote by $T_w$ the element of $\CC[\fS_n]$ corresponding to $w \in \fS_n$.  
We also denote by $\sH_{n,r}^{q=1}$ the Ariki-Koike algebra over $\CC$ with parameters 
$q=1$ and $\bQ =(Q_0,Q_1,\dots, Q_{r-1}) \in (\CC^{\times})^r$. 

\para 
The algebra $U_{[0]} = U(L\Fsl_m)$ acts on $\wh{V}$ by 
\begin{align*}
&E_{i,t} \cdot (v_j \otimes x^k) = \d_{i+1, j} v_{j-1} \otimes x^{k+t}, 
\\
&F_{i,t} \cdot (v_j \otimes x^k) = \d_{i, j} v_{j+1} \otimes x^{k+t}, 
\\
&H_{i,t} \cdot (v_j \otimes x^k) = (\d_{i,j} - \d_{i+1,j}) v_j \otimes x^{k+t}, 
\end{align*}
and $U_{[0]}$ also act on the tensor space $\wh{V}^{\otimes n}$ through the coproduct $\D$. 

On the other hand, 
$\wh{\sH}_n^{q=1}$ acts on $\wh{V}^{\otimes n}$ from right  by 
\begin{align*}
&\Big( (v_{j_1} \otimes v_{j_2} \otimes \dots \otimes v_{j_n}) \otimes x_1^{k_1} x_2^{k_2} \dots x_n^{k_n} \Big) 
	\cdot T_w X_1^{l_1} X_2^{l_2} \dots X_n^{l_n} 
	\\
	& \quad 
	= ( v_{j_{w^{-1}(1)}} \otimes v_{j_{w^{-1}(2)}} \otimes \dots \otimes v_{j_{w^{-1}(n)}} ) \otimes 
		x_1^{k_{w^{-1}(1)} + l_1 } x_2^{k_{w^{-1}(2)} +l_2 } \dots x_n^{k_{w^{-1}(n)} + l_n}. 
\end{align*}
Then, the actions of $U_{[0]}$ and of $\wh{\sH}_n^{q=1}$ on $\wh{V}^{\otimes n}$ commute with each other.  


\para 
\label{para setting Ubm} 
We take $\bb_{\bm} =(b_1,b_2,\dots, b_{m-1}) \in \{0,1\}^{m-1}$ as in \eqref{choice bi}. 
We also take 
$\Be =(\eta_1, \eta_2,\dots, \eta_{m-1}) \in (\CC^{\times})^{m-1}$ as 
\begin{align*}
\eta_i = \begin{cases}
 		1 & \text{ unless } \xi^{-1}(i) = (m_k,k) \text{ for some } k,  
		\\
		Q_k & \text{ if } \xi^{-1}(i) = (m_k,k) \text{ for some } k. 
	\end{cases} 
\end{align*}
We consider the shifted loop Lie algebra $\fg_{\Be, [\bb_{\bm}]}$ 
and its universal enveloping algebra 
$U_{[\bb_{\bm}]}= U(\fg_{\Be, [\bb_{\bm}]})$. 

\para 
We consider the one-dimensional $U_{[\bb_{\bm}]}$-module 
$L = \CC w$ defined by 
\begin{align*}
E_{i,t} \cdot w = F_{i,t} \cdot w = H_{i,t} \cdot w =0 
\end{align*}
for $1\leq i \leq m-1$ and $ t \in \ZZ$ 
(see Remark \ref{Remark bbm wtU}). 

Put 
\begin{align*}
V_{n,r} = L \otimes ( \wh{V}^{\otimes n} \otimes_{\wh{\sH}_n^{q=1}} \sH_{n,r}^{q=1}),
\end{align*}
and $U_{[\bb_{\bm}]}$ acts on $V_{n,r}$ through the homomorphism 
$\D_{\bb_{\bm}}^{\sharp} : U_{[\bb_{\bm}]} \ra U_{[\bb_{\bm}]} \otimes U_{[0]}$. 
We have the weight space  decomposition 
\begin{align*}
V_{n,r} = \bigoplus_{\mu \in \vL_n(m)} V_{n,r}^{\mu}, 
\quad 
V_{n,r}^{\mu} = \{ v \in V_{n,r} \mid H_{i,0} \cdot v = (\mu_i - \mu_{i+1}) v\}, 
\end{align*}
and this is also a decomposition of  $\sH_{n,r}^{q=1}$-modules. 
Moreover, we see that 
$V_{n,r}^{\mu}$ is generated by $v_{\mu}$ as the $\sH_{n,r}^{q=1}$-module,  
and there exists an isomorphism of $\sH_{n,r}^{q=1}$-modules 
\begin{align}
\label{iso wtMmu Vmunr 2} 
\wt{M}^{\mu} \ra V_{n,r}^{\mu}, \quad x_{\mu} \mapsto v_{\mu} 
\end{align}
in a similar reason to Lemma \ref{Lemma iso wtMmu wt Vnrmu}. 
By direct calculations,  we have 
\begin{align}
\label{act Ubm Vnrmu} 
\begin{split} 
& E_{i,0} \cdot v_{\mu} = \d_{(\mu_{i+1}\not=0)} v_{\mu + \a_i} 
	\cdot \big( 1 +  \sum_{p=1}^{\mu_{i+1} -1} \TT_{N_i^{\mu} +1, \, N_i^{\mu} +p} \big), 
\\
& F_{i,0} \cdot v_{\mu} = \d_{(\mu_i \not=0)} 
	\begin{cases}\dis 
	v_{\mu - \a_i} \cdot  \big(  1 + \sum_{p=1}^{\mu_i -1} \wt{\TT}_{N_i^{\mu}-1, \, N_i^{\mu} -p} \big) 
		\hspace{-20em}
		\\
		& \text{ unless } \xi^{-1} (i) =(m_k,k) \text{ for some }k, 
	\\ \dis 
	- Q_k^{-1} v_{\mu - \a_i} \cdot ( L_{N_i^{\mu}} - Q_k ) 
		\big(  1 + \sum_{p=1}^{\mu_i -1} \wt{\TT}_{N_i^{\mu}-1, \, N_i^{\mu} -p} \big) 
		\hspace{-15em}
		\\
		& \text{ if } \xi^{-1} (i) =(m_k,k) \text{ for some }k, 
	\end{cases} 
\\
& H_{i,0} \cdot v_{\mu} = (\mu_i - \mu_{i+1}) v_{\mu}, 
\\
& H_{i,1} \cdot v_{\mu} = (\mu_i - \mu_{i+1}) v_{\mu} \cdot 
	\big( \sum_{p=1}^{\mu_i} L_{N_{i-1}^{\mu}+p} + \sum_{p=1}^{\mu_{i+1}} L_{N_i^{\mu} +p} \big), 
\\
& H_{i,-1} \cdot v_{\mu} = (\mu_i - \mu_{i+1}) v_{\mu} \cdot 
	\big( \sum_{p=1}^{\mu_i} L_{N_{i-1}^{\mu}+p}^{-1} + \sum_{p=1}^{\mu_{i+1}} L_{N_i^{\mu} +p}^{-1} \big). 
\end{split} 
\end{align}

\para 
Let $\cI_{V_{n,r}}$ be the $\sH_{n,r}^{q=1}$-submodule of $V_{n,r}$ generated by 
\begin{align*}
\{ v_{\mu} \cdot L_j^{\lan c_j^{\mu} \ran} \mid \mu \in \vL_n (m), \, 1 \leq j \leq n \}, 
\end{align*}
where we remark that   
$L_j^{\lan k \ran} = (L_j - Q_0) (L_j- Q_1) \dots (L_j - Q_{k-1})$ in the case where $q=1$. 
Then, we have the isomorphism of $\sH_{n,r}^{q=1}$-modules 
\begin{align}
\label{iso Vnr I VNr Mmu}
\begin{split}
&V_{n,r} / \cI_{V_{n,r}} = \bigoplus_{\mu \in \vL_n(m)} V_{n,r}^{\mu} / (\cI_{V_{n,r}} \cap V_{n,r}^{\mu}) 
\ra \bigoplus_{\mu \in \vL_{n,r}(\bm)} M^{\mu}, 
\\
&v_{\mu} + \cI_{V_{n,r}} \mapsto m_{\mu} \,\, (\mu \in \vL_n(m)) 
\end{split} 
\end{align}
by the isomorphism \eqref{iso wtMmu Vmunr 2} and Theorem \ref{Thm def rel mmu}. 
Moreover, we can check that 
$\cI_{V_{n,r}}$ is also a $U_{[\bb_{\bm}]}$-submodule by using \eqref{act Ubm Vnrmu} 
in a similar way to Corollary \ref{Cor cIVnr Uqbm sub}. 
Then, the space $V_{n,r} / \cI_{V_{n,r}}$ has the $(U_{[\bb_{\bm}]}, \sH_{n,r}^{q=1})$-bimodule structure,  
and we can prove 
the statements of Theorem \ref{Thm SW SQA AK} at $q=1$ in a similar manner.  

\para 
In the case where $q=1$, 
we can replace the space $V_{n,r} / \cI_{V_{n,r}}$ with a certain tensor space $\VV^{\otimes n}$ as follows. 

For $\bm=(m_1, m_2,\dots, m_r) \in \ZZ_{>0}^r$ and $1\leq j \leq m$, 
we determin the integer $c_j$ by 
\begin{align*}
c_j =k  \text{ if } (m_0 + m_1+ \dots + m_{k-1}) < j \leq (m_0+m_1 + \dots + m_{k-1} + m_k), 
\end{align*}
where we put $m_0=0$. 

We consider the natural action of $\CC[x^{\pm}]$ on $\wh{V} = V \otimes \CC[x^{\pm}]$.
Let $\cI_{\bQ}$ be the $\CC [x^{\pm}]$-submodule of $\wh{V}$ generated by 
\begin{align*}
\{ v_j \otimes (x -Q_0) (x - Q_1) \dots (x - Q_{c_j-1}) \mid 1 \leq j \leq m \}, 
\end{align*}
and take the quotient space $\VV = \wh{V} / \cI_{\bQ}$. 
For an element $v \in \wh{V}$, 
we denote the element $v + \cI_{\bQ} \in \wh{V} / \cI_{\bQ}$ by $v$ simply. 
Then, the $\CC$-vector space $\VV$ has a basis 
\begin{align*}
\{ v_j \otimes x^{p} \mid 1 \leq j \leq m, \, 0 \leq p \leq c_j -1 \}. 
\end{align*}

We can define the action of $U_{[\bb_{\bm}]}$ on $\VV$ by 
\begin{align*}
&E_{i,t} \cdot (v_j \otimes x^p) = \d_{i+1, j} v_{j-1} \otimes x^{p+t}, 
\\
&F_{i,t} \cdot (v_j \otimes x^p) = \d_{i, j} 
	\begin{cases}
	v_{j+1} \otimes x^{p+t} & \text{ unless } \xi^{-1}(i) = (m_k,k) \text{ for some }k, 
	\\
	v_{j+1}\otimes (x^{p+t}  - Q_k^{-1} x^{p+t+1}) & \text{ if } \xi^{-1}(i) = (m_k,k) \text{ for some }k, 
	\end{cases} 
\\
&H_{i,t} \cdot (v_j \otimes x^p) = (\d_{i,j} - \d_{i+1,j}) v_j \otimes x^{p+t}. 
\end{align*}

We consider the action of $U_{[\bb_{\bm}]}$ on the tensor space $\VV^{\otimes n}$ 
through the coproduct $\D$ of $U_{[\bb_{\bm}]}$. 
We also define the right action of $\sH_{n,r}^{q=1}$ on $\VV^{\otimes n}$ by  
\begin{align*}
&\big( (v_{j_1} \otimes x^{p_1}) \otimes (v_{j_2} \otimes x^{p_2}) \otimes \dots \otimes (v_{j_n} \otimes x^{p_n})  \big) 
	\cdot T_w L_1^{l_1} L_2^{l_2} \dots L_n^{l_n} 
	\\
	&\quad 
	= \big( (v_{j_{w^{-1}(1)}} \otimes x^{p_{w^{-1}(1)} +l_1}) \otimes (v_{j_{w^{-1}(2)} } \otimes x^{p_{w^{-1}(2)} +l_2}) 
		\otimes \dots \otimes (v_{j_{w^{-1}(n)}} \otimes x^{p_{w^{-1}(n)} + l_n})  \big). 
\end{align*}
Then, we can easily check that 
the actions of $U_{[\bb_{\bm}]}$ and of  $\sH_{n,r}^{q=1}$ on $\VV^{\otimes n}$ commute with each other. 

\para 
We have 
\begin{align*}
\VV^{\otimes n} = \bigoplus_{\mu \in \vL_n(m)} \VV^{\otimes n}_{\mu}, 
\quad 
\VV_{\mu}^{\otimes n} = \{ v \in \VV^{\otimes n} \mid H_{i,0} \cdot v = (\mu_i - \mu_{i+1}) v \quad (1 \leq i \leq m-1)\}, 
\end{align*}
and this is also a decomposition of $\sH_{n,r}^{q=1}$-modules. 
Moreover, we see that 
\begin{align}
\label{basis VVnmu} 
\left\{  (v_{j_1} \otimes x^{p_1}) \otimes (v_{j_2} \otimes x^{p_2}) \otimes \dots \otimes (v_{j_n} \otimes x^{p_n})  
	\mid \begin{matrix}
		\mu_i = \sharp \{ k \mid j_k = i \} \, (1\leq i \leq m), 
		\\
		0 \leq p_k \leq c_{j_k} -1 \quad (1\leq k \leq n) 
		\end{matrix}
		\right\}  
\end{align}
gives a basis of $\VV_{\mu}^{\otimes n}$. 
For $\mu \in \vL_n (m)$, put 
\begin{align*}
v_{\mu}' = \underbrace{(v_1 \otimes 1) \otimes \dots \otimes (v_1 \otimes 1)}_{\mu_1} 
	\otimes \dots \otimes \underbrace{ (v_m \otimes 1) \otimes \dots \otimes (v_m \otimes 1)}_{\mu_m}, 
\end{align*}
and we see that 
$V_{\mu}^{\otimes n}$ is generated by $v_{\mu}'$ as the $\sH_{n,r}^{q=1}$-module. 
We can also check that 
\begin{align*}
v_{\mu}' \cdot T_i = v_{\mu}'  \text{ if } s_i \in \fS_{\mu} 
\text{ and }
v_{\mu}' \cdot L_j^{\lan c_j^{\mu} \ran} =0 \text{ for } 1 \leq j \leq n. 
\end{align*}
Then, we have the surjective homomorphism of $\sH_{n,r}^{q=1}$-modules 
\begin{align}
\label{iso Mmu Vmun}
M^{\mu} \ra \VV_{\mu}^{\otimes n}, 
\quad 
m_{\mu} \mapsto v_{\mu}'
\end{align}
by Theorem \ref{Thm def rel mmu}, 
and we can check that this is an isomorphism by compairing their dimensions using Proposition \ref{Prop basis Mmu} (\roii) and 
\eqref{basis VVnmu}.
Moreover, we have 
\begin{align}
\label{act Ubm VVmu} 
\begin{split} 
& E_{i,0} \cdot v_{\mu}' = \d_{(\mu_{i+1}\not=0)} v_{\mu + \a_i}' 
	\cdot \big( 1 +  \sum_{p=1}^{\mu_{i+1} -1} \TT_{N_i^{\mu} +1, \, N_i^{\mu} +p} \big), 
\\
& F_{i,0} \cdot v_{\mu}' = \d_{(\mu_i \not=0)} 
	\begin{cases}\dis 
	v_{\mu - \a_i}' \cdot  \big(  1 + \sum_{p=1}^{\mu_i -1} \wt{\TT}_{N_i^{\mu}-1, \, N_i^{\mu} -p} \big) 
		\hspace{-20em}
		\\
		& \text{ unless } \xi^{-1} (i) =(m_k,k) \text{ for some }k, 
	\\ \dis 
	- Q_k^{-1} v_{\mu - \a_i}' \cdot ( L_{N_i^{\mu}} - Q_k ) 
		\big(  1 + \sum_{p=1}^{\mu_i -1} \wt{\TT}_{N_i^{\mu}-1, \, N_i^{\mu} -p} \big) 
		\hspace{-15em}
		\\
		& \text{ if } \xi^{-1} (i) =(m_k,k) \text{ for some }k, 
	\end{cases} 
\\
& H_{i,0} \cdot v_{\mu}' = (\mu_i - \mu_{i+1}) v_{\mu}', 
\\
& H_{i,1} \cdot v_{\mu}' = (\mu_i - \mu_{i+1}) v_{\mu}' \cdot 
	\big( \sum_{p=1}^{\mu_i} L_{N_{i-1}^{\mu}+p} + \sum_{p=1}^{\mu_{i+1}} L_{N_i^{\mu} +p} \big), 
\\
& H_{i,-1} \cdot v_{\mu}' = (\mu_i - \mu_{i+1}) v_{\mu}' \cdot 
	\big( \sum_{p=1}^{\mu_i} L_{N_{i-1}^{\mu}+p}^{-1} + \sum_{p=1}^{\mu_{i+1}} L_{N_i^{\mu} +p}^{-1} \big). 
\end{split} 
\end{align}

We denote the representations corresponding to the actions of $U_{[\bb_{\bm}]}$ and of $\sH_{n,r}^{q=1}$ 
on $\VV^{\otimes n}$ by 
\begin{align*}
\rho_{n,r} : U_{q,[\bb_{\bm}]} \ra \End_{\CC}( \VV^{\otimes n} ), 
\quad 
\s_{n,r} : \sH_{n,r}^{\opp} \ra \End_{\CC}( \VV^{\otimes n} ) 
\end{align*}
respectively. 
Then we have the following theorem. 


\begin{thm}
\label{Thm SW SLA AK} 
Assume that $m_i \geq n$ for all $i=1,2, \dots, r$. 
\begin{enumerate}
\item 
$\VV^{\otimes n} \cong V_{n,r} / \cI_{V_{n,r}}$ as $(U_{[\bb_{\bm}]}, \sH_{n,r}^{q=1})$-bimodules. 
\item 
$\VV^{\otimes n} \cong \bigoplus_{\mu \in \vL_{n,r} (\bm)} M^{\mu}$ as right $\sH_{n,r}^{q=1}$-modules. 

\item 
$\Im \rho_{n,r} \cong \sS_{n,r}^{q=1}(\bm)$ as algebras, 
and these are quasi-hereditary algebras. 
We also have $ \Im \s_{n,r} \cong \sH_{n,r}^{q=1} $ as algebras. 

\item 
The algebras $U_{[\bb_{\bm}]}$ and $\sH_{n,r}^{q=1}$ 
satisfy the double centralizer property through the bimodule $\VV^{\otimes n}$, namely we have 
\begin{align*}
\Im \rho_{n,r} = \End_{(\sH_{n,r}^{q=1})^{\opp}} ( \VV^{\otimes n} ), 
\quad 
\Im \s_{n,r} = \End_{U_{[\bb_{\bm}]}} ( \VV^{\otimes n} ). 
\end{align*}
\end{enumerate}
\begin{proof}
(\roi). 
The isomorphisms \eqref{iso Vnr I VNr Mmu} and \eqref{iso Mmu Vmun} imply 
the isomorphism $\VV^{\otimes n} \cong V_{n,r} / \cI_{V_{n,r}}$ of $\sH_{n,r}^{q=1}$-modules, 
and we can check that it is also an isomorphism of $U_{[\bb_{\bm}]}$-modules 
by \eqref{act Ubm Vnrmu} and \eqref{act Ubm VVmu}, 
where we note that  $U_{[\bb_{\bm}]}$ is generated by 
$E_{i,0}$, $F_{i,0}$, $H_{i,0}$, $H_{i, \pm 1}$ ($1 \leq i \leq m-1$). 

The statement (\roii) follows from 
the isomorphism \eqref{iso Mmu Vmun}, 
and we can prove the statement (\roiii) and (\roiv) 
in a similar arguments as in the proof of Theorem \ref{Thm SW SQA AK}. 
\end{proof}
\end{thm} 

\remark 
By Theorem \ref{Thm SW SLA AK} (\roi), 
we can regard the $(U_{q,[\bb_{\bm}]}, \sH_{n,r})$-bimodule $V_{n,r} / \cI_{V_{n,r}}$ 
in \S \ref{Section SW} as a $q$-analogue of the tensor space $\VV^{\otimes n}$ in this section. 
However, we can not construct the tensor space $\VV^{\otimes n}$ as the $U_{q,[\bb_{\bm}]}$-module directly 
since the shifted quantum affine algebra $U_{q,[\bb_{\bm}]}$ does not have a Hopf algebra structure although 
it is a coideal subalgebra of $U_{q,[0]}$.  
Therefore, we need  to construct the bimodule $V_{n,r} / \cI_{V_{n,r}}$ as in \S \ref{Section SW} 
to establish the Schur-Weyl duality, 
but it seems that the bimodule $V_{n,r} / \cI_{V_{n,r}}$ is natural one in the representation theory of shifted quantum affine algebras 
thanks to Theorem \ref{Thm SW SLA AK} (\roi). 



\appendix
\section{Proofs of Lemma \ref{Li lan k ran}, \ref{Lemma TT Lk},  \ref{Lemma Li Ljlan k ran} and \ref{Lemma Lj Qk LJlankran}}
\label{Appendix proof technical lemmas}

\para 
{\it (Proof of Lemma  \ref{Li lan k ran}).} 
\begin{proof}
We prove Lemma  \ref{Li lan k ran} by the induction on $i$. 

In the case where $i=1$, 
note that $L_1 =T_0$, and 
we have 
\begin{align*}
L_1^{\lan k \ran} 
&= (T_0 - Q_0) (T_0 - Q_1) \dots (T_0 - Q_{k-1}) 
\\
&= L_1^k + \sum_{p=0}^{k-1} (-1)^{k-p} \big( \sum_{ 0\leq i_1<i_2<\dots< i_{k-p} \leq k-1} Q_{i_1} Q_{i_2} \dots Q_{i_{k-p}} \big) L_1^p.
\end{align*}

In the case where $i>1$, 
we have $L_i^{\lan k \ran} = T_{i-1} L_{i-1}^{\lan k \ran} T_{i-1}$. 
Applying the inductive hypothesis,   
we have 
\begin{align*}
L_i^{\lan k \ran} 
= T_{i-1} \big( L_{i-1}^k + \sum_{0 \leq p_1, p_2,\dots, p_{i-1} \leq k-1} L_1^{p_1} L_2^{p_2} \dots L_{i-1}^{p_{i-1}} 
	h_{(p_1,p_2,\dots, p_{i-1})} \big) T_{i-1}, 
\end{align*}
where $h_{(p_1,p_2,\dots, p_{i-1})} \in \sH_{[1,i-1]}$. 
Moreover, 
by Lemma \ref{Li Lj}, 
we have 
\begin{align*}
T_{i-1} L_{i-1}^k T_{i-1} 
&= \big( L_i^k T_{i-1} - (q-q^{-1}) \sum_{s=0}^{k-1} L_{i-1}^s L_i^{k-s} \big) T_{i-1} 
\\
&= L_i^k - (q-q^{-1}) \sum_{s=1}^{k-1} L_{i-1}^s L_i^{k-s} T_{i-1}, 
\end{align*}
where we note that $T_{i-1}^2 - (q-q^{-1}) T_{i-1} =1$. 
We also have 
\begin{align*}
&T_{i-1} L_1^{p_1} L_2^{p_2} \dots L_{i-1}^{p_{i-1}} 
\\
&= L_1^{p_1} L_2^{p_2} \dots L_{i-2}^{p_{i-2}} 
	\big( L_i^{p_{i-1}} T_{i-1} - (q-q^{-1}) \sum_{s=0}^{p_{i-1} -1} L_{i-1}^s L_i^{p_{i-1} -s} \big) 
\\
&= L_1^{p_1} L_2^{p_2} \dots L_{i-2}^{p_{i-2}} L_i^{p_{i-1}} T_{i-1}
	- (q-q^{-1}) \sum_{s=0}^{p_{i-1} -1} L_1^{p_1} L_2^{p_2} \dots L_{i-2}^{p_{i-2}}  L_{i-1}^s L_i^{p_{i-1} -s}. 
\end{align*}
These equations imply Lemma \ref{Li lan k ran}. 
\end{proof}


\para 
{\it (Proof of Lemma  \ref{Lemma TT Lk}).} 
\begin{proof}
(\roi). 
In the case where $l <i$, this is clear from the defining relations of $\sH_{n,r}$. 

In the case where $l > i+1$, we have 
\begin{align*}
T_i L_l^{\lan k \ran} = T_i \wt{\TT}_{l-1,i} L_i^{\lan k \ran} \TT_{i, l-1} =T_i \wt{\TT}_{l-1, i+2} T_{i+1} T_i L_i^{\lan k \ran} \TT_{i,l-1}.
\end{align*} 
By the defining relation of $\sH_{n,r}$, 
we see that $T_i  \wt{\TT}_{l-1, i+2} T_{i+1} T_i = \wt{\TT}_{l-1, i+2} T_{i+1} T_i T_{i+1}$. 
On the other hand, we have already proven $T_{i+1} L_i^{\lan k \ran} = L_i^{\lan k \ran} T_{i+1}$. 
Then we see that 
\begin{align*}
T_i L_l^{\lan k \ran} = \wt{\TT}_{l-1, i+2} T_{i+1} T_i L_i^{\lan k \ran} T_{i+1} \TT_{i, l-1} 
= \wt{\TT}_{l-1,i} L_i^{\lan k \ran} \TT_{i, l-1} T_i = L_l^{\lan k \ran} T_i,
\end{align*}
where we note that 
$ T_{i+1} \TT_{i,l-1} = T_{i+1} T_i T_{i+1} \TT_{i+2, l-1} = T_i T_{i+1} T_i \TT_{i+2,l-1} = \TT_{i,l-1} T_i$. 

(\roii). 
$ L_{i+1}^{\lan k \ran} T_i = T_i L_i^{\lan k \ran} T_i T_i = T_i L_i^{\lan k \ran} \big( 1+ (q-q^{-1}) T_i \big) 
	= T_i L_i^{\lan k \ran} + (q-q^{-1}) L_{i+1}^{\lan k \ran}  $. 
	
(\roiii) is similar to (\roii). 

(\roiv). 
In the case where $l > j+1$ or $i > l$, it is clear from (\roi). 

We prove the case where $l=j+1$ by the induction on $j-i$. 
If $j-i=0$, then we have the formula by (\roiii). 
Suppose that $j- i >0$, 
and we have 
\begin{align*}
\TT_{i,j} L_{j+1}^{\lan k \ran} 
	= \TT_{i,j-1} T_j L_{j+1}^{\lan k \ran} 
	= \TT_{i, j-1} \big( L_j^{\lan k \ran} T_j + (q-q^{-1}) L_{j+1}^{\lan k \ran} \big) 
\end{align*} 
by (\roiii). 
Applying the inductive hypothesis and (\roi) to the right-hand side, 
we have 
\begin{align*}
\TT_{i,j} L_{j+1}^{\lan k \ran} 
&= \big( L_i^{\lan k \ran} \TT_{i,j-1} + (q-q^{-1}) \sum_{p=1}^{j-1-i+1} L_{i+p}^{\lan k \ran} \TT_{i, i + p -2} \TT_{i+p,j-1} \big) T_j 
	+ (q-q^{-1}) L_{j+1}^{\lan k \ran} \TT_{i,j-1}
\\
&= \big( L_i^{\lan k \ran} \TT_{i,j} + (q-q^{-1}) \sum_{p=1}^{j-i+1} L_{i+p}^{\lan k \ran} \TT_{i, i + p -2} \TT_{i+p,j} \big). 
\end{align*}

We prove the case where $j+1 > l \geq i$. 
By using (\roi), we have 
\begin{align*}
\TT_{i,j} L_l^{\lan k \ran} 
&= \TT_{i,l} \TT_{l+1,j} L_l^{\lan k \ran} 
=\TT_{i,l} L_l^{\lan k \ran} \TT_{l+1,j} 
= \TT_{i, l-1} T_l L_l^{\lan k \ran} T_l T_l^{-1} \TT_{l+1,j}
\\
&= \TT_{i,l-1} L_{l+1}^{\lan k \ran} \big( T_l - (q-q^{-1}) \big) \TT_{l+1,j} 
\\
&= L_{l+1}^{\lan k \ran} \TT_{i, l-1} T_l \TT_{l+1,j} - (q-q^{-1}) L_{l+1}^{\lan k \ran} \TT_{i, l-1} \TT_{l+1,j} 
\\
&= L_{l+1}^{\lan k \ran} \big( \TT_{i,j} - (q-q^{-1}) \TT_{i,l-1} \TT_{l+1,j} \big).  
\end{align*}

(\rov). 
In the case where $l > j+1$ or $i > l$, it is clear from (\roi). 

We prove the case where $j +1 \geq l >i$. 
By using (\roi) and (\roiii), we have 
\begin{align*}
\wt{\TT}_{j,i} L_l^{\lan k \ran} 
&= \wt{\TT}_{j,l-1} \wt{\TT}_{l-2,i} L_l^{\lan k \ran} 
= \wt{\TT}_{j,l-1} L_l^{\lan k \ran} \wt{\TT}_{l-2,i} 
= \wt{\TT}_{j,l} T_{l-1} L_l^{\lan k \ran} \wt{\TT}_{l-2,i} 
\\
&= \wt{\TT}_{j,l}  \big( L_{l-1}^{\lan k \ran} T_{l-1} + (q-q^{-1}) L_{l}^{\lan k \ran} \big) \wt{\TT}_{l-2,i}  
\\
&= L_{l-1}^{\lan k \ran} \wt{\TT}_{j,l} T_{l-1}\wt{\TT}_{l-2,i}   
	+ (q-q^{-1}) \wt{\TT}_{j,l} L_l^{\lan k \ran} \TT_{l,j} (\TT_{l,j})^{-1} \wt{\TT}_{l-2,i} 
\\
&= L_{l-1}^{\lan k \ran} \wt{\TT}_{j,i} + (q-q^{-1}) L_{j+1}^{\lan k \ran} (\TT_{l,j})^{-1} \wt{\TT}_{l-2,i}. 
\end{align*}

In the case where $l=i$, we have 
$\wt{\TT}_{j,i} L_i^{\lan k \ran} = \wt{\TT}_{j,i} L_i^{\lan k \ran} \TT_{i,j} (\TT_{i,j})^{-1} = L_{j+1}^{\lan k \ran} (\TT_{i,j})^{-1}$. 
\end{proof}

\para 
{\it (Proof of Lemma  \ref{Lemma Li Ljlan k ran}).} 
\begin{proof}
(\roi)-(a). 
If $i >j$, we have $L_i L_j^{\lan k \ran} = L_j^{\lan k \ran} L_i$ by Lemma \ref{Li Lj} (\roii). 

We prove (\roi)-(b) by the induction on $j$. 
In the case where $j=1$, it is clear from definitions. 
Suppose that $j >1$. 
By Lemma \ref{Li Lj} (\roii), we have 
\begin{align*}
L_j L_j^{\lan k \ran} 
= L_j T_{j-1} L_{j-1}^{\lan k \ran} T_{j-1} 
= \big( T_{j-1} L_{j-1} + (q-q^{-1}) L_j \big) L_{j-1}^{\lan k \ran} T_{j-1}. 
\end{align*}
Applying the induction hypothesis, we have  
\begin{align*}
L_j L_j^{\lan k \ran} 
&= T_{j-1}  \big\{ L_{j-1}^{\lan k \ran} L_{j-1} 
	\\ & \qquad 
	- (q-q^{-1}) \sum_{p=1}^{j-2} (L_{j-1}^{\lan k \ran} L_{j-p-1} - L_{j-p-1}^{\lan k \ran} L_{j-1} ) 
		( \TT_{j-p, j-2})^{-1} \TT_{j-p-1, j-2} \big\} T_{j-1} 
	\\ & \quad 
	+ (q-q^{-1}) L_{j-1}^{\lan k \ran} L_j T_{j-1}, 
\end{align*}
where we note that $L_j L_{j-1}^{\lan k \ran} = L_{j-1}^{\lan k \ran} L_j$ by (\roi)-(a). 
Note that $T_{j-1}^2 - (q-q^{-1}) T_{j-1} =1$ by the defining relations of $\sH_{n,r}$, 
and we have 
\begin{align*}
L_j L_j^{\lan k \ran} 
&= T_{j-1} L_{j-1}^{\lan k \ran} \big( T_{j-1}^2 - (q-q^{-1}) T_{j-1} \big) L_{j-1} T_{j-1} 
	\\ & \quad 
	- (q-q^{-1})  \sum_{p=1}^{j-2} \big\{ 
		T_{j-1} L_{j-1}^{\lan k \ran} (T_{j-1} T_{j-1}^{-1}) L_{j- p-1} 
		- L_{j- p-1}^{\lan k \ran} T_{j-1} L_{j-1} (T_{j-1} T_{j-1}^{-1}) \big\} 
	\\ & \qquad \times ( \TT_{j-p, j-2})^{-1} \TT_{j-p-1, j-2} T_{j-1}  
	+ (q-q^{-1}) L_{j-1}^{\lan k \ran} L_j T_{j-1}
\\
&= L_j^{\lan k \ran} L_j - (q-q^{-1}) L_j^{\lan k \ran} L_{j-1} T_{j-1} 
	\\ & \quad 
	- (q-q^{-1}) \sum_{p=1}^{j-2} 
		\big\{ L_j^{\lan k \ran} L_{j-p-1}  - L_{j-p-1}^{\lan k \ran} L_j \big\} T_{j-1}^{-1} ( \TT_{j-p, j-2})^{-1} \TT_{j-p-1, j-2} T_{j-1}  
	\\ & \qquad 
	+ (q-q^{-1}) L_{j-1}^{\lan k \ran} L_j T_{j-1} 
\\
&= L_j^{\lan k \ran} L_j 
	- (q-q^{-1}) \sum_{p=1}^{j-1} \big( L_j^{\lan k \ran} L_{j-p} - L_{j-p}^{\lan k \ran} L_j \big) (\TT_{j-p+1, j-1})^{-1} \TT_{j-p, j-1}. 
\end{align*}

We prove (\roi)-(c) by the induction on $j-i$. 
In the case where $j-i=1$, by using  Lemma \ref{Li Lj} (\roii), we have 
$L_i L_{i+1}^{\lan k \ran} 
= L_i T_i L_i^{\lan k \ran} T_i 
= \big( T_i L_{i+1} - (q-q^{-1}) L_{i+1} \big) L_i^{\lan k \ran} T_i$. 
Applying (\roi) to the right-hand side, we have 
$L_i L_{i+1}^{\lan k \ran} = T_i L_i^{\lan k \ran} L_{i+1} T_i  - (q-q^{-1}) L_i^{\lan k \ran} L_{i+1} T_i $, 
and, by using  Lemma \ref{Li Lj} (\roii) again, we have 
\begin{align*}
L_i L_{i+1}^{\lan k \ran} 
&= T_i L_i^{\lan k \ran} \big( T_i L_i + (q-q^{-1}) L_{i+1} \big) - (q-q^{-1}) L_i^{\lan k \ran} L_{i+1} T_i 
\\
&= L_{i+1}^{\lan k \ran} L_i + (q-q^{-1}) \big( L_{i+1}^{\lan k \ran} L_i - L_i^{\lan k \ran} L_{i+1} \big) T_i, 
\end{align*}
where we note that $T_i L_i^{\lan k \ran} L_{i+1} = T_i L_i^{\lan k \ran} (T_i T_i^{-1}) (T_i L_{i} T_i) 
	= L_{i+1}^{\lan k \ran} L_i T_i$. 
Suppose that $j-i>1$, 
we have 
$L_i L_j^{\lan k \ran} = L_i T_{j-1} L_{j-1}^{\lan k \ran} T_{j-1} = T_{j-1} L_i L_{j-1}^{\lan k \ran} T_{j-1}$ 
by Lemma \ref{Li Lj} (\roii). 
Applying the induction hypothesis, we have 
\begin{align*}
L_i L_j^{\lan k \ran} 
&= T_{j-1} \big( L_{j-1}^{\lan k \ran} L_i 
	+ (q-q^{-1}) ( L_{j-1}^{\lan k \ran} L_i - L_i^{\lan k \ran} L_{j-1} ) (\TT_{i+1, j-2})^{-1} \TT_{i, j-2} \big) T_{j-1} 
\\
&= T_{j-1} L_{j-1}^{\lan k \ran} L_i T_{j-1} 
	+ (q-q^{-1}) \big\{ T_{j-1} L_{j-1}^{\lan k \ran} ( T_{j-1} T_{j-1}^{-1} ) L_i 
		\\ & \qquad 
		- L_i^{\lan k \ran} T_{j-1} L_{j-1} (T_{j-1} T_{j-1}^{-1}) \big\} 
		(\TT_{i+1, j-2})^{-1} \TT_{i, j-2}  T_{j-1}, 
\end{align*}
where we note that $T_{j-1} L_i^{\lan k \ran} = L_i^{\lan k \ran} T_{j-1}$ by (\roi). 
Then, by Lemma \ref{Li Lj} (\roii), we have 
\begin{align*}
L_i L_j^{\lan k \ran} 
&= T_{j-1} L_{j-1}^{\lan k \ran}  T_{j-1} L_i 
	\\ & \quad 
	+ (q-q^{-1}) \big\{ T_{j-1} L_{j-1}^{\lan k \ran}  T_{j-1} L_i 
		- L_i^{\lan k \ran} T_{j-1} L_{j-1} T_{j-1}  \big\} 
		T_{j-1}^{-1} (\TT_{i+1, j-2})^{-1} \TT_{i, j-2} T_{j-1}  
\\
&= L_j^{\lan k \ran} L_i + (q-q^{-1}) \big( L_j^{\lan k \ran} L_i - L_i^{\lan k \ran} L_j \big) (\TT_{i+1, j-1})^{-1} \TT_{i, j-1}.
\end{align*}

(\roii)-(a) follows form (\roi)-(a). 

We prove (\roii)-(b) by the induction on $j$. 
In the case where $j=1$, 
it is clear that $L_1^{-1} L_1^{\lan k \ran} = L_1^{\lan k \ran} L_1^{-1}$ by definitions. 
Suppose that $j >1$, 
we have 
\begin{align*}
L_j^{-1} L_j^{\lan k \ran} 
= (T_{j-1} L_{j-1} T_{j-1})^{-1} ( T_{j-1} L_{j-1}^{\lan k \ran} T_{j-1} ) 
= T_{j-1}^{-1} L_{j-1}^{-1} L_{j-1}^{\lan k \ran} T_{j-1}.
\end{align*}
Applying the induction hypothesis to the right-hand side, we have  
\begin{align*}
L_j^{-1} L_j^{\lan k \ran} 
&= T_{j-1}^{-1} \big\{ 
	L_{j-1}^{\lan k \ran} L_{j-1}^{-1} \wt{\TT}_{j-2,1} \TT_{1,j-2} 
		\\ & \qquad 
		- (q-q^{-1}) \sum_{p=1}^{j-2} L_{j-p-1}^{\lan k \ran} L_{j-p-1}^{-1} (\TT_{j-p, j-2})^{-1} \wt{\TT}_{j-p-2,1} \TT_{1,j-2} 
	\big\} T_{j-1}. 
\end{align*}
By using Lemma \ref{Li Lj} (\roii) and Lemma \ref{Lemma TT Lk} (\roi), (\roiii) together with 
$T_{j-1}^{-1} L_j^{\lan k \ran} T_{j-1}^{-1} = L_{j-1}^{\lan k \ran}$, we have 
\begin{align*}
L_j^{-1} L_j^{\lan k \ran}  
&= \big( L_j^{\lan k \ran} T_{j-1}^{-1} - (q-q^{-1}) L_{j-1}^{\lan k \ran} \big) L_{j-1}^{-1} \wt{\TT}_{j-2,1} \TT_{1, j-2} T_{j-1} 
	\\ & \quad 
	- (q-q^{-1}) \sum_{p=1}^{j-2} L_{j-p-1}^{\lan k \ran} L_{j-p-1}^{-1} 
		T_{j-1}^{-1}  (\TT_{j-p, j-2})^{-1} \wt{\TT}_{j-p-2,1} \TT_{1,j-2} T_{j-1}
\\
&= L_j^{\lan k \ran} T_{j-1}^{-1} L_{j-1}^{-1} (T_{j-1}^{-1} T_{j-1}) \wt{\TT}_{j-2,1} \TT_{1,j-2} T_{j-1}  
	\\ & \quad 
	- (q-q^{-1}) L_{j-1}^{\lan k \ran} L_{j-1}^{-1} \wt{\TT}_{j-2,1} \TT_{1,j-2} T_{j-1} 
	\\ & \quad 
	- (q-q^{-1}) \sum_{p=1}^{j-2} L_{j-p-1}^{\lan k \ran} L_{j-p-1}^{-1} 
		T_{j-1}^{-1}  (\TT_{j-p, j-2})^{-1} \wt{\TT}_{j-p-2,1} \TT_{1,j-2} T_{j-1}
\\
&= L_j^{\lan k \ran} L_j^{-1} \wt{\TT}_{j-1,1} \TT_{1,j-1} 
	- (q-q^{-1}) \sum_{p=1}^{j-1} L_{j-p}^{\lan k \ran} L_{j-p}^{-1} (\TT_{j-p+1, j-1})^{-1} \wt{\TT}_{j-p-1,1} \TT_{1,j-1}. 
\end{align*}

We prove (\roii)-(c). 
In the case where $j=i+1$, we have 
\begin{align*}
L_i^{-1} L_{i+1}^{\lan k \ran} 
= L_i^{-1} T_i L_i^{\lan k \ran} T_i 
= \big( T_i L_{i+1}^{-1} + (q-q^{-1}) L_i^{-1} \big) L_i^{\lan k \ran} T_i 
\end{align*}
by Lemma \ref{Li Lj} (\roii). 
Applying (\roii)-(a) and (\roii)-(b) to the right-hand side, we have 
\begin{align*}
L_i^{-1} L_{i+1}^{\lan k \ran} 
&= T_i L_i^{\lan k \ran} L_{i+1}^{-1} T_i 
	\\ & \quad 
	+ (q-q^{-1}) \big\{ L_i^{\lan k \ran} L_i^{-1} \wt{\TT}_{i-1,1} \TT_{1,i-1} 
		\\ & \qquad 
		- (q-q^{-1}) \sum_{p=1}^{i-1} L_{i-p}^{\lan k \ran} L_{i-p}^{-1} (\TT_{i-p+1, i-1})^{-1} \wt{\TT}_{i-p-1,1} \TT_{1,i-1}
	\big\} T_i. 
\end{align*}
By Lemma \ref{Lemma TT Lk} (\roii), we have 
\begin{align*}
L_i^{-1} L_{i+1}^{\lan k \ran} 
	&= \big( L_{i+1}^{\lan k \ran} T_i - (q-q^{-1}) L_{i+1}^{\lan k \ran} \big) L_{i+1}^{-1} T_i 
	+ (q-q^{-1}) L_i^{\lan k \ran} L_i^{-1} \wt{\TT}_{i-1,1} \TT_{1,i-1} T_i 
	\\ & \quad 
	- (q-q^{-1})^2 \sum_{p=1}^{i-1} L_{i-p}^{\lan k \ran} L_{i-p}^{-1} (\TT_{i-p+1, i-1})^{-1} \wt{\TT}_{i-p-1,1} \TT_{1,i-1}T_i
\\
&= L_{i+1}^{\lan k \ran} L_i^{-1} - (q-q^{-1}) L_{i+1}^{\lan k \ran} L_{i+1}^{-1} T_i 
	+ (q-q^{-1}) L_i^{\lan k \ran} L_i^{-1} \wt{\TT}_{i-1,1} \TT_{1,i} 
	\\ & \quad 
	- (q-q^{-1})^2 \sum_{p=1}^{i-1} L_{i-p}^{\lan k \ran} L_{i-p}^{-1} \wt{\TT}_{i-p-1,1} (\TT_{i-p+1, i-1})^{-1}  \TT_{1,i}.
\end{align*}
Suppose that $j > i+1$, we have 
\begin{align*}
L_i^{-1} L_j^{\lan k \ran} = 
L_i^{-1} \wt{\TT}_{j-1, i+1} L_{i+1}^{\lan k \ran} \TT_{i+1,j-1} 
= \wt{\TT}_{j-1, i+1} L_i^{-1} L_{i+1}^{\lan k \ran} \TT_{i+1,j-1} 
\end{align*}
by Lemma \ref{Li Lj} (\roii). 
Applying the above calculation of the case where $j=i+1$, we have 
\begin{align*}
L_i^{-1} L_j^{\lan k \ran} 
&= \wt{\TT}_{j-1,i+1} \big\{ 
	L_{i+1}^{\lan k \ran} L_i^{-1} - (q-q^{-1}) L_{i+1}^{\lan k \ran} L_{i+1}^{-1} T_i 
	+ (q-q^{-1}) L_i^{\lan k \ran} L_i^{-1} \wt{\TT}_{i-1,1} \TT_{1,i} 
	\\ & \quad 
	- (q-q^{-1})^2 \sum_{p=1}^{i-1} L_{i-p}^{\lan k \ran} L_{i-p}^{-1} \wt{\TT}_{i-p-1,1} (\TT_{i-p+1, i-1})^{-1}  \TT_{1,i}
	\big\} \TT_{i+1, j-1}. 
\end{align*}
By using Lemma \ref{Li Lj} (\roii) and Lemma \ref{Lemma TT Lk} (\rov), we have  
\begin{align*}
&L_i^{-1} L_j^{\lan k \ran} 
\\
&= \wt{\TT}_{j-1,i+1} L_{i+1}^{\lan k \ran} \TT_{i+1,j-1} L_i^{-1}  
	\\ & \quad  
	- (q-q^{-1}) \wt{\TT}_{j-1,i+1} L_{i+1}^{\lan k \ran} ( \TT_{i+1,j-1} (\TT_{i+1,j-1})^{-1} ) L_{i+1}^{-1} 
		( \wt{\TT}_{j-1,i+1})^{-1} \wt{\TT}_{j-1,i+1} ) T_i \TT_{i+1,j-1} 
	\\ & \quad 
	+ (q-q^{-1}) L_i^{\lan k \ran} L_i^{-1} \wt{\TT}_{j-1, i+1} \wt{\TT}_{i-1,1} \TT_{1,i} \TT_{i+1,j-1}  
	\\ & \quad 
	- (q-q^{-1})^2 \sum_{p=1}^{i-1} L_{i-p}^{\lan k \ran} L_{i-p}^{-1} 
	\wt{\TT}_{j-1, i+1} \wt{\TT}_{i-p-1,1} (\TT_{i-p+1, i-1})^{-1}  \TT_{1,i} \TT_{i+1,j-1}  
\\
&=  L_{j}^{\lan k \ran}  L_i^{-1}  
	- (q-q^{-1})  L_{j}^{\lan k \ran}   L_{j}^{-1} \wt{\TT}_{j-1,i+1}  \TT_{i,j-1} 
	\\ & \quad 
	+ (q-q^{-1}) L_i^{\lan k \ran} L_i^{-1} \wt{\TT}_{j-1, i+1} \wt{\TT}_{i-1,1} \TT_{1,j-1} 
	\\ & \quad 
	- (q-q^{-1})^2 \sum_{p=1}^{i-1} L_{i-p}^{\lan k \ran} L_{i-p}^{-1} 
	\wt{\TT}_{j-1, i+1} \wt{\TT}_{i-p-1,1} (\TT_{i-p+1, i-1})^{-1}  \TT_{1,j-1}.
\end{align*}
\end{proof}


\para 
{\it (Proof of Lemma  \ref{Lemma Lj Qk LJlankran}).} 
\begin{proof} 
We  prove that 
\begin{align}
\label{Lj wtTTj-1i}
L_j \wt{\TT}_{j-1,i} = \wt{\TT}_{j-1,i} L_i + (q-q^{-1}) \sum_{p=i}^{j-1} \wt{\TT}_{j-1, p+1} \wt{\TT}_{p-1,i} L_{p+1}
\quad (j > i \geq 1)
\end{align} 
by the induction on $j-i$. 
If $j-i=1$, this is just the formula in Lemma \ref{Li Lj} (\roii). 
Suppose that $j-i>1$, 
and we have 
\begin{align*}
L_j \wt{\TT}_{j-1,i} = L_{j} T_{j-1} \wt{\TT}_{j-2,i} = ( T_{j-1} L_{j-1} + (q-q^{-1}) L_j ) \wt{\TT}_{j-2,i}
\end{align*}
by Lemma \ref{Li Lj} (\roii). 
Applying the induction hypothesis and Lemma \ref{Li Lj} (\roii), we have 
\begin{align*}
L_j \wt{\TT}_{j-1,i} 
&= T_{j-1} \big( \wt{\TT}_{j-2, i} L_i + (q-q^{-1}) \sum_{p=i}^{j-2} \wt{\TT}_{j-2, p+1} \wt{\TT}_{p-1,i} L_{p+1} \big) 
	+ (q-q^{-1}) \wt{\TT}_{j-2,i} L_j 
\\
&= \wt{\TT}_{j-1,i} L_i + (q-q^{-1}) \sum_{p=i}^{j-2} \wt{\TT}_{j-1,p+1} \wt{\TT}_{p-1,i} L_{p+1} + (q-q^{-1}) \wt{\TT}_{j-2,i} L_{j}. 
\end{align*}
Then, the formula \eqref{Lj wtTTj-1i} implies that 
\begin{align*}
&(L_j - Q_k) L_j^{\lan k \ran} 
\\
&= (L_j - Q_k) \wt{\TT}_{j-1,1} L_1^{\lan k \ran} \TT_{1,j-1} 
\\
&= \big\{ \wt{\TT}_{j-1,1} L_1 + (q-q^{-1}) \sum_{p=1}^{j-1} \wt{\TT}_{j-1,p+1}\wt{\TT}_{p-1,1} L_{p+1} \big\} 
	L_1^{\lan k \ran} \TT_{1,j-1} 
	- Q_k \wt{\TT}_{j-1,1} L_1^{\lan k \ran} \TT_{1,j-1}
\\
&= \wt{\TT}_{j-1,1} (L_1 - Q_k) L_1^{\lan k \ran} \TT_{1,j-1} 
	+ (q-q^{-1}) \sum_{p=1}^{j-1} \wt{\TT}_{j-1,p+1}\wt{\TT}_{p-1,1} L_{p+1} L_1^{\lan k \ran} \TT_{1,j-1}.  
\end{align*}
By using  Lemma \ref{Li Lj} (\roii) and Lemma  \ref{Lemma Li Ljlan k ran}, 
we have  
\begin{align*}
&(L_j - Q_k) L_j^{\lan k \ran} 
\\
&= \wt{\TT}_{j-1,1} L_1^{\lan k+1 \ran} \TT_{1,j-1} 
	+ (q-q^{-1}) \sum_{p=1}^{j-1} \wt{\TT}_{j-1,p+1} \wt{\TT}_{p-1,1} L_1^{\lan k \ran} \TT_{1,p-1} L_{p+1} \TT_{p, j-1} 
\\
&= L_j^{\lan k +1 \ran} + (q-q^{-1}) \sum_{p=1}^{j-1} \wt{\TT}_{j-1, p+1} L_p^{\lan k \ran} L_{p+1}  \TT_{p,j-1}. 
\end{align*}
Apply Lemma \ref{Lemma TT Lk} (\rov) to the right-hand side, 
and we obtain the formula in Lemma \ref{Lemma Lj Qk LJlankran}. 
\end{proof}




\end{document}